\documentclass[11pt]{amsart}
\usepackage{amsmath}
\usepackage{amsthm,amscd,graphicx}
\usepackage{dsfont}
\usepackage{amssymb}
\usepackage{amsfonts}
\usepackage{mathrsfs}
\usepackage{mdframed}
\usepackage{colortbl}
\usepackage{bm}
\usepackage{tikz}
\usepackage{tikz-cd}
\usetikzlibrary{calc,matrix}
\usepackage{enumitem}
\usepackage{float}
\usepackage[all]{xy}
\usetikzlibrary{decorations.pathmorphing,decorations.markings,decorations.pathreplacing,arrows,shapes}
\usetikzlibrary{knots}

\let\oldtocsection=\tocsection
\let\oldtocsubsection=\tocsubsection
\let\oldtocsubsubsection=\tocsubsubsection
\renewcommand{\tocsection}[2]{\hspace{0em}\oldtocsection{#1}{#2}}
\renewcommand{\tocsubsection}[2]{\hspace{2em}\oldtocsubsection{#1}{#2}}
\renewcommand{\tocsubsubsection}[2]{\hspace{3em}\oldtocsubsubsection{#1}{#2}}

\usepackage[bookmarks=true, bookmarksopen=true,%
bookmarksdepth=3,bookmarksopenlevel=2,%
colorlinks=true,%
linkcolor=black,%
citecolor=blue,%
filecolor=blue,%
menucolor=blue,%
urlcolor=blue]{hyperref}

\newtheorem{thm}{Theorem}[section]
\newtheorem{prop}[thm]{Proposition}
\newtheorem{cor}[thm]{Corollary}
\newtheorem{lem}[thm]{Lemma}

\theoremstyle{definition}
\newtheorem{defn}[thm]{Definition}

\newtheorem{rmk}[thm]{Remark}
\newtheorem{exmp}[thm]{Example}

\newtheorem{assumption}[thm]{Assumption}

\makeatletter
\let\c@equation\c@thm
\makeatother
\numberwithin{equation}{section}

\newcommand{\rank}{\mathrm{rank}}
\newcommand{\conf}{\mathrm{Conf}}

\newcommand{\Gr}{\mathrm{Gr}}
\newcommand{\GL}{\mathrm{GL}}

\newcommand{\SL}{\mathrm{SL}}

\newcommand{\DT}{\mathrm{DT}}

\renewcommand{\sc}{\mathrm{sc}}

\newcommand{\doverline}[1]{\overline{\overline{#1}}}

\newcommand{\uf}{\mathrm{uf}}

\newcommand{\spec}{\mathrm{Spec} \ }
\newcommand{\Aug}{\mathrm{Aug}}

\newcommand{\id}{\mathrm{id}}

\newcommand{\Triv}{\mathrm{Triv}}

\newcommand{\tb}{\mathrm{tb}}

\newcommand{\A}{\mathsf{A}}
\newcommand{\B}{\mathsf{B}}

\newcommand{\G}{\mathsf{G}}
\newcommand{\K}{\mathsf{K}}

\newcommand{\Q}{\mathsf{Q}}
\newcommand{\T}{\mathsf{T}}
\newcommand{\U}{\mathsf{U}}
\newcommand{\R}{\mathsf{R}}
\newcommand{\W}{\mathsf{W}}

\newcommand{\Br}{\mathsf{Br}}

\newcommand{\cA}{\mathcal{A}}

\newcommand{\cC}{\mathcal{C}}

\newcommand{\cM}{\mathcal{M}}
\newcommand{\cR}{\mathcal{R}}
\newcommand{\cT}{\mathcal{T}}
\newcommand{\cX}{\mathcal{X}}
\newcommand{\cP}{\mathcal{P}}
\newcommand{\cI}{\mathcal{I}}

\newcommand{\bbR}{\mathbb{R}}
\newcommand{\bbC}{\mathbb{C}}
\newcommand{\bbZ}{\mathbb{Z}}

\renewcommand{\vec}[1]{\mathbf{#1}}
\newcommand{\inprod}[2]{\left\langle#1,#2\right\rangle}

\newcommand{\del}{\partial}
\newcommand{\ep}{\epsilon}

\setenumerate[1]{label={(\arabic*)}}
\setlist[enumerate]{itemsep=0mm}

\allowdisplaybreaks

\title{Augmentations, Fillings, and Clusters}

\author{Honghao Gao}
\address{Yau Mathematical Sciences Center, Jingzhai, Tsinghua University, Haidian District, Beijing, 100084, China}
\email{gaohonghao@tsinghua.edu.cn}

\author{Linhui Shen}
\address{Department of Mathematics, Michigan State University, 619 Red Cedar Road, East Lansing, MI 48824, USA}
\email{linhui@math.msu.edu}

\author{Daping Weng}
\address{Department of Mathematics, University of California, Davis, One Shields Avenue, Davis, CA 95616, USA}
\email{dweng@ucdavis.edu}

\date{}

\begin{document}

\maketitle

\begin{abstract}

We investigate positive braid Legendrian links via a Floer-theoretic approach and prove that their augmentation varieties are cluster $\mathrm{K}_2$ (aka. $\mathcal{A}$-) varieties. Using the exact Lagrangian cobordisms of Legendrian links in \cite{EHK12}, we prove that a large family of exact Lagrangian fillings of positive braid Legendrian links correspond to cluster seeds of their augmentation varieties. We solve the infinite-filling problem for positive braid Legendrian links; i.e., whenever a positive braid Legendrian link is not of type ADE, it admits infinitely many exact Lagrangian fillings up to Hamiltonian isotopy.

\end{abstract}

\tableofcontents

\section{Introduction}

This paper is the first attempt to relate Floer theory and cluster algebras in the context of contact manifolds and Legendrian knots.
Starting from \cite{STWZ}, and together with subsequent \cite{SWflag, CZ}, the relations between microlocal sheaf theory and cluster Poisson (aka $\cX$)-varieties have been studied for several classes of Legendrian links.
It is natural to ask whether the cluster structure exists on the celebrated Floer-theoretic invariant, namely the Chekanov-Eliashberg dga and its augmentations \cite{Che02}. Despite the famous augmentation-sheaf correspondence for Legendrian links \cite{NRSSZ}, which suggests a similar cluster structure on the augmentation moduli, it is to our surprise that we obtain an intrinsic cluster $\mathrm{K}_2$ (aka $\cA$)-structure on the augmentation variety of any positive braid Legendrian link. This paper further utilizes this new cluster ${\rm K}_2$ structure to build invariants for exact Lagrangian fillings. As an application, we prove that positive braids that do not underline finite type quivers admit infinitely many Lagrangian fillings. To our knowledge, this is by far the largest family of Legendrian links satisfying the infinite filling properties.

\subsection{Context}
In the standard contact three-space  $(\bbR^3,\xi_{st} = \ker \alpha)$ with $\alpha={\rm d}z-y {\rm d}x$, a Legendrian link $\Lambda$ is a smooth one-dimensional submanifold where $\alpha|_\Lambda =0$. The Chekanov-Eliashberg differential graded algebra (CE dga) is the first non-classical algebraic invariant for Legendrian links \cite{Che02}. An exact Lagrangian cobordism between Legendrian links functorially induces an algebraic map between the dgas \cite{EHK12}. Following this functoriality, each exact Lagrangian filling $L$ gives rise to an embedding of the decorated $\GL_1$-character variety\footnote{Here the decoration means a specific trivialization of the line bundle near the boundary of the surface. See Definition \ref{decorated.character} for a precise description.} on $L$ into the augmentation variety. The images of these morphisms are invariants that distinguish exact Lagrangian fillings.

Cluster algebras are a class of commutative algebras introduced by Fomin and Zelevinsky \cite{FZI}.
Since its inception, the theory of cluster algebras has found tremendous applications in diverse areas of mathematics and physics. Fock and Goncharov \cite{FGensemble} introduce a pair $(\mathscr{X}, \mathscr{A})$ of log Calabi-Yau varieties, which are a geometric enrichment of the cluster algebras. The variety $\mathscr{X}$ carries a natural Poisson structure and is referred to as a {\it cluster Poisson variety}. The variety $\mathscr{A}$ carries a canonical class in the Milnor $\mathrm{K}_2$ group of its function field and is referred to as a {\it cluster ${\rm K}_2$ variety}. See Section 6.2 of {\it loc.cit.} for the construction of such a canonical class. The duality between $\mathscr{A}$ and $\mathscr{X}$, conjectured by Fock and Goncharov, has been realized by Gross, Hacking, Keel, and Kontsevich \cite{GHKK} under the framework of scattering diagrams and mirror symmetry. Despite such a duality, the geometries of $\mathscr{X}$ and $\mathscr{A}$ are rather different. For the convenience of the reader, we recall the definition of cluster varieties in Appendix \ref{app.A}.

\smallskip 

This paper focuses on certain representatives of positive braid Legendrian links with maximum Thurston-Bennequin (tb) numbers. It follows from \cite[Theorem 3.4]{EV} that a positive braid has a unique Legendrian representative with maximal tb. We include a construction of these Legendrian representatives in Section \ref{sec2.1.positive-braid}.
We prove that their augmentation varieties carry natural cluster $\mathrm{K}_2$ structures. We consider a large family of exact Lagrangian fillings and prove that each filling induces a cluster seed of the augmentation variety. As an application, we prove that all positive braid Legendrian links, except those underlying ADE Dynkin-type quivers, admit infinitely many non-Hamiltonian isotopic exact Lagrangian fillings.

The classification of exact Lagrangian fillings is a central but rather difficult problem. Except for the unique filling for unknot \cite{EP}, most subsequent works focus on giving a lower bound on the number of distinct fillings. The existence of infinitely many exact Lagrangian fillings was not  known until the year 2020. Within the year, several methods emerged concurrently and each successfully solved this problem for a certain class. Two proceeding results are:
\begin{itemize}
    \item Casals-Gao \cite{CasalsGao} proved that any positive torus $(n,m)$-link, $(n,m)\neq (2,m)$, $(3,3)$, $(3,4)$, and $(3,5)$, admits infinitely many fillings. The proof uses Legendrian loops, microlocal sheaves, and cluster structures on Grassmannians.
    \item Casals-Zaslow \cite{CZ} proved that the rainbow closure of a class of 3-strand positive braids admit infinitely many fillings. The proof uses Legendrian weaves and cluster Poisson structures on moduli space of microlocal sheaves. 
\end{itemize}
The present paper investigates the infinite-filling problem for all positive braid closures, covering all examples of \cite{CasalsGao, CZ} as special cases. 

This paper is based on a Floer theoretical approach. In particular, our proof uses the Ekholm-Honda-K\'{a}lm\'{a}n (EHK) functor \cite{EHK12} instead of the microlocal sheaves in \cite{CasalsGao, CZ}. In this paper, we show for the first time that the augmentation varieties are cluster $\mathrm{K}_2$ varieties. It is an interesting direction for future research to compare with cluster structures arising from sheaves.

We would like to remark that, shortly after our result, Casals-Ng \cite{CasalsNg} proved the existence of infinitely many fillings for certain Legendrian links that are not positive braid closures. They use holomorphic curves but without cluster theory. 

\subsection{\texorpdfstring{Cluster ${\rm K}_2$ Structures on Augmentation Varieties}{}}

For any positive braid word $\beta$, we construct a quiver $Q_\beta$ via the following three-step procedure:
\begin{enumerate}
    \item[Step 1] Plot $\beta$ on $\mathbb{R}^2$ horizontally. Put a vertex in each region of the diagram sandwiched by strands (including the leftmost and the rightmost ``half-open" regions).
    \item[Step 2] At each crossing, draw the following arrow pattern among neighboring vertices (dashed arrows are of weight $1/2$):
    \[
    \begin{tikzpicture}[baseline=-12]
    \draw [lightgray] (-1,-0.5) to [out=0,in=180] (1,0.5);
    \draw [lightgray] (-1,0.5) to [out=0,in=180] (1,-0.5);
    \node (w) at (-0.75,0) [] {$\bullet$};
    \node (e) at (0.75,0) [] {$\bullet$};
    \node (s) at (0,-0.5) [] {$\bullet$};
    \node (n) at (0,0.5) [] {$\bullet$};
    \draw [->] (w) -- (e);
    \draw [->,dashed] (n) -- (w);
    \draw [->,dashed] (s) -- (w);
    \draw [->,dashed] (e) -- (n);
    \draw [->,dashed] (e) -- (s);
    \end{tikzpicture}.
    \]
    \item[Step 3] Sum up the arrows between each pair of vertices. Freeze the vertices in the rightmost regions. Delete the leftmost vertices and their incident arrows.
\end{enumerate}

\begin{exmp} Below is a positive braid word $\beta=s_1^6s_2s_1^3s_2$ and its quiver $Q_\beta$. 
\begin{figure}[h]
    \begin{tikzpicture}
       \draw[lightgray, dashed] (0, 1.2)--(0,-.2);
       \draw[lightgray, dashed] (11.5, 1.2)--(11.5,-.2);
    \draw [lightgray] (0,1) to [out=0,in=180] (1,0.5) to [out=0,in=180] (2,1) to [out=0,in=180] (3,0.5) to [out=0,in=180] (4,1) to [out=0,in=180] (5,0.5) to [out=0,in=180] (6,1) -- (7,1) to [out=0,in=180] (8,0.5) to [out=0,in=180] (9,1) to [out=0,in=180] (10,0.5) to [out=0,in=180] (11,0)--(11.5,0);
    \draw [lightgray] (0,0.5) to [out=0,in=180] (1,1) to [out=0,in=180] (2,0.5)  to [out=0,in=180] (3,1) to [out=0,in=180] (4,0.5) to [out=0,in=180] (5,1) to [out=0,in=180] (6,0.5) to [out=0,in=180] (7,0) -- (10,0) to [out=0,in=180] (11,0.5)--(11.5,0.5);
    \draw [lightgray] (0,0) -- (6,0)  to [out=0,in=180] (7,0.5) to [out=0,in=180] (8,1) to [out=0,in=180] (9,0.5) to [out=0,in=180] (10,1) -- (11,1)--(11.5,1);
    \foreach \i in {1,...,5,8,9}
    {
    \node (\i) at (\i,0.75) [] {$\bullet$};
    }
    \node (6) at (6.5,0.75) [] {$\bullet$};
    \node (7) at (8.5,0.25) [] {$\bullet$};
    \node (10) at (10.2,0.75) [] {$\square$};
    \node (11) at (11.2,0.25) [] {$\square$};
    \draw [->] (1) -- (2);
    \draw [->] (2) -- (3);
    \draw [->] (3) -- (4);
    \draw [->] (4) -- (5);
    \draw [->] (5) -- (6);
    \draw [->] (7) -- (6);
    \draw [->] (6) -- (8);
    \draw [->] (8) -- (9);
    \draw [->] (9) -- (10);
    \draw [->] (7) -- (11);
    \draw [->] (10) -- (7);
    \draw [->,dashed] (11) -- (10);
    \end{tikzpicture}
    \caption{An $\mathrm{E}_9$ quiver with two frozen vertices}
    \label{E_9}
\end{figure}
\end{exmp}

\vskip 1mm

Let $\Lambda_\beta$ be the positive braid Legendrian link associated to  $\beta$ as in Section \ref{sec2.1.positive-braid}. Let   $\Aug\left(\Lambda_\beta\right)$ be the augmentation variety of $\Lambda_\beta$ defined over an algebraically closed field $\mathbb{F}$ of characteristic 2 as in Definition \ref{def aug var}. Our first main result is as follows.

\begin{thm}[Theorem \ref{3.6}, Corollary \ref{all Reeb chords are cluster variables}, and Proposition \ref{prop initial seed}]\label{Mainthm1} The augmentation variety $\Aug(\Lambda_\beta)$ is a cluster $\mathrm{K}_2$ variety associated to the quiver $Q_\beta$. The degree zero Reeb chords of  $\Lambda_\beta$ are  cluster variables that generate the coordinate ring of $\Aug(\Lambda_\beta)$.
\end{thm}

\begin{rmk}
As defined in Definition \ref{clusterK2}, a cluster $\mathrm{K}_2$ variety is the spectrum of an {\it upper cluster algebra.} Up to codimension 2, each cluster $\mathrm{K}_2$ variety is obtained by gluing a collection of algebraic tori. The transition maps between different algebraic tori are given by particular relations called cluster mutations.
\end{rmk}

\begin{rmk}
Associated with each quiver is a {\it cluster algebra} $\mathcal{A}$ generated by cluster variables, and an {\it upper cluster algebra} $\mathcal{U}$ that is the intersection of the ring of Laurent polynomials for each seed \cite{BFZ}. The Laurent phenomenon of cluster variables implies that $\mathcal{A}\subset\mathcal{U}$, but in general $\mathcal{A}\neq \mathcal{U}$.
The problem when $\mathcal{A}=\mathcal{U}$ is a fundamental question in cluster theory. See \cite{GLS} for its application on the quantization of cluster algebras. 
As an application of  Theorem \ref{Mainthm1}, the upper cluster algebra $\mathcal{U}$ associated to the quiver $Q_\beta$ coincides the coordinate ring $\mathcal{O}({\rm Aug}(\Lambda_\beta))$. Meanwhile, the $\mathcal{O}({\rm Aug}(\Lambda_\beta))$ is generated by the Reeb chords as cluster variables. Therefore we get $\mathcal{A}=\mathcal{U}$ for the quiver $Q_\beta$. 
\end{rmk}

\subsection{From Fillings to Cluster Seeds}
Our second main result establishes a natural correspondence between a large family of exact Lagrangian fillings, which we call ``admissible fillings", of positive braid closures and the  seeds of the cluster $\mathrm{K}_2$ structure on their augmentation varieties. 

\begin{defn} An exact Lagrangian cobordism $L$ between positive braid Legendrian links is \emph{admissible} if it is a concatenation of the following exact Lagrangian cobordisms:
\begin{enumerate}
    \item saddle cobordism, which resolves a crossing inside the positive braid;
    \item braid move, also known as a Legendrian Reidemeister III move;
    \item cyclic rotation, which changes $\Lambda_{\delta s_i}\leftrightarrow \Lambda_{s_i\delta}$ for any positive braid $\delta$ and any elementary braid $s_i$;
    \item minimum cobordism, which is the unique filling of a maximal tb unknot.
\end{enumerate}
An \emph{admissible filling} of $\Lambda_\beta$ is an admissible cobordism from the empty set to $\Lambda_\beta$.
\end{defn}

We refer readers to Section 3.1 for the definition of exact Lagrangian cobordisms and fillings. A notable property of exact Lagrangian cobordisms is that they are directed. While a smooth cobordism surface can be reversed to interchange the two end, the same operation does not apply to the exact Lagrangian setting due to the directionality of the Liouville vector field, which is $\del_t$ in the symplectization $(\bbR_t\times \bbR^3_{xyz}, d(e^t\alpha))$. For instance, a Lagrangian cobordism $L: \Lambda_- \rightarrow \Lambda_+$ must satisfy $tb(\Lambda_+) - tb(\Lambda_-) = -\chi(L)$, where $tb(\Lambda)$ is the Thurston-Bennequin number of the Legendrian $\Lambda$ and $\chi(L)$ is the Euler characteristic of the surface $L$ \cite{chantraine2010}. Even in the case of Lagrangian concordance, which means the cobordism surface is smoothly a union of cylinders, the relation is still not symmetric \cite{chantraine2015}.

Let $L$ be an exact Lagrangian filling of $\Lambda$ with a collection $\mathcal{T}$ of marked points on $\Lambda$. Following Fock and Goncharov \cite{FGteich}, we consider the moduli space $\mathscr{A}(L, \mathcal{T})$ of decorated ${\rm GL}_1$-local systems (See Definition \ref{decorated.character}). Applying the Ekholm-Honda-K\'alm\'an functor, we obtain an open embedding 
\[
\alpha_L: \mathscr{A}(L, \mathcal{T}) \longrightarrow \Aug(\Lambda_\beta).
\]
where $\mathscr{A}(L, \mathcal{T})$ is isomorphic to an algebraic torus.

\begin{thm}[Theorem \ref{thm correspondence} and Corollary \ref{cor 3.40}]\label{Mainthm2} For any admissible filling $L$ of $\Lambda_\beta$, the image of  $\alpha_L$ is an open torus which determines a cluster seed of $\Aug(\Lambda_\beta)$. Admissible fillings determining distinct cluster seeds are non-Hamiltonian isotopic. 
\end{thm}

Theorem \ref{Mainthm2} gives a new method to distinguish non-Hamiltonian isotopic admissible fillings of positive braid closures via computing their corresponding cluster seeds. This theorem generalizes the methods in \cite{EHK12, Pan_2017}, which use the set of augmentations induced from a filling as an invariant and discovered a Catalan number worth of fillings for torus $(2,n)$ links. A different approach using sheaves can be found in \cite{STWZ}. The set of augmentations is the chart in the corresponding cluster seed.

The machinery of cluster theory allows us to develop an efficient method to compute the induced toric chart of augmentations from an admissible filling. Instead of keeping track of holomorphic curves bounded by the filling (following the recipe of \cite{EHK12}), one can compute the sequence of cluster mutations, find the induced cluster variables, and their non-vanishing loci give the desired result. We summarize the algorithm in Section \ref{sec algorithm} and build a program to implement the computation. Every cluster seed admits a complete combinatorial invariant called the $g$-matrix. Our algorithm presents an efficient method to explicitly compute the $g$-matrix as an invariant for admissible fillings.

\begin{rmk} 
The CE dga is defined over $\bbZ_2$ with contributions from marked points. Theorems \ref{Mainthm1} and \ref{Mainthm2} can be enhanced to characteristic $0$ by including the spin structure as suggested by \cite{ENS02,karlsson2017note}. Nevertheless, by Proposition \ref{distinguishing charts},  the cluster structure in characteristic $0$ will \emph{not} distinguish more fillings than characteristic $2$. For the purpose of building an invariant for Lagrangian fillings from cluster theory, it is enough to consider characteristic $2$.
Meanwhile, for Proposition \ref{distinguishing charts} to apply, augmentations must be defined over an \emph{algebraically closed} field. 
\end{rmk}

\subsection{Finite Type Classifications}

Recall that an ADE quiver is a directed graph whose underlying graph is one of the ADE Dynkin diagrams. The quivers $Q_\beta$ for different words $\beta$ of $[\beta]$ are mutation equivalent, leading us to the following definition. 
\begin{defn}\label{finite type braid def} 
A positive braid $[\beta]$ is of \emph{finite type} if the unfrozen part of $Q_\beta$ is mutation equivalent to a disjoint union of ADE quivers for one (equivalently any) word $\beta$ of  $[\beta]$. Otherwise,  $[\beta]$ is  of \emph{infinite type}.
\end{defn}

\begin{defn}\label{stad.links}  The Legendrian links $\Lambda_\beta$ associated with the positive braid words $\beta$ in the following table are called the \emph{standard ADE links}.
\begin{center}
    \begin{tabular}{|c|c|c|c|c|} \hline
    \rule{0pt}{2.3ex} $\Br_2^+$ & \multicolumn{4}{c|}{$\Br_3^+$} \\[0.045cm] \hline
    \rule{0pt}{2.3ex} $\mathrm{A}_r$     &  $\mathrm{D}_r$ & $\mathrm{E}_6$  & $\mathrm{E}_7$ & $\mathrm{E}_8$ \\[0.045cm] \hline
    \rule{0pt}{2.3ex} $s_1^{r+1}$ &     $ s_1^{r-2}s_2s_1^2s_2$ & $s_1^3s_2s_1^3s_2$ & $s_1^4s_2s_1^3s_2$ & $s_1^5s_2s_1^3s_2$ \\[0.045cm] \hline
    \end{tabular}
\end{center}
\end{defn}

\begin{rmk}
The underlying smooth links of the above standard Legendrian ADE links are the same as links of certain plane curves singularities as in \cite{Arn}. Namely, they coincide with the intersections $B_\ep(0,0)\cap V_f$, where $V_f$ is the vanishing locus of $f(x,y): \bbC^2\rightarrow \bbC$ given by
$$\mathrm{A}_r: x^{r+1} + y^2,\quad D_r: x^2y+y^{r-1}, \quad \mathrm{E}_6: x^3+y^4, \quad \mathrm{E}_7: x^3+xy^3,\quad \mathrm{E}_8: x^3+y^5.$$
Another topological description of this class is prime positive braid links with positive-definite symmetric Seifert forms \cite{Baa}.
\end{rmk}

The next result provides several characterizations of positive braids of finite type. 

\begin{thm} Let $\Lambda_\beta$ be the Legendrian link associated with a positive braid $\beta$. The following statements are equivalent.
\begin{enumerate}
    \item $[\beta]$ is of finite type.
    \item ${\rm Aug}(\Lambda_\beta)$ has a finite number of cluster seeds.
    \item $\Lambda_\beta$ is Legendrian isotopic to a split union of unknots and connect sums of standard $\mathrm{ADE}$ links.
    \item The symmetric Seifert form of $\Lambda_\beta$ is positive definite. 
\end{enumerate}
\end{thm}

\begin{proof}
The equivalence ``(1) $\Leftrightarrow$ (2)" follows from the finite classification of cluster algebras \cite{FZII} and Theorem \ref{Mainthm1}. The implications ``(3) $\Rightarrow$ (4) $\Rightarrow$ (1)" are a result of \cite{Baa}. We prove ``(1) $\Rightarrow$ (3)" in Theorem \ref{4.6}.
\end{proof}

\subsection{Infinitely Many Exact Lagrangian Fillings}
Our last main result is as follows.
\begin{thm}[Theorem \ref{MainStep1}] \label{thm.main3} If $[\beta]$ is of infinite type, then $\Lambda_\beta$ admits infinitely many non-Hamiltonian isotopic exact Lagrangian fillings. 
\end{thm}

The proof of Theorem \ref{thm.main3} uses the aperiodicity of some cluster Donaldson-Thomas transformations \cite{SWflag} and a trichotomy of the frieze variety \cite{lee2018frieze}. The DT transformation is not generally aperiodic (Remark \ref{periodicDT}), and we employ sophisticated combinatorial arguments to solve the problem.

As a topological consequence, this theorem yields that most positive braid links admit infinitely many non-Hamiltonian isotopic Lagrangian fillings. Hence it is reasonable to conjecture that Legendrian links with infinitely many fillings exist more broadly than those with finitely many fillings, whenever fillings are unobstructed. This theorem also motivates and proves a major class in the conjecture of ADE classification of Lagrangian fillings proposed in a later paper by Casals \cite{casals2020lagrangian}. It also provides interesting examples of concordance monoids, group of Legendrian loops, and Weinstein manifolds, following the framework of \cite{CasalsGao}.

\subsection*{Acknowledgements} We would like to thank Roger Casals, Matthew Hedden, Cecilia Karlsson, Lenhard Ng, Dan Rutherford, Kevin Sackel, and Eric Zaslow for the helpful discussion on Floer theory, contact geometry, and Legendrian links. We would like to thank Alexander Goncharov, Bernhard Keller, and Michael Shapiro for the helpful discussion on the cluster theoretical aspects of the paper. We would like to thank Vladimir Retakh for the helpful discussion on quasi-determinants and non-commutative cluster algebras. HG thanks Tao Su for discussing a perspective from Morse complex families. We would also like to thank Roger Casals, Kevin Sackel, and Bernhard Keller for their valuable suggestions and comments on the first draft of this paper. LS is partially supported by the Collaboration Grant for Mathematicians from the Simons Foundation ($\#$711926) and the NSF grant DMS-2200738.

\bigskip

\section{Cluster \texorpdfstring{${\rm K}_2$ Structure on Augmentation Varieties}{}}
\subsection{Positive braid Legendrian Links} 
\label{sec2.1.positive-braid}
Artin's braid group on $n$ strands is 
$$\Br_n = \langle s_1^{\pm 1}, \dotsb,s_{n-1}^{\pm 1} \;|\; s_is_{i+1}s_i = s_{i+1}s_i s_{i+1} \text{, and } s_j s_k = s_k s_j \text{ if } |j-k|\geq 2 \rangle. $$
The \emph{positive braid semigroup} $\Br_n^+$ is the sub-semigroup inside $\Br_n$ generated by the $s_i$'s. 
The positive braid
$w_0 = (s_{1}\dotsb s_{n-1})(s_{1}\dotsb s_{n-2})\dotsb(s_{1}s_{2})(s_{1})$
is called the \emph{half twist}, 
and its square $w_0^2$ is the \emph{full twist}.  Under the quotient map from $\Br_n$ to the symmetric group $\mathsf{S}_n$, $w_0$ becomes the element of the longest Coxeter length.

We denote a word of a positive braid by $\beta$, and its equivalence class by $[\beta]$.
Every positive braid word $\beta$ uniquely determines a Legendrian link $\Lambda_\beta$ with maximal Thurston-Bennequin number in its smooth isotopy class \cite[Theorem 3.4]{EV}. The Legendrian $\Lambda_\beta$ can be obtained via a satellite construction, that is, the braid closure of $w_0 \beta w_0$ satellited along the standard unknot, $|w_0\beta w_0|\subset J^1(S^1)\subset \bbR^3$, produces the Legendrian embedding $\Lambda_\beta$. Alternatively, the front projection of $\Lambda_\beta$ is given via the rainbow closure construction \cite{STZ}. We apply Ng's resolution {\cite[Proposition 2.2]{NG2003}} to obtain its Lagrangian projection $\pi_L(\Lambda_\beta)$ as follows, where the left cusps are smoothed out and the right cusps are resolved to a crossing attaching a teardrop loop. Note that the Lagrangian projection is not drawn in scale -- the teardrop loop should be drawn much larger, so that the signed area on the left of a resolved crossing equals to the area of the teardrop.

\begin{figure}[H]
    \centering
    \begin{tikzpicture}[baseline=30,scale =0.5]
    \draw (0,0) rectangle node [] {$\beta$} (3,1.5);
    \draw  [decoration={markings,mark=at position 0.52 with {\arrow{>}}},postaction={decorate}]  (3,0.25)  to [out=0,in=180] (5.5,2) to [out=180,in=0] (3, 3.75) -- (0,3.75) to [out=180,in=0] (-2.5,2) to [out=0,in=180] (0,0.25);
    \draw (3,1) to [out=0,in=180] (4.5,2) to [out=180,in=0] (3, 3) -- (0,3) to [out=180,in=0] (-1.5,2) to [out=0,in=180] (0,1);
    \draw (3,1.25)  to [out=0,in=180] (4,2) to [out=180,in=0] (3, 2.75) -- (0,2.75) to [out=180,in=0] (-1,2) to [out=0,in=180] (0,1.25);
    \node at (1.5,3.5) [] {$\vdots$};
    \node at (1.5,-0.75) [] {$\pi_F(\Lambda_\beta)$};
    \end{tikzpicture}
        $\quad \rightsquigarrow \quad$
    \begin{tikzpicture}[baseline=30,scale =0.5]
    \draw (0,0) rectangle node [] {$\beta$} (3,1.5);
    \draw (3,1.25) to [out=0,in=-135] (4,2) to [out=45,in=90] (4.5,2) to [out=-90,in=-45] (4,2) to [out=135,in=0] (3,2.75) -- (0,2.75) to [out=180,in=90] (-1,2) to [out=-90,in=180] (0,1.25);
    \draw (3,1) to [out=0,in=-135] (5,2) to [out=45,in=90] (5.5,2) to [out=-90,in=-45] (5,2) to [out=135,in=0] (3,3) -- (0,3) to [out=180,in=90] (-1.25,2) to [out=-90,in=180] (0,1);
    \draw [decoration={markings,mark=at position 0.6 with {\arrow{>}}},postaction={decorate}] (3,0.25) to [out=0,in=-135] (6.5,2) to [out=45,in=90] (7,2) to [out=-90,in=-45] (6.5,2) to [out=135,in=0] (3,3.75) -- (0,3.75) to [out=180,in=90] (-2.25,2) to [out=-90,in=180] (0,0.25);
    \node at (6,2) [] {$\cdots$};
    \node at (-1.75,2) [] {$\cdots$};
    \node at (1.5,3.5)[] {$\vdots$};
    \node at (1.5,-0.75) [] {$\pi_L(\Lambda_\beta)$};
    \end{tikzpicture}.
    \caption{Ng's resolution}
     \label{rainbow}
\end{figure}

\subsection{Augmentation Varieties for Positive Braid Legendrian Links}

In this section we compute the CE dga $\cA\left(\Lambda_\beta\right)$ and the augmentation variety ${\rm Aug}(\Lambda_\beta)$ for a positive braid word $\beta=s_{i_1}\dotsb s_{i_l}$ with $n$ strands. We refer the readers to \cite{ENsurvey} for the definition of CE dga for general Legendrians.  For postive braids, K\'alm\'an \cite{Kalman} had explicitly computed the dga with $\bbZ_2$-coefficient. We recover the computation with a different method using the boarded dga  in \cite{Siv10}. The coefficients are enhanced to include contributions from marked points.

Let $\pi_L(\Lambda_\beta)$ be  the Lagrangian projection of $\Lambda_\beta$ as in Figure \ref{rainbow}. The Reeb chords of $\Lambda_\beta$ correspond to crossings in $\pi_L(\Lambda_\beta)$. We equip $\Lambda_\beta$ with a binary Maslov potential $\{0,1\}$, which determines degrees for the Reeb chords. The crossings in the braid have degree $0$ and are denoted by $b_1,\ldots, b_l$. The crossings located at resolved right cusps have degree $1$ and are denoted by $a_1,\ldots, a_n$.  We decorate $\Lambda_\beta$ by placing a marked point $t_i$ next to each crossing $a_i$, located on the resolved teardrop loop. Let $\cT$ be the set of marked points. The dga $\cA\left(\Lambda_\beta\right)$ is generated by the Reeb chords and the formal variables $t_{i}^{\pm 1}$.  The non-trivial differentials of $\cA\left(\Lambda_\beta\right)$ are the $\partial a_k$'s, which we shall describe.

For any noncommutative formal variable $b$ and  $1\leq i <n$, we define an $n\times n$ matrix
\begin{equation}\label{SwitchMat}
Z_i(b): = 
\begin{pmatrix}
1  & &  & & &\\
& \ddots  & && &\\
& &b & 1 & &\\
& &1 & 0 & &\\
& & & & \ddots &\\
& & & & & 1 \\
\end{pmatrix}, 
\end{equation}
where the $2\times 2$ sub-matrix sits at the $i$th and $(i+1)$st rows and columns. For reference, this matrix is called the path matrix in \cite{Kalman06}.
For $\beta=s_{i_1}\ldots s_{i_l}$, let us set $M^{(1)}:=Z_{i_1}\left(b_1\right)\cdots Z_{i_l}\left(b_l\right)$. Define the matrices $M^{(k)}=\left(M^{(k+1)}_{ij}\right)_{k\leq i,j\leq n}$ recursively by
\begin{equation}
\label{qcdondoqw}
M^{(k+1)}_{ij}=M^{(k)}_{ij}+M^{(k)}_{ik}t_kM^{(k)}_{kj}.
\end{equation}

\begin{prop}\label{partial a} 
The differential of the CE-dga $\mathcal{A}(\Lambda_\beta)$ has the following compact form:
\[
\partial a_k=M^{(k)}_{kk}+t_k^{-1}, \hskip 7mm \forall 1\leq k\leq n.
\]
\end{prop}
\begin{proof} Borrowing the idea of the bordered dga \cite{Siv10}, we consider the diagram
\[
\begin{tikzpicture}[baseline=15, scale=.8]
\draw (0,0.25) rectangle node [] {$\beta$} (4,1.3);
\draw (4,1.25) to [out=0,in=-135] (5,1.75) to [out=45,in=90] (5.5,1.75) to [out=-90,in=-45] (5,1.75) to [out=135,in=0] (4,2.25) -- (0,2.25) to [out=180,in=90] (-1,1.75) to [out=-90,in=180] (0,1.25);
\draw (4,1) to [out=0,in=-135] (6,1.75) to [out=45,in=90] (6.5,1.75) to [out=-90,in=-45] (6,1.75) to [out=135,in=0] (4,2.5) -- (0,2.5) to [out=180,in=90] (-1.25,1.75) to [out=-90,in=180] (0,1);
\draw [decoration={markings,mark=at position 0.5 with {\arrow{>}}},postaction={decorate}] (4,0.5) to [out=0,in=-135] (7.5,1.75) to [out=45,in=90] (8,1.75) to [out=-90,in=-45] (7.5,1.75) to [out=135,in=0] (4,3.25) -- (0,3.25) to [out=180,in=90] (-2.25,1.75) to [out=-90,in=180] (0,0.5);
\node at (7,1.75) [] {$\cdots$};
\node at (-1.75,1.75) [] {$\cdots$};
\node at (2,3)[] {$\vdots$};
\draw [dashed] (-3,1.25) node [left] {\footnotesize{$1$st}} -- (6,1.25);
\draw [dashed] (-3,1) node [left] {\footnotesize{$2$nd}} -- (6,1);
\draw [dashed] (-3,0.5) node [left] {\footnotesize{$n$th}} -- (6,0.5);
\node at (-2.5,0.9) [] {$\vdots$};
\draw [dashed] (-3,2.25) node [left] {\footnotesize{$1$st}} -- (6,2.25);
\draw [dashed] (-3,2.5) node [left] {\footnotesize{$2$nd}} -- (6,2.5);
\draw [dashed] (-3,3.25) node [left] {\footnotesize{$n$th}} -- (6,3.25);
\node at (-2.5,3) [] {$\vdots$};
\draw [dashed] (4.25,0) node [below] {$M^{(1)}$} -- (4.25,3.5);
\draw [dashed] (6.65,0) node [below] {$M^{(k)}$} -- (6.65,3.5);
\end{tikzpicture}.
\]
Let us label a dashed line between the braid region and the right cusps so that each disk contributing to the differential can be divided into two parts. On the left, each disk boundary will travel along a strand on the top, making many or no turns in the braid region, and then hit the dashed line. In general, the disk configuration near a resolved right cusp can be one of the following: 
\[
\begin{tikzpicture}
\path [fill=lightgray] (-0.707,-1) rectangle (1,1);
\draw (135:1) to [out=-45,in=135] (0,0) to [out=-45,in=180] (0.5,-0.25) to [out=0,in=0] (0.5,0.25) to [out=180,in=45] (0.1,0.1); 
\draw (-0.1,-0.1) to [out=-135,in=45] (-135:1);
\end{tikzpicture}
\quad \quad
\begin{tikzpicture}
\path [fill=lightgray] (-0.707,-1) -- (-135:1) to [out=45,in=-135] (0,0) to [out=-45,in=-180] (0.5,-0.25) to [out=0,in=0] (0.5,0.25) to [out=180,in=45] (0,0) to [out=135,in=-45] (135:1) -- (-0.707,1) -- (1,1) -- (1,-1) --cycle;
\draw (135:1) to [out=-45,in=135] (0,0) to [out=-45,in=180] (0.5,-0.25) to [out=0,in=0] (0.5,0.25) to [out=180,in=45] (0.1,0.1); 
\draw (-0.1,-0.1) to [out=-135,in=45] (-135:1);
\end{tikzpicture}
\quad \quad
\begin{tikzpicture}
\path [fill=lightgray] (-0.707,-1) -- (-0.707,-0.757) to [out=45,in=-135] (0,-0.05) to [out=-45,in=-180] (0.5,-0.25) to [out=0,in=0] (0.5,0.25) to [out=180,in=45] (0,0.05) -- (-0.707,-0.657) -- (-0.707,1) -- (1,1) --  (1,-1) --cycle;
\draw (135:1) to [out=-45,in=135] (0,0) to [out=-45,in=180] (0.5,-0.25) to [out=0,in=0] (0.5,0.25) to [out=180,in=45] (0.1,0.1); 
\draw (-0.1,-0.1) to [out=-135,in=45] (-135:1);
\end{tikzpicture}
\quad \quad
\begin{tikzpicture}
\path [fill=lightgray] (-0.707,-1) -- (-0.707,0.657) to [out=-45,in=135] (0,-0.05) to [out=-45,in=-180] (0.5,-0.25) to [out=0,in=0] (0.5,0.25) to [out=180,in=45] (0,0.05) to [out=135,in=-45] (-0.707,0.757) -- (-0.707,1) -- (1,1) --  (1,-1) --cycle;
\draw (135:1) to [out=-45,in=135] (0,0) to [out=-45,in=180] (0.5,-0.25) to [out=0,in=0] (0.5,0.25) to [out=180,in=45] (0.1,0.1); 
\draw (-0.1,-0.1) to [out=-135,in=45] (-135:1);
\end{tikzpicture}
\quad \quad
\begin{tikzpicture}
\path [fill=gray] (-135:1) -- (0,0) -- (135:1) --cycle;
\path [fill=lightgray] (-0.707,-1) -- (-135:1) to [out=45,in=-135] (0,0) to [out=-45,in=-180] (0.5,-0.25) to [out=0,in=0] (0.5,0.25) to [out=180,in=45] (0,0) to [out=135,in=-45] (135:1) -- (-0.707,1) -- (1,1) -- (1,-1) --cycle;
\draw (135:1) to [out=-45,in=135] (0,0) to [out=-45,in=180] (0.5,-0.25) to [out=0,in=0] (0.5,0.25) to [out=180,in=45] (0.1,0.1); 
\draw (-0.1,-0.1) to [out=-135,in=45] (-135:1);
\end{tikzpicture}.
\]
In our setup, only the first and the last configurations occur.

Let us start with the left part. Suppose the $(i,j)$-th entry of $M$ counts disks that are bounded by top level $i$ and bottom level $j$ near the dashed line. It can be computed inductively on crossings from left to right. Before the braiding region, there is a unique pairing between the strands, giving the identity matrix. For an arbitrary crossing $i_k$, let $N$ (resp. $N'$) be the disk counting matrix before (resp. after) scanning across $i_k$. Then 
\begin{itemize}
    \item $N'_{i_k+1,j}=N_{i_k,j}$ due to (a); and $N'_{i_k,j}=N_{i_k,j}b_{k}+N_{i_k+1,j}$ due to (b) and (c).
\end{itemize}
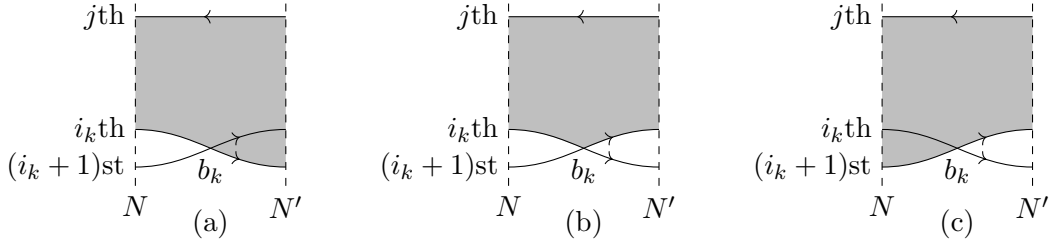
\[
\begin{tikzpicture}[scale=.8]
\path [fill=lightgray] (0,0.5) to [out=0,in=180] (2,0) -- (2,2) -- (0,2) --cycle;
\draw [decoration={markings,mark=at position 0.5 with {\arrow{<}}},postaction={decorate}] (0,2) node [left] {$j$th} -- (2,2);
\draw [decoration={markings,mark=at position 0.7 with {\arrow{>}}},postaction={decorate}] (0,0) node [left] {$(i_k+1)$st} to [out=0,in=180] (2,0.5);
\draw [decoration={markings,mark=at position 0.7 with {\arrow{>}}},postaction={decorate}] (0,0.5) node [left] {$i_k$th} to [out=0,in=180] (2,0);
\node at (1,0.25) [below] {$b_k$};
\draw [dashed] (0,-0.25) node [below] {$N$} -- (0,2.25);
\draw [dashed] (2,-0.25) node [below] {$N'$} -- (2,2.25);
\node at (1, -1.25) {(a)};
\end{tikzpicture}
\quad \quad
\begin{tikzpicture}[scale=.8]
\path [fill=lightgray] (0,0.5) to [out=0,in=160] (1,0.25) to [out=20,in=180] (2,0.5) -- (2,2) -- (0,2) --cycle;
\draw [decoration={markings,mark=at position 0.5 with {\arrow{<}}},postaction={decorate}] (0,2) node [left] {$j$th} -- (2,2);
\draw [decoration={markings,mark=at position 0.7 with {\arrow{>}}},postaction={decorate}] (0,0) node [left] {$(i_k+1)$st} to [out=0,in=180] (2,0.5);
\draw [decoration={markings,mark=at position 0.7 with {\arrow{>}}},postaction={decorate}] (0,0.5) node [left] {$i_k$th} to [out=0,in=180] (2,0);
\node at (1,0.25) [below] {$b_k$};
\draw [dashed] (0,-0.25) node [below] {$N$} -- (0,2.25);
\draw [dashed] (2,-0.25) node [below] {$N'$} -- (2,2.25);
\node at (1, -1.25) {(b)};
\end{tikzpicture}
\quad \quad
\begin{tikzpicture}[scale=.8]
\path [fill=lightgray] (0,0) to [out=0,in=180] (2,0.5) -- (2,2) -- (0,2) --cycle;
\draw [decoration={markings,mark=at position 0.5 with {\arrow{<}}},postaction={decorate}] (0,2) node [left] {$j$th} -- (2,2);
\draw [decoration={markings,mark=at position 0.7 with {\arrow{>}}},postaction={decorate}] (0,0) node [left] {$(i_k+1)$st} to [out=0,in=180] (2,0.5);
\draw [decoration={markings,mark=at position 0.7 with {\arrow{>}}},postaction={decorate}] (0,0.5) node [left] {$i_k$th} to [out=0,in=180] (2,0);
\node at (1,0.25) [below] {$b_k$};
\draw [dashed] (0,-0.25) node [below] {$N$} -- (0,2.25);
\draw [dashed] (2,-0.25) node [below] {$N'$} -- (2,2.25);
\node at (1, -1.25) {(c)};
\end{tikzpicture}
\]
In other words, $N'=NZ_{i_k}\left(b_k\right)$. By induction, we have $M=M^{(1)}=Z_{i_1}\left(b_1\right)\dotsb Z_{i_l}\left(b_l\right)$.

Similarly, we place a dashed line between each pair of right cusps. Let $M^{(k)}$ be the matrix associated to the dashed line between $a_{k-1}$ and $a_k$. There is no disk between any two top strands or between any two bottom strands near $M^{(1)}$ dashed line, and will be inductively true for any other dashed lines. 
\[
\begin{tikzpicture}[baseline=20]
\path[fill=lightgray] (0,0) -- (2,0) -- (2,2) -- (0,2) -- cycle;
\draw [decoration={markings,mark=at position 0.5 with {\arrow{>}}},postaction={decorate}] (0,0) node [left] {$j$th} -- (2,0);
\draw [decoration={markings,mark=at position 0.4 with {\arrow{>}}},postaction={decorate}] (0,0.5) node [left] {$k$th} to [out=0,in=-135] (1,1) to [out=45,in=90] (1.5,1) to [out=-90,in=-45] (1,1) to [out=135,in=0] (0,1.5) node [left] {$k$th};
\draw [decoration={markings,mark=at position 0.5 with {\arrow{<}}},postaction={decorate}] (0,2) node [left] {$i$th} -- (2,2);
\node at (1,1) [left] {$a_k$};
\node at (1.5,1) [] {$\ast$};
\node at (1.5,1) [right] {$t_k$};
\end{tikzpicture}
\quad\quad\quad\quad
\begin{tikzpicture}[baseline=20]
\path[fill=lightgray] (0,0) -- (2,0) -- (2,2) -- (0,2) -- (0,1.5) to [out=0,in=135] (1,1) to [out=45,in=90] (1.5,1) to [out=-90,in=-45] (1,1) to [out=-135,in=0] (0,0.5) -- cycle;
\path[fill=gray] (0,1.5) to [out=0,in=135] (1,1) to [out=-135,in=0] (0,0.5) --cycle;
\draw [decoration={markings,mark=at position 0.5 with {\arrow{>}}},postaction={decorate}] (0,0) node [left] {$j$th} -- (2,0);
\draw [decoration={markings,mark=at position 0.4 with {\arrow{>}}},postaction={decorate}] (0,0.5) node [left] {$k$th} to [out=0,in=-135] (1,1) to [out=45,in=90] (1.5,1) to [out=-90,in=-45] (1,1) to [out=135,in=0] (0,1.5) node [left] {$k$th};
\draw [decoration={markings,mark=at position 0.5 with {\arrow{<}}},postaction={decorate}] (0,2) node [left] {$i$th} -- (2,2);
\node at (1,1) [left] {$a_k$};
\node at (1.5,1) [] {$\ast$};
\node at (1.5,1) [right] {$t_k$};
\end{tikzpicture}
\]
Enumerating the two local situations above, we have
\[
M^{(k+1)}_{ij} = M^{(k)}_{ij}+M^{(k)}_{ik}t_kM^{(k)}_{kj}.
\]

We are ready to compute $\partial a_k$. It counts two types of disks. One consists of those disks that hit the dashed line labeled by $M^{(k)}$, and the other one consists of only one disk given by the teardrop loop with no negative punctures. Hence,
$
\partial a_k=M^{(k)}_{kk}+t_k^{-1}. $
\end{proof}

Gelfand and Retakh \cite{GR} introduced the {\it quasi-determinant} as a replacement for the determinant for matrices with noncommutative entries. Let $M_{1,2,\ldots, k-1, i}^{1,2,...,k-1,j}$ be the $k\times k$ submatrix of $M=M^{(1)}$ consisting of rows $1,2, \ldots, k-1,i$ and columns $1,2,\ldots, k-1,j$. The next proposition establishes a connection between $M_{ij}^{(k)}$ and the quasi-determinants.

\begin{prop} 
\label{cnjona}
If $\partial a_k=0$ for $1\leq k\leq n$, then $M_{ij}^{(k)}$ is the quasi-determinant $\left|M_{1,2,\dots, k-1,i}^{1,2,\dots, k-1,j}\right|_{ij}$. 
\end{prop}
\begin{proof}
The assumption $\partial a_k=0$ implies that $t_k=-\left(M_{kk}^{(k)}\right)^{-1}$. Then \eqref{qcdondoqw} becomes
\[
M_{ij}^{(k+1)}=M_{ij}^{(k)}-M_{ik}^{(k)}\left(M_{kk}^{(k)}\right)^{-1}M_{kj}^{(k)}.
\]
Inductively, the RHS equals $\left|M_{1,\dots, k-1,i}^{1,\dots,k-1,j} \right|_{ij}-\left|M_{1,\dots, k-1,i}^{1,\dots,k-1,k} \right|_{ik}\left|M_{1,\dots, k-1,k}^{1,\dots,k-1,k} \right|_{kk}^{-1}\left|M_{1,\dots, k-1,k}^{1,\dots,k-1,j} \right|_{kj}$, which yields $\left|M_{1,\dots, k,i}^{1,\dots,k,j} \right|_{ij}$ by the Sylvester's identity for quasi-determinants (Proposition 1.5 of \cite{GR}).
\end{proof}

The dga $\cA(\Lambda_\beta)$ is concentrated at non-negative degrees. Its homology $H_0\left(\cA\left(\Lambda_\beta\right)\right)$ is a non-commutative algebra.  Let us write $M_k:=M_{kk}^{(k)}$ for short. 
The following result is a direct consequence of Propositions \ref{partial a}  and \ref{cnjona}.
\begin{cor}\label{H0=quasi-det} As non-commutative algebras over $\mathbb{Z}_2$, we have \[H_0\left(\cA\left(\Lambda_\beta\right)\right)\cong \displaystyle\frac{\mathbb{Z}_2\left\langle b_1,\dots, b_l, t_1^{\pm 1},\dots, t_n^{\pm 1}\right\rangle}{\left(M_k=t_k^{-1}\right)}.\]
\end{cor}

Now let us fix an algebraic closed field $\mathbb{F}$ of characteristic 2. 
\begin{defn}\label{def aug var}  An \emph{augmentation} of $\cA(\Lambda_\beta)$ is a unital dga homomorphism
$$\varepsilon: (\cA(\Lambda_\beta),\del) \rightarrow (\mathbb{F},0).$$
The \emph{augmentation variety} $\Aug\left(\Lambda_\beta\right)$ is the moduli space of augmentations of $\cA\left(\Lambda_\beta\right)$.
\end{defn}

\begin{rmk}
The augmentation variety is different from the moduli stack of objects in the unital augmentation category introduced in \cite{NRSSZ}. Therefore it is not isomorphic to the moduli space of microlocal rank one sheaves associated to $\Lambda$ in general.
\end{rmk}

\begin{lem}\label{aug affine} For any positive braid word $\beta$, let $\cA(\Lambda_\beta)^c$ be the abelianization of $\cA(\Lambda_\beta)$. Then $\Aug(\cA(\Lambda_\beta))$ is an affine variety whose coordinate ring is $H_0\left(\cA(\Lambda_\beta)^c,\mathbb{F}\right)$.
\end{lem}
\begin{proof}
By definition, the augmentations $\varepsilon$ preserve the degree. In particular,  $\varepsilon(a) =0$ for any generator $a$ of non-zero degree. Hence,  $\varepsilon$
is uniquely determined by its evaluations at the degree zero Reeb chords and the formal variables, and the evaluations are subject to the conditions $\varepsilon\circ \partial (a)=0$ for any degree 1 Reeb chord $a$. As a consequence, the augmentation varieties are affine varieties.

The field $\mathbb{F}$ is commutative. Thus, the augmentations for $\cA(\Lambda_\beta)$ and $\cA(\Lambda_\beta)^c$ coincide. Let $\del_i$ be the $i$-th degree of $\del$. Since $\cA(\Lambda_\beta)^c$ is concentrated in non-negative degrees, $\ker \del_0$ is the free algebra generated by the formal variables and the degree 0 Reeb chords, and $\mathrm{im}\, \del_1$ is an ideal generated by $(\partial a)$ for all degree 1 Reeb chords $a$. Hence, $H_0(\cA(\Lambda_\beta)^c, \mathbb{F}) = \ker \del_0/ \mathrm{im}\, \del_1$ is the coordinate ring of $\Aug(\cA(\Lambda_\beta))$. 
\end{proof}

\begin{defn} \label{cfqhbvo}
Let $N$ be an $n\times n$ matrix over $\mathbb{F}$. The $m$th \emph{principal minor} of $N$, denoted by $\Delta_m(N)$, is the determinant of the $m\times m$ submatrix of $N$ formed by the first $m$ rows and columns.
\end{defn}

\begin{prop}\label{non-vanishing} The coordinate ring of $\Aug(\Lambda_\beta)$ is  
\[
\mathbb{F}\big[b_1,\dotsb, b_l,t_1^{\pm 1}, \dotsb, t_n^{\pm 1}\big]\Big/\cI,
\]
where the ideal $\cI$ is generated by 
        \begin{equation}\label{aug equ simplify}
            \Delta_m\left(Z_{i_1}(b_1)\ldots Z_{i_l}(b_l)\right)=\prod_{k=1}^m t_k^{-1}, \quad\quad 1\leq m\leq n.
        \end{equation}
The $\Aug\left(\Lambda_\beta\right)$ is the non-vanishing locus of the polynomial $\prod_{m=1}^n\Delta_m\left(Z_{i_1}(b_1)\ldots Z_{i_l}(b_l)\right)$ inside the ambient affine space $\mathbb{F}^l_{b_1,\dots, b_l}$.
\end{prop}
\begin{rmk}
Note that the Reeb chords $b_i$ can be regarded as coordinate functions on $\Aug\left(\Lambda_\beta\right)$.
 We call them the \emph{Reeb coordinates}.
\end{rmk}
\begin{proof} 
By Lemma \ref{aug affine}, we have $\Aug\left(\Lambda_\beta\right)=\spec H_0(\mathcal{A}\left(\Lambda_\beta\right)^c,\mathbb{F})$, where
\[
H_0\left(\mathcal{A}\left(\Lambda_\beta\right)^c,\mathbb{F}\right)=\left.\mathbb{F}\big[b_i,t_j^{\pm 1}\big]\right/\left(\partial^c (a_k) =0\right).
\]
Therefore the defining equations of $\Aug\left(\Lambda_\beta\right)$ are $\partial^c (a_k) =0$  for $k=1,\dotsb, n$.

In the commutative setting, the quasi-determinant reduces to ratio of determinants
\begin{equation}
\label{cabdiinw}
\left|N\right|_{ij}=(-1)^{i+j}\frac{\det N}{\det N^{ij}},
\end{equation}
where $N^{ij}$ is the minor that results from deleting  row $i$ and  column $j$ from $N$. We ignore the signs in the setting of characteristic 2. 
Using Lemma \ref{cnjona} and \eqref{cabdiinw} inductively,  $\partial^c(a_k)=0$ is equivalent to \eqref{aug equ simplify}, which concludes the proof of the first part.

 Note that $t_j^{\pm 1}$ are invertible. Using \eqref{aug equ simplify} recursively, each $t_k$ is can be expressed in terms of principal minors, and hence in terms of the coordinates $b_1, \ldots, b_l$. After eliminating all the formal variables $t_k^{\pm 1}$, we end up with the desired equation.
\end{proof}

\begin{cor}\label{aug vty smooth}
The augmentation variety $\Aug\left(\Lambda_\beta\right)$ is smooth.
\end{cor}
\begin{proof} By Proposition \ref{non-vanishing}, it is the non-vanishing locus of a polynomial function.
\end{proof}

\begin{cor}\label{prod t = 1} Any augmentation $\varepsilon$ on $\cA\left(\Lambda_\beta\right)$ satisfies $\prod_{k=1}^n \varepsilon\left(t_k\right)=1$.
\end{cor}
\begin{proof} Take $m=n$ in \eqref{aug equ simplify}. Then $\prod_{k=1}^nt_k^{-1}=\Delta_n(M)=\det (M)$. Since each commutative $Z_{i_k}\left(b_k\right)$ has determinant $1$, we have $\det(M)=1$. Therefore 
$
\prod_{k=1}^nt_k=1.$\qedhere
\end{proof}

\subsection{\texorpdfstring{Cluster $\mathrm{K}_2$ structures on augmentation varieties}{}}
A double Bott-Samelson cell $\conf^e_\beta(\mathcal{C})$ is a cluster ${\rm K}_2$ variety introduced in \cite{SWflag}. We recall its definition and cluster structure in  Appendix \ref{Sec 3}.
In this section, we construct a natural isomorphism between  $\Aug\left(\Lambda_\beta\right)$ and  $\conf^e_\beta(\mathcal{C})$, which endows $\Aug\left(\Lambda_\beta\right)$ with a cluster $\mathrm{K}_2$ structure. 

By Proposition \ref{caosno}, $\conf^e_\beta(\mathcal{C})$ is a scheme over $\mathbb{Z}$. We have

\begin{thm}\label{3.6} Let $\G=\SL_n$ and let $\beta$ be a positive braid word of $n$ strands. After a base-change of $\conf^e_\beta(\mathcal{C})$ to $\mathbb{F}$, there is a natural isomorphism as $\mathbb{F}$-varieties:
\[
\Aug\left(\Lambda_\beta\right) \overset{\gamma}{\longrightarrow} \conf^e_\beta(\mathcal{C}).
\]
The pull-back of the cluster $\mathrm{K}_2$ structure on $\conf^e_\beta(\mathcal{C})$ equips $\Aug(\Lambda_\beta)$ with a cluster $\mathrm{K}_2$ structure.
\end{thm}
\begin{proof}
Let $\left(b_1,\dots, b_l\right)$ be the Reeb coordinates of $\Aug\left(\Lambda_\beta\right)$  and let $\left(q_1,\dots, q_l\right)$ be the affine coordinates of $\conf^e_\beta(\mathcal{C})$ as  in Proposition \ref{caosno}. 
Let $\gamma$ be the isomorphism of the ambient affine spaces $\mathbb{F}^l_{b_1,\dots, b_l}$ and $\mathbb{F}^l_{q_1,\dots, q_l}$ such that  $q_k=b_k$ for $1\leq k \leq l$. Since $\mathbb{F}$ is of characteristic 2, the matrix $Z_{i_k}\left(b_k\right)$ in \eqref{SwitchMat} equals $R_{i_k}\left(q_k\right)$ in \eqref{R def}. 
Hence the non-vanishing locus of $\prod_{1\leq i\leq n}\Delta_i\left(R_{i_1}\left(q_1\right)\cdots R_{i_l}\left(q_l\right)\right)$ coincides with that of $\prod_{1\leq i\leq n}\Delta_i\left(Z_{i_1}\left(b_1\right)\cdots Z_{i_l}\left(b_l\right)\right)$. By Propositions \ref{caosno} and \ref{non-vanishing}, these two non-vanishing loci are $\conf^e_\beta(\mathcal{C})$  and $\Aug\left(\Lambda_\beta\right)$  respectively. Therefore $\gamma$ restricts to an isomorphism between the two $\mathbb{F}$-varieties. 
\end{proof}

\bigskip

\section{From fillings to clusters}

\subsection{Exact Lagrangian Cobordisms and Enhanced EHK Functors}

 Ekholm, Honda, and K\'{a}lm\'{a}n \cite{EHK12} introduced a contravariant functor from the exact Lagrangian cobordism category of Legendrian links to the category of dga's. In this section, we discuss an enhancement of the EHK functor that includes decorations on exact Lagrangian cobordisms.  

Recall the standard contact $\mathbb{R}^3_{xyz}$ with the contact 1-form $\alpha={\rm d}z-y{\rm d}x$. Let $\mathbb{R}^4_{txyz}:=\mathbb{R}_t
\times \mathbb{R}^3_{xyz}$ be its symplectization with the symplectic form $\omega=d\left(e^t\alpha\right)$.

\begin{defn} Let $\Lambda_+$ and $\Lambda_-$ be two Legendrian links in  $\mathbb{R}^3_{xyz}$. An \emph{exact Lagrangian cobordism} $L:\Lambda_-\rightarrow \Lambda_+$ is an embedded oriented Lagrangian submanifold $L$ of $\mathbb{R}^4_{txyz}$ such that for some $N>0$,
\begin{enumerate}
    \item $L\cap\left((-\infty,-N]\times \mathbb{R}^3\right)=(-\infty,-N]\times \Lambda_-$ and $L\cap \left([N,\infty)\times \mathbb{R}^3\right)=[N,\infty)\times \Lambda_+$;
    \item there is a function $f$ of $L$, constant on $(-\infty,-N]\times \Lambda_-$ and $[N,\infty)\times \Lambda_+$, such that $df=\omega|_L$.
\end{enumerate}
An \emph{exact Lagrangian filling} of $\Lambda$ is an exact Lagrangian cobordism from $\emptyset$ to $\Lambda$. An \emph{exact Lagrangian concordance} is an exact Lagrangian cobordism that is topologically a cylinder.
\end{defn}

Exact Lagrangian fillings are central objects in contact and symplectic topology \cite{NZ,Nadler,Sylvan,GPSsectorial,EL}. These fillings induce augmentations \cite{EGH,ENsurvey}. Many, but not all, augmentations can be obtained from fillings. Note that the exact fillings of a Legendrian link have the same genus \cite{chantraine2010}. It is expected that their induced charts of augmentations have the same dimension.

 A \emph{$t$-minimum} on an exact Lagrangian cobordism $L$ is a point which achieves a local minimum for the coordinate function $t$ restricted on $L$. Denote by  $\cT_{\min}$ the set of $t$-minima.
Up to a Morse type perturbation, we will always assume that $L$ has finitely many isolated $t$-minima in the rest of this paper. 

\begin{defn} Let $\left(\Lambda_+, \cT_+\right)$ and $\left(\Lambda_-, \cT_-\right)$ be two decorated Legendrian links. A \emph{decorated exact Lagrangian cobordism} 
\[
(L,\cP):\left(\Lambda_-,\cT_-\right)\rightarrow \left(\Lambda_+,\cT_+\right)
\]
is an exact Lagrangian cobordism $L:\Lambda_-\rightarrow \Lambda_+$, together with a \emph{decoration} $\cP$, that is, a set of generic oriented marked curves $\cP=\left\{p_1,\dots, p_m\right\}$ on $L$, such that
\begin{enumerate}
    \item each $p_i$ is either a closed 1-cycle or an oriented curve that begins and ends at $\cT_+\cup \cT_-\cup \cT_{\min}$;
    \item intersections between these marked curves are transverse and isolated;
    \item each marked point in $\cT_+\cup \cT_-$ is the restriction of a unique marked curve $p_i$ to $\Lambda_-\sqcup \Lambda_+$.
\end{enumerate}
\end{defn}

Recall that the dga $\cA\left(\Lambda\right)$ of a decorated Legendrian $\left(\Lambda, \cT\right)$ is a non-commutative $\mathbb{Z}_2$-algebra freely generated by the set of Reeb chords $\cR$ and formal variables $\cT^{\pm 1}$. Given a decorated exact Lagrangian cobordism $(L,\cP):\left(\Lambda_-,\cT_-\right)\rightarrow \left(\Lambda_+,\cT_+\right)$, we define a pair of non-commutative $\mathbb{Z}_2$-algebras $\cA\left(\Lambda_\pm,\cP\right)$, each of which is generated by the respective set of Reeb chords $\cR_\pm$ and the formal variables $\cP^{\pm 1}$, modulo the following relations:
\begin{enumerate}
    \item$p_ip_j=p_jp_i$ if $p_i$ and $p_j$ intersect;
    \item near each $t$-minimum $\tau$ of $L$, let $\gamma$ be a small oriented loop around $\tau$, intersecting a collection of oriented marked curves cyclically, say $p_{i_1}, p_{i_2},\dots, p_{i_l}$; then 
    \begin{equation}\label{monodromy condition}
    p_{i_1}^{\inprod{\gamma}{p_{i_1}}}p_{i_2}^{\inprod{\gamma}{p_{i_2}}}\cdots p_{i_l}^{\inprod{\gamma}{p_{i_l}}}=1
    \end{equation}
    where $\inprod{\cdot}{\cdot}$ denotes the intersection number with respect to the orientation of $L$.
\end{enumerate}
The degrees of the Reeb chord generators of $\cA\left(\Lambda_\pm,\cP\right)$ are the same as those of $\cA\left(\Lambda_\pm\right)$. The degree of  $p_i^{\pm 1}$ is set to be 0. It makes $\cA\left(\Lambda_\pm,\cP\right)$ graded $\mathbb{Z}_2$-algebras.

We further define a pair of graded algebra homomorphisms
\[
\phi_{\pm}^*:\cA\left(\Lambda_\pm\right)\rightarrow \cA\left(\Lambda_\pm,\cP\right),
\]
each of which sends the Reeb chord generators to themselves and sends the marked points $t$  to
\[
\phi_\pm^*(t)=\left\{\begin{array}{ll}
    p^{\pm 1} & \text{if $p$ starts from $t$}, \\
    p^{\mp 1} & \text{if $p$ ends at $t$}. 
\end{array}\right.
\]
Let $\partial_\pm$ be the differentials on  $\cA\left(\Lambda_\pm\right)$. By defining the differentials on $\cA\left(\Lambda_\pm,\cP\right)$ to be $\phi_{\pm}^*\circ \partial_\pm$, we make $\cA\left(\Lambda_\pm,\cP\right)$ into a pair of dga's over $\mathbb{Z}_2$. 

\begin{defn} We call $\cA\left(\Lambda_\pm,\cP\right)$ the \emph{enhanced CE dga's} for $\Lambda_\pm$ with respect to the decorated exact Lagrangian cobordism $(L,\cP)$. 
\end{defn}

\begin{rmk} For an exact Lagrangian concordance $L:\Lambda_-\rightarrow \Lambda_+$ coming from a Legendrian isotopy, with a decoration $\cP$ coming from the trace of marked points, the dga homomorphisms $\phi_\pm^*$ are isomorphisms between $\cA\left(\Lambda_\pm\right)$ and $\cA\left(\Lambda_\pm,\cP\right)$. 
\end{rmk}

Two decorations $\cP$ and $\cP'$ on the same exact Lagrangian cobordism $L$ are \emph{equivalent} if the two sets of oriented marked curves can be related by a sequence of path homotopy and orientation reversing. Note that if $p\in \cP$ and $p'\in \cP'$ have the same underlying path but opposite orientation, the change of variable $p\leftrightarrow p'^{-1}$ gives rise to a natural dga isomorphism $\cA\left(\Lambda_\pm,\cP\right)\cong \cA\left(\Lambda_\pm,\cP'\right)$. Therefore, for the rest of this paper, we no longer distinguish equivalent decorations on the same exact Lagrangian cobordism.

Given two decorated exact Lagrangian cobordisms
\[
\left(\Lambda_0,\cT_0\right)\xrightarrow{\left(L_{01},\cP_{01}\right)} \left(\Lambda_1,\cT_1\right)\xrightarrow{\left(L_{12},\cP_{12}\right)} \left(\Lambda_2,\cT_2\right),
\]
we can compose them by concatenation (possibly with orientation reversing on some elements of the decorations) and get a decorated exact Lagrangian cobordism 
\[
(L,\cP):\left(\Lambda_0,\cT_0\right)\rightarrow \left(\Lambda_2,\cT_2\right).
\]
In particular, the resulting decoration $\cP$ is unique up to equivalence of decorations.

Let us now describe the enhancement of the EHK functor for the enhanced CE dga.

Let $(L,\cP):\left(\Lambda_-,\cT_-\right)\rightarrow \left(\Lambda_+,\cT_+\right)$ be a decorated exact Lagrangian cobordism. Let $J$ be a generic compatible tame almost complex structure on the symplectization $\mathbb{R}^4_{txyz}$. For $a\in \cR_+$ and $b_1,\dots, b_n\in \cR_-$, we define the moduli space $\cM\left(a;b_1,\dots, b_n\right)$
to be the set of bi-holomorphic equivalence classes of $J$-holomorphic curves, each with a positive puncture asymptotic to the strip over the Reeb chord $a$ at $+\infty$ and a negative puncture asymptotic to the strip over the Reeb chord $b_i$ at $-\infty$ for each $b_i$, appearing in the counterclockwise order along the boundary of the curve. For generic $J$, the moduli space $\cM\left(a,b_1,\dots, b_n\right)$ is a manifold of dimension $|a|-\sum_i|b_i|$ (see \cite[Lemma 3.7]{EHK12}).

For any $u\in \cM\left(a;b_1,\dots, b_n\right)$, the image of the disk boundary $\partial u$ is the disjoint union of $n+1$ oriented paths $\eta_0,\dots, \eta_n$ in the Lagrangian surface $L$. Suppose the path $\eta_i$ crosses oriented marked curves $p_{j_1}, p_{j_2},\dots, p_{j_l}$ in this particular order. We  define
\begin{equation}\label{w(u)}
p\left(\eta_i\right):=p_{j_1}^{\inprod{\eta_i}{p_{j_1}}}p_{j_2}^{\inprod{\eta_i}{p_{j_2}}}\cdots p_{j_l}^{\inprod{\eta_i}{p_{j_l}}} \quad \text{and}  \quad
w(u):=p\left(\eta_0\right)b_1p\left(\eta_1\right)b_2\cdots b_np\left(\eta_n\right).
\end{equation}
Following \cite{EHK12}, we define the dga homomorphism $\Phi_L^*:\cA\left(\Lambda_+,\cP\right)\rightarrow \cA\left(\Lambda_-,\cP\right)$ such that 
\[
\Phi_L^*(a) =\sum_{\substack{b_1,\dots, b_n\in \cR_-\\ \dim \cM\left(a;b_1,\dots, b_n\right)=0}}\sum_{u\in \cM\left(a;b_1,\dots, b_n\right)}w(u) \quad  \forall a\in \cR_+, \]
and
\[\Phi_L^*(p) = p \qquad \forall p\in \cP\]
If $(L,\cP):\left(\Lambda_0,\cT_0\right)\rightarrow \left(\Lambda_2,\cT_2\right)$ is the composition of $\left(L_{01},\cP_{01}\right):\left(\Lambda_0,\cT_0\right)\rightarrow \left(\Lambda_1,\cT_1\right)$ and $\left(L_{12},\cP_{12}\right):\left(\Lambda_1,\cT_1\right)\rightarrow \left(\Lambda_2,\cT_2\right)$, then the functorial homomorphisms can be composed via the following commutative diagram:
\[
\xymatrix{\cA\left(\Lambda_2\right) \ar[d] & & \cA\left(\Lambda_1\right)\ar[dl]\ar[dr]  & & \cA\left(\Lambda_2\right) \ar[d] \\
\cA\left(\Lambda_2,\cP_{12}\right)\ar[r]^{\Phi^*_{L_{12}}} \ar[d] &
\cA\left(\Lambda_1,\cP_{12}\right)\ar[dr] & & \cA\left(\Lambda_1,\cP_{01}\right)\ar[r]^{\Phi^*_{L_{01}}} \ar[dl] &
\cA\left(\Lambda_0,\cP_{01}\right)\ar[d] \\
\cA\left(\Lambda_2,\cP\right) \ar[rr]_{\Phi^*_{L_{12}}}  \ar@/_25pt/[rrrr]_{\Phi^*_{L}} & & \cA\left(\Lambda_1,\cP\right) \ar[rr]_{\Phi^*_{L_{01}}} & & \cA\left(\Lambda_0,\cP\right)
}.
\]
The  dga homomorphism $\Phi_L^*$ satisfies the following important property.

\begin{thm}[{\cite[Lemma 3.13]{EHK12}}]\label{thm 3.8} Suppose $L$ and $L'$ are Hamiltonian isotopic exact Lagrangian cobordisms from $\left(\Lambda_-,\cT_-\right)$ to $\left(\Lambda_+,\cT_+\right)$ and their decorations can be identified via the underlying isotopy (up to equivalence of decorations). Denote both decorations by $\cP$. Then the Hamiltonian isotopy induces a dga homotopy $\Phi_L^*\cong \Phi_{L'}^*$.
\end{thm}

Let $(L, \cP)$ be a decorated 
exact Lagrangian filling of $\Lambda$. Dualizing the homomorphism
\[ 
\cA\left(\Lambda\right)\xrightarrow{\phi_+^*} \cA\left(\Lambda, \cP\right) \xrightarrow{\Phi_L^*}
\cA\left(\emptyset, \cP\right),
\]
we obtain a morphism of algebraic varieties
\begin{equation}
\label{alphaLP}
\alpha_{L, \mathcal{P}}=\phi_+\circ \Phi_L:~ \Aug\left(\emptyset,\cP\right)\xrightarrow{\Phi_L} \Aug\left(\Lambda,\cP\right)\xrightarrow{\phi_+}\Aug\left(\Lambda\right).
\end{equation}
A decoration $\mathcal{P}$ is {\it sufficient} if its complement  $L-\mathcal{P}$ is a disjoint union of simply-connected regions. If each component of $\Lambda$ contains at least one marked point, then such a sufficient decoration $\mathcal{P}$ exists. We study the image of $\alpha_{L, \mathcal{P}}$ for sufficient $\mathcal{P}$.

\begin{defn} \label{decorated.character} Let $L$ be a compact oriented surface with boundary. Let  $\mathcal{T}$ be a collection of marked points on the boundary of $L$. We assume that each boundary component contains at least one marked point.
The decorated character variety $\mathscr{A}(L, \mathcal{T})$ parametrizes the data $(\mathcal{L}, \{v_i\})$, where
\begin{itemize}
    \item $\mathcal{L}$ is a line bundle over $L$ with flat connection,
    \item for every boundary interval $i$ in $\partial L-\mathcal{T}$, the data $v_i$ is a nontrivial flat section of $\mathcal{L}$ over $i$.
\end{itemize}
\end{defn}
The space $\mathscr{A}(L, \mathcal{T})$ is a special case of the moduli space of decorated G-local systems introduced by Fock and Goncharov in \cite{FGteich}. Let $[c]$ be the homotopy class of a oriented curve connecting two boundary components $j$ and $k$.
For each $(\mathcal{L}, \{v_i\})\in \mathscr{A}(L, \mathcal{T})$,  we parallel transport the section $v_j$ along $[c]$, obtaining a flat section $v_i'$ over $k$. It gives rise to a function $g_{[c]}$ of $\mathscr{A}(L, \mathcal{T})$ such that 
\[g_{[c]}=\frac{v_k}{v_j'}.\]

\medskip

Now suppose $L$ is an exact Lagrangian filling of a Legendrian link $\Lambda$. Let $\mathcal{T}$ be a collection of marked points on $\Lambda$, with each component of $\Lambda$ contains at least one marked point. 

\begin{lem} 
\label{8.28.11.27}
There is natural morphism
\[
\pi: {\rm Aug}(\emptyset, \mathcal{P}) \longrightarrow \mathscr{A}(L, \mathcal{T}).
\]
If $\mathcal{P}$ is sufficient, then $\pi$ is surjective.
\end{lem}
\begin{proof}
Following the proof of Lemma \ref{triv=augp}, ${\rm Aug}(\emptyset, \mathcal{P})$ is naturally isomorphic to the moduli space of trivilizations of ${\rm GL}_1$-local systems on $L$, with a choice of a vector on each connected region of $L-\mathcal{P}$. By forgetting vectors assigned to regions that are not connected to boundaries of $L$, we obtain a morphism $\pi$ from ${\rm Aug}(\emptyset, \mathcal{P})$ to $\mathscr{A}(L, \mathcal{T})$. By the definition of sufficiency of $\mathcal{P}$, the map $\pi$ is surjective. 
\end{proof}

\begin{lem} 
\label{8.28.11.28}
For every exact Lagrangian filling $L$ of $\Lambda$, there is a natural morphism \[\kappa_L: ~\mathscr{A}(L, \mathcal{T}) \longrightarrow {\rm Aug}(\Lambda).\]
The composition $\kappa_L\circ \pi$ coincides with the morphism $\alpha_{L,\mathcal{P}}$ in \eqref{alphaLP}.
\end{lem}
\begin{proof}
Let $a$ be a Reeb chord of $\Lambda$ with degree 0. Recall the moduli space $\mathcal{M}(a)$ of $J$-holomorphic disks such that the boundary of each disk $u\in \mathcal{M}(a)$ is $a$ together with a path $c(u)$ in $L$.  We set 

\begin{equation}
\label{ep1}
\epsilon(a)= \sum_{u\in \mathcal{M}(a)} g_{[c(u)]}.
\end{equation}
For every marked point $t\in \mathcal{T}$, let $c(t)$ be the unique path connecting the neighbored boundary intervals of $t$ such that $c(t)$ can be retracted to $t$. We set
\begin{equation}
\label{ep2}
\epsilon(t)=g_{[c(t)]}.
\end{equation}
 By definition, for each $(L, \{v_i\}) \in \mathscr{A}(L, \mathcal{T})$, its image under  \eqref{ep1} and \eqref{ep2} is an augmentation of $\Lambda$, which gives rise to the morphism $\kappa_L$. 
The identity $\alpha_L\circ \pi = \alpha_{L, \mathcal{P}}$ follows by a comparison of definitions. 
\end{proof}

As a consequence of Lemma \ref{8.28.11.27} and \ref{8.28.11.28}, when $\mathcal{P}$ is sufficient, the image of $\alpha_{L, \mathcal{P}}$ coincides with the image of $\kappa_L$. In this case, the image, denoted by ${\bf Im}(\alpha_L)$, is independent of the sufficient decoration $\mathcal{P}$ chosen. Combining with Theorem \ref{thm 3.8}, we get

\begin{cor}\label{sameimage}
Suppose the exact Lagrangian fillings $L$ and $L'$ are Hamiltonian isotopic, then 
${\bf Im}(\alpha_L)={\bf Im}(\alpha_{L'})$.
\end{cor}

\smallskip

\subsection{EHK Functorial Morphisms for Admissible Cobordisms}
In this section, we present explicit computations of the EHK morphisms associated with four basic types of exact Lagrangian cobordisms, i.e., saddle cobordisms, cyclic rotations, braid moves, and minimum cobordisms, for positive braid Legendrian links. Compositions of such cobordisms are called {\it admissible cobordisms}. An admissible cobordism from $\emptyset$ to $\Lambda$ is called an {\it admissible filling}.

\medskip
\noindent \textbf{(I) Saddle Cobordism.} Let $\left(\Lambda_+,\cT_+\right)$ be a decorated Legendrian link. Let $b$ be  a contractible Reeb chord as in Definition 6.12 of \cite{EHK12}. As in Figure \ref{saddle figure}, we contract $b$ via a saddle cobordism $S$, obtaining a new Legendrian link $\Lambda_-$. 
The holomorphic disk, represented by the gradient flow tree traced out by the contraction of $b$, is called the \emph{basic disk} associated with $b$ and denoted by $u_b$.
\begin{figure}[H]
\begin{tikzpicture}[baseline=45]
\draw (-2,-0.5) to [out=0,in=180] (-1,0) to [out=180,in=0] (-2,0.5);
\draw (2,-0.5) to [out=180,in=0] (1,0) to [out=0,in=180] (2,0.5);
\draw (-1,0) to [out=90,in=180] (0,0.5) to [out=0,in=90] (1,0);
\draw (-2,0.5) to [out=90,in=180] (0,1.5) to [out=0,in=90] (2,0.5);
\draw (-2,-0.5) to [out=90,in=-135] (-1.5,0.25);
\draw (2,-0.5) to [out=90,in=-45] (1.5,0.25);
\draw [dashed] (-1.5,0.25) to [out=45, in=180] (0,0.75) to [out=0,in=135] (1.5,0.25);
\path [fill=lightgray] (0,1.5) to [out=-120,in=65] (-0.25,0.5) to (0,0.75);
\draw [->] (0,0.75) -- node [right] {$b$} (0,1.5);
\draw [<-, thick] (1,0) to [out=90,in=0] (0,0.5) to [out=180,in=90] (-1,0) ;
\node at (0.5,0.25) [] {$p$};
\end{tikzpicture}
\quad \quad \quad \quad 
\begin{tikzpicture}[baseline=0, scale =0.8]
\draw [decoration={markings,mark=at position 0.7 with {\arrow{<}}},postaction={decorate}] (-1,0.5) to [out=-20,in=180] (0,0.25) to [out=0,in=-160] (1,0.5);
\draw [decoration={markings,mark=at position 0.7 with {\arrow{<}}},postaction={decorate}] (1,-0.5) to [out=160,in=0] (0,-0.25) to [out=180,in=20] (-1,-0.5);
\draw [<-] (-1,-1.5) to [out=0,in=180] (-0.25,-2) to [out=180,in=0] (-1,-2.5);
\draw [<-] (1,-2.5) to [out=180,in=0] (0.25,-2) to [out=0,in=180] (1,-1.5);
\node at (0,-3.5) [] {front projection};
\draw [->] (0,-0.25) -- node [right] {$b$} (0,0.25);
\node at (-0.25,-2) [] {$\ast$};
\node at (0.25,-2) [] {$\ast$};
\end{tikzpicture}
\quad \quad 
\begin{tikzpicture}[baseline=0, scale=0.8]
\draw [->] (1,0.5) -- (-1,-0.5);
\draw (-1,0.5) -- (-0.1,0.05);
\draw [->] (0.1,-0.05) -- (1,-0.5);
\draw [decoration={markings,mark=at position 0.75 with {\arrow{>}}},postaction={decorate}] (-1,-1.5) arc (90:-90:0.5);
\draw [decoration={markings,mark=at position 0.75 with {\arrow{>}}},postaction={decorate}] (1,-1.5) arc (90:270:0.5);
\node at (0,0) [below] {$b$};
\node at (0,-3.5) [] {Lagrangian projection};
\node at (-0.5,-2) [] {$\ast$};
\node at (0.5,-2) [] {$\ast$};
\end{tikzpicture}
\caption{Saddle Cobordism}\label{saddle figure}
\end{figure}

We decorate a saddle cobordism $S$ as follows. First, each marked point on $\Lambda_+$ traces out an oriented path that goes from $\Lambda_+$ to $\Lambda_-$. Second, the unstable manifold of the saddle defines a new path $p$, both of whose endpoints are on $\Lambda_-$, and we orient it so that the homological intersection of $\partial u_b$ and $p$ is 1. The induced decoration on $\Lambda_-$ is $\cT_-:=\cT_+\sqcup\left\{p^{\pm 1}\right\}$.

When the contractible Reeb chord $b$ is simple (\cite[Definition 6.15]{EHK12}), the recipe for the dga homomorphism $\Phi_S^*$ stated in \textit{loc. cit.} can be modified slightly to incorporate the enhancement of coefficients. For any Reeb chord $a\neq b$, let $\cM\left(a,b;c_1,\dots, c_n\right)$ be the moduli space of holomorphic disks that map into $\left(\mathbb{R}^4_{txyz},\mathbb{R}_t\times \Lambda_+\right)$, with one positive puncture at each of $a$ and $b$, and one negative puncture at each of the $c_i$'s. We define
\[
\left(\Phi_S^*\right)^0(d)=\left\{\begin{array}{ll}
    d & \text{if $d\neq b$}, \\
    p & \text{if $d=b$},
\end{array}\right. 
\]
\[
\left(\Phi_S^*\right)^1(d)=\left\{\begin{array}{cl}
    \displaystyle\sum_{\substack{c_1,\dots, c_n\in \cR_+ \\ \dim \cM\left(d,b;c_1,\dots, c_n\right)=1}}\sum_{u\in \cM\left(d,b;c_1,\dots,c_n\right){{/\bbR}}}  w(u){\Big|}_{b=p^{-1}}  & \text{if $d\neq b$},\\
    0 & \text{if $d=b$},
\end{array}\right.
\]
where $w(u)$ is defined in the same way as \eqref{w(u)}. The homomorphism $\Phi_S^*$ is 
\begin{equation}\label{Phi_S}
\Phi_S^*=\left(\Phi_S^*\right)^0+\left(\Phi_S^*\right)^1.
\end{equation}

\begin{rmk} 
In $(\Phi_S^{*})^0$, $b$ is mapped to $p$, whereas in $(\Phi_S^{*})^1$, $b$ is substituted by $p^{-1}$. Their difference can be understood via holomorphic disk degeneration. Note the cusp edge in Figure \ref{saddle figure}. The term in $(\Phi_S^{*})^1$ comes from a negative degeneration of an \emph{end} at the cusp edge, and the term in $(\Phi_S^{*})^1$ comes from a positive degeneration of a \emph{switch} at the cusp edge.\footnote{These singular points were first introduced in \cite{ekholm2007}. Pictures are available in Figure 3 of \cite{EENS}.} Their contributions with respect to the marked curve are reciprocal.
\end{rmk}

Degree 0 Reeb chords in a positive braid Legendrian link are contractible but not necessarily simple in general. For contractible Reeb chords that are not simple, Ekholm, Honda, and K\'{a}lm\'{a}n stated that these cases can be reduced to the simple cases by implementing a collection of ``dippings'' \cite[Figure 17]{EHK12}, a notion introduced in \cite{F} and also appeared in \cite{Sab, FR}. 

For positive braid Legendrian links, it turns out that two dippings will suffice. Below is a depiction of local moves on the Lagrangian projection for a saddle cobordism that pinches a degree 0 Reeb chord $b_k$ of a positive braid Legendrian link $\Lambda_\beta$.
\begin{figure}[H]
\begin{tikzpicture}[baseline=-3, scale = 0.8]
\draw (-1,-0.5) -- (-0.75,-0.5) to [out=0,in=-135] (-0.1,-0.1);
\draw (0.1,0.1) to [out=45,in=180] (0.75,0.5) -- (1,0.5);
\draw (-1,0.5) -- (-0.75,0.5) to [out=0,in=180] (0.75,-0.5) -- (1,-0.5);
\node at (0,0) [above] {$b_k$};
\node at (0,-1.5) [] {$\Lambda_+$};
\end{tikzpicture} 
\, $\xleftarrow{D}$ \,
\begin{tikzpicture}[baseline=-3, scale = 0.8]
\draw (-2.5,-0.5) -- (-1.75,-0.5) --(-1.75,1) -- (-1,1) -- (-1,-0.5) -- (-0.75,-0.5) to [out=0,in=-135] (-0.1,-0.1);
\draw (-2.5,0.5) -- (-1.85,0.5);
\draw (-1.65,0.5) -- (-1.1,0.5);
\draw (-0.9,0.5) -- (-0.75,0.5) to [out=0,in=180] (0.75,-0.5) -- (0.9,-0.5);
\draw (1.1, -0.5) -- (1.65,-0.5);
\draw (1.85,-0.5) -- (2.5,-0.5);
\draw (0.1,0.1) to [out=45,in=180] (0.75,0.5)--(1,0.5) -- (1,-1) -- (1.75,-1) -- (1.75,0.5) -- (2.5,0.5);
\node at (0,0) [above] {$b_k$};
\node at (-1,0.5) [above right] {$x_L$};
\node at (-1.75,0.5) [above left] {$y_L$};
\node at (1,-0.5) [below left] {$x_R$};
\node at (1.75,-0.5) [below right] {$y_R$};
\node at (0,-1.5) [] {$\Lambda'$};
\end{tikzpicture}  
\, $\xleftarrow{S}$ \,
\begin{tikzpicture}[baseline=-3, scale = 0.8]
\draw (-2.5,-0.5) -- (-1.75,-0.5) --(-1.75,1) -- (-1,1) -- (-1,-0.5) -- (0.9,-0.5);
\draw (-2.5,0.5) -- (-1.85,0.5);
\draw (-1.65,0.5) -- (-1.1,0.5);
\draw (1.1, -0.5) -- (1.65,-0.5);
\draw (1.85,-0.5) -- (2.5,-0.5);
\draw (-0.9,0.5) --(1,0.5) -- (1,-1) -- (1.75,-1) -- (1.75,0.5) -- (2.5,0.5);
\node at (-1,0.5) [above right] {$x_L$};
\node at (-1.75,0.5) [above left] {$y_L$};
\node at (1,-0.5) [below left] {$x_R$};
\node at (1.75,-0.5) [below right] {$y_R$};
\node at (0,0.5) [] {$\ast$};
\node at (0,-0.5) [] {$\ast$};
\node at (0,0.5) [above right] {$p_k$};
\node at (0,-0.5) [below left] {$p_k^{-1}$};
\node at (0,-1.5) [] {$\Lambda''$};
\end{tikzpicture}
\,$\xleftarrow{D^{-1}}$\,
\begin{tikzpicture}[baseline=-3, scale = 0.8]
\draw (-0.75,0.5) -- (0.75,0.5);
\draw (-0.75,-0.5) -- (0.75,-0.5);
\node at (0,0.5) [] {$\ast$};
\node at (0,-0.5) [] {$\ast$};
\node at (0,0.5) [above] {$p_k$};
\node at (0,-0.5) [below] {$p_k^{-1}$};
\node at (0,-1.5) [] {$\Lambda_-$};
\end{tikzpicture}
\caption{Dipping and Saddle Cobordism in the Lagrangian Projection}
\end{figure}
Among the three steps, $D$ and $D^{-1}$ are compositions of Legendrian Reidemeister II moves, and hence we can compute $\Phi_{D}^*$ and $\Phi_{D^{-1}}^*$ by following  \cite[\S 8.4]{Che02};  $S$ is a simple saddle cobordism, allowing us to employ \eqref{Phi_S} to compute $\Phi_S^*$. 

\begin{prop}\label{pinch theorem} Let $S_k$ be the saddle cobordism contracting the Reeb chord $b_k$ of $\Lambda_\beta$. The functorial dga homomorphism $\Phi_{S_k}^*:\cA\left(\Lambda_+,\cP\right)\rightarrow \cA\left(\Lambda_-,\cP\right)$ maps the degree 0 Reeb chords as follows:
\begin{equation}
\label{formula for Phi}
\Phi_{S_k}^*\left(b_s\right)=\left\{\begin{array}{ll}
    \displaystyle b_s+\sum_{\partial b_s=\sum uy_Lv}up_k^{-1}\Phi_{S_k}^*(v) & \text{if $s<k$},  \\
    p_k & \text{if $s=k$},\\
    \displaystyle b_s+\sum_{\partial b_s=\sum uy_Rv} \Phi_{S_k}^*(u)p_k^{-1}v & \text{if $s>k$}.
\end{array}\right.
\end{equation}
Here the summation index $\partial b_s=\sum uy_Lv$ and $\partial b_s=\sum uy_Rv$ are computed on the Lagrangian projection after the dipping $D$.
\end{prop}
\begin{proof} Chekanov \cite{Che02} constructed a pair of tame dga isomorphisms $\psi_+:\cA\left(\Lambda'\right)\rightarrow S\cA\left(\Lambda_+\right)$ and $\psi_-:\cA\left(\Lambda''\right)\rightarrow S\cA\left(\Lambda_-\right)$, where $S$ denotes a stablization of the dga. By \cite[Lemma 6.7, 6.8, Remark 6.9]{EHK12}, we know that the dga homomorphisms $\Phi_D^*$ and $\Phi_{D^{-1}}^*$ are given by
\[
\Phi_D^*:\cA\left(\Lambda_+\right) \hookrightarrow S\cA\left(\Lambda_+\right)\xrightarrow{\psi_+^{-1}} \cA\left(\Lambda'\right),
\]
\[
\Phi_{D^{-1}}^*:\cA\left(\Lambda''\right)\xrightarrow{\psi_-} S\cA\left(\Lambda_-\right)\twoheadrightarrow\cA\left(\Lambda_-\right).
\]
By following Chekanov's recipe, we see that for any degree 0 Reeb chord $b_s$ of $\Lambda_\beta$,
\[
\Phi_D^*\left(b_s\right)=\left\{\begin{array}{ll}
    \displaystyle b_s+\sum_{\partial b_s=\sum uy_Lv} ux_L\Phi_D^*(v)& \text{if $s<k$}, \\
    b_k & \text{if $s=k$}, \\
    \displaystyle b_s+\sum_{\partial b_s=\sum uy_Rv} \Phi_D^*(u)x_Rv & \text{if $s>k$}.
\end{array}\right.
\]
and $\Phi_{D^{-1}}^*$ annihilates all occurrences of $x_L$ and $x_R$.

On the other hand, between the dipping and undipping cobordisms, we have a simple saddle cobordism $S$, and by \eqref{Phi_S} we see that
\[
\Phi_S^*\left(b_k\right)=p_k, \quad \quad \Phi_S^*\left(x_L\right)=x_L+p_k^{-1}, \quad \quad \Phi_S^*\left(x_R\right)=x_R+p_k^{-1}.
\]
By composing $\Phi_{D^{-1}}^*\circ \Phi_S^*\circ \Phi_D^*$, we get the formula stated in the proposition.
\end{proof}

The recursive nature of Formula \eqref{formula for Phi} suggests an algorithm, termed \emph{matrix scanning}, to compute $\Phi_{S_k}^*$ for degree 0 Reeb chords of $\Lambda_\beta$. This algorithm starts at the $k$th crossing and scans the left and the right portions of the braid using two family of matrices, which keep track of all possible incomplete disks sandwiched between levels.

Let us describe in details the family of matrices $\left\{U^{(s)}\right\}_{k+1\leq s \leq l}$, which we use to scan the braid word $s_{i_{k+1}}s_{i_{k+2}}\cdots s_{i_l}$. Each $U^{(s)}$ is an $n\times n$ upper triangular matrix, and he $(i,j)$-entry of $U^{(s)}$ counts partial disks between strands $i<j$ right before scanning through the crossing $i_s$. Following this idea, we see that the initial matrix $U^{(k+1)}$ must have all entries 0 except the $\left(i_k,i_k+1\right)$-entry, which is $p_k^{-1}$.

Inductively for $s>k$, we scan through the crossing $i_s$ and perform two actions. First, we record
\begin{equation}\label{right-scanning}
\Phi_{S_k}^*\left(b_s\right)=b_s+U^{(s)}_{i_s,i_s+1}.
\end{equation}
Second, we define $U^{(s+1)}$ in terms of $U^{(s)}$; note that these two matrices differ only at entries whose rows or columns are equal to $i_s$ or $i_s+1$:
\[
\begin{array}{ll}
U^{(s+1)}_{i,i_s}=U_{i,i_s+1}^{(s)}+U_{i,i_s}^{(s)}b_s   \quad \quad & U_{i,i_s+1}^{(s+1)}=U_{i,i_s}^{(s)} \\[4pt]
U_{i_s,j}^{(s+1)}=U_{i_s+1,j}^{(s)} & U_{i_s+1,j}^{(s+1)}=U_{i_s,j}^{(s)}+\Phi_{S_k}^*\left(b_s\right)U_{i_s+1,j}^{(s)}.
\end{array} \quad  
\begin{tikzpicture}[baseline=15, scale=.8]
\draw (-1,0) node [left] {$j$th} -- (1,0);
\draw (-1,0.5) node [left] {$\left(i_s+1\right)$st} to [out=0,in=180] (1,1);
\draw (-1,1) node [left] {$i_s$th} to [out=0,in=180] (1,0.5);
\draw (-1,1.5) node [left] {$i$th} -- (1,1.5);
\node at (0,0.75) [below] {$b_s$};
\end{tikzpicture}
\]

To describe this transformation more compactly, we introduce the following \emph{triangular truncations} for  matrices:
\begin{equation}\label{triangular truncation}
M_{ij}^+:=\left\{\begin{array}{ll} M_{ij} & \text{if $i<j$,}\\
0 & \text{otherwise},
\end{array}\right.\quad \quad 
M_{ij}^-:=\left\{\begin{array}{ll} M_{ij} & \text{if $i>j$},\\
0 & \text{otherwise}.
\end{array}\right.
\end{equation}
Recall the matrix $Z_{i_s}$ from \eqref{SwitchMat}. Then $U^{(s+1)}$ is defined in terms of $U^{(s)}$ as
\begin{equation}\label{right-scanning2}
U^{(s+1)} = \left(Z_{i_s}\left(\Phi_{S_k}^*\left(b_s\right)\right)^{-1}\cdot U^{(s)}\cdot Z_{i_s}\left(b_s\right)\right)^+.
\end{equation}

The left-scanning family of matrices $\left\{L^{(s)}\right\}_{1\leq s \leq k-1}$ works similarly. Each matrix $L^{(s)}$ is an $n\times n$ lower triangular matrix and the $(i,j)$-entry of $L^{(s)}$ counts partial disks between strands $i>j$ right before scanning through the crossing $i_s$. The initial matrix $L^{(k-1)}$ has all entries $0$ except that its  $\left(i_k+1,i_k\right)$-entry is $p_k^{-1}$. Inductively for $s<k$, we perform the following two actions when scanning through a crossing $i_s$. First, we record 
\begin{equation}\label{left-scanning}
\Phi_{S_k}^*\left(b_s\right)=b_s+L_{i_s+1,i_s}^{(s)}.
\end{equation}
Second, we define $L^{(s-1)}$ in terms of $L^{(s)}$ according to
\begin{equation}\label{left-scanning2}
L^{(s-1)} := \left(Z_{i_s}\left(b_s\right)\cdot L^{(s)}\cdot Z_{i_s}\left(\Phi_{S_k}^*\left(b_s\right)\right)^{-1}\right)^-.
\end{equation}

Note that \eqref{right-scanning} and \eqref{left-scanning}, together with $\Phi_{S_k}^*\left(b_k\right)=p_k$, completely describe the image of all degree 0 Reeb chords in $\Lambda_\beta$ under the functorial homomorphism $\Phi_{S_k}^*$.

\begin{rmk} For a degree $0$ Reeb chord in a positive braid Legendrian link, Proposition \ref{pinch theorem} yields an explicit description $\Phi^*_{S_k}\left(b_s\right):=\sum_{m=0}^\infty (\Phi_{S_k}^{*})^{m}\left(b_s\right)$, where for $m\geq 2$,
\begin{equation}
(\Phi_{S_k}^{*})^m\left(d\right):=
\left\{
\begin{array}{cc}
     \displaystyle \sum_{\substack{c_1,\dotsb,c_n \in \cR_+ \\ \dim \cM\left(d,b_k^m;c_1,\dots, c_n\right)=m-1 }}\sum_{u\in \cM\left(d,b_k^m;c_1,\dots, c_n\right)/\mathbb{R}^{m-1}}w(u){\Big|}_{b_k=p_k^{-1}}& \textrm{if } d\neq b_k, \\
     0 & \textrm{if } d=b_k,
\end{array}
\right.
\end{equation}
where $\cM\left(d,b_k^m;c_1,\dots, c_t\right)$ is the moduli space of immersed disks in $(\bbR^2_{xy},\pi_L(\Lambda_+))$ with positive quadrants at $b_s$ and $b_k$, and remaining negative quadrants, where a negative quadrant is allowed to be a $(-+-)$ triple quadrant.
\end{rmk}

\medskip 

\noindent \textbf{(II) Cyclic Rotation.} A cyclic rotation is a Legendrian isotopy from  $\Lambda_{\beta s_i}$ to $\Lambda_{s_i \beta}$, illustrated by the moves on the front projection of Legendrian links in Figure \ref{cyclic rotation figure}.
\begin{figure}[H]
\begin{tikzpicture}[baseline=17,scale =0.35]
    \draw (0,0) rectangle node [] {$\beta$} (-1.5,1.5);
    \draw  (0,0.25) to [out=0,in=180] (2.25,2) to [out=180,in=0] (0, 3.75) -- (-2.5,3.75) to [out=180,in=0] (-4.75,2) to [out=0,in=180] (-2.5,0.25) to [out=0,in=180] (-1.5, 0.75);
    \draw (0,0.75) to [out=0,in=180] (1.5,2) to [out=180,in=0] (0, 3.25) -- (-2.5,3.25) to [out=180,in=0] (-4,2) to [out=0,in=180] (-2.5,0.75) to [out=0,in=180] (-1.5, 0.25);
    \draw (0,1.25)  to [out=0,in=180] (0.75,2) to [out=180,in=0] (0, 2.75) -- (-2.5,2.75) to [out=180,in=0] (-3.25,2) to [out=0,in=180] (-2.5,1.25) -- (-1.5,1.25);
    \node at (-2, 0.5) [below] {$s_i$};
\end{tikzpicture} 
 $\xleftarrow{\alpha}$
\begin{tikzpicture}[baseline=17,scale =0.35]
    \draw (0,0) rectangle node [] {$\beta$} (-1.5,1.5);
    \draw  (0,0.25) to [out=0,in=180] (2.25,2) to [out=180,in=0] (0, 3.75) -- (-1.5,3.75) to [out=180,in=0]  (-2.5,3.25) to [out=180,in=0] (-4.75,1.65) to [out=0,in=180] (-2.5,0.25) to [out=0,in=180] (-1.5, 0.75);
    \draw (0,0.75) to [out=0,in=180] (1.5,2) to [out=180,in=0] (0, 3.25) -- (-1.5,3.25) to[out=180,in=0] (-2.5,3.75) to [out=180,in=0] (-4.75,2.35) to [out=0,in=180] (-2.5,0.75) to [out=0,in=180] (-1.5, 0.25);
    \draw (0,1.25)  to [out=0,in=180] (0.75,2) to [out=180,in=0] (0, 2.75) -- (-2.5,2.75) to [out=180,in=0] (-3.5,2) to [out=0,in=180] (-2.5,1.25) -- (-1.5,1.25);
\end{tikzpicture} 
$\xleftarrow{\beta}$
\begin{tikzpicture}[baseline=17,scale=0.35]
\draw (0,0) rectangle node [] {$\beta$} (1.5,1.5);
\draw  (0,0.25) to [out=180,in=0] (-2.25,2) to [out=0,in=180] (0, 3.75) -- (0.25,3.75) to[out=0,in=180] (1.25,3.25) --(1.5,3.25) to [out=0,in=180] (3,2) to [out=180,in=0] (1.5,0.75);
\draw (0,0.75)  to [out=180,in=0] (-1.5,2) to [out=0,in=180] (0,3.25) -- (0.25,3.25) to [out=0,in=180] (1.25,3.75) -- (1.5,3.75) to [out=0,in=180] (3.75,2) to [out=180,in=0] (1.5,0.25);
\draw (0,1.25) to [out=180,in=0] (-0.75,2) to[out=0,in=180] (0,2.75) -- (1.5,2.75) to [out=0,in=180] (2.25,2) to[out=180,in=0] (1.5,1.25);
\node at (0.75,3.5) [above] {$c$};
\end{tikzpicture} 
$\xleftarrow{\gamma}$ 
\begin{tikzpicture}[baseline=17,scale =0.35]
    \draw (0,0) rectangle node [] {$\beta$} (1.5,1.5);
    \draw  (0,0.25) to [out=180,in=0] (-2.25,2) to [out=0,in=180] (0, 3.75) -- (1.5,3.75) to [out=0,in=180]  (2.5,3.25) to [out=0,in=180] (4.75,1.65) to [out=180,in=0] (2.5,0.25) to [out=180,in=0] (1.5, 0.75);
    \draw (0,0.75) to [out=180,in=0] (-1.5,2) to [out=0,in=180] (0, 3.25) -- (1.5,3.25) to[out=0,in=180] (2.5,3.75) to [out=0,in=180] (4.75,2.35) to [out=180,in=0] (2.5,0.75) to [out=180,in=0] (1.5, 0.25);
    \draw (0,1.25)  to [out=180,in=0] (-0.75,2) to [out=0,in=180] (0, 2.75) -- (2.5,2.75) to [out=0,in=180] (3.5,2) to [out=180,in=0] (2.5,1.25) -- (1.5,1.25);
    \node at (2,3.5) [above] {$c$};
\end{tikzpicture} 
$\xleftarrow{\delta}$
\begin{tikzpicture}[baseline=17,scale =0.35]
    \draw (0,0) rectangle node [] {$\beta$} (1.5,1.5);
    \draw  (0,0.25) to [out=180,in=0] (-2.25,2) to [out=0,in=180] (0, 3.75) -- (2.5,3.75) to [out=0,in=180] (4.75,2) to [out=180,in=0] (2.5,0.25) to [out=180,in=0] (1.5, 0.75);
    \draw (0,0.75) to [out=180,in=0] (-1.5,2) to [out=0,in=180] (0, 3.25) -- (2.5,3.25) to [out=0,in=180] (4,2) to [out=180,in=0] (2.5,0.75) to [out=180,in=0] (1.5, 0.25);
    \draw (0,1.25)  to [out=180,in=0] (-0.75,2) to [out=0,in=180] (0, 2.75) -- (2.5,2.75) to [out=0,in=180] (3.25,2) to [out=180,in=0] (2.5,1.25) -- (1.5,1.25);
    \node at (2, 0.5) [below] {$s_i$};
\end{tikzpicture}
\caption{Cyclic Rotation} \label{cyclic rotation figure}
\end{figure}

We denote the  exact Lagrangian concordance corresponding to a cyclic rotation as 
\begin{equation}\label{def rho}
\rho:\Lambda_{\beta s_i}\rightarrow \Lambda_{s_i\beta}.
\end{equation}
We decorate $\rho$ with oriented marked curves tracing the marked points on either end.

\begin{prop}\label{3.15} Let $\beta=s_{i_1}\cdots s_{i_l}$ be a positive braid word. Set \[r(\beta):=s_{i_l}s_{i_1}\cdots s_{i_{l-1}}, \hskip 10mm l(\beta):= s_{i_2}\cdots s_{i_l}s_{i_1}.\]
Recall $M^{(k)}$  in Proposition \ref{partial a}.  The dga homomorphism $\Phi_\rho^*:\cA\left(\Lambda_{r(\beta)}\right)\rightarrow \cA\left(\Lambda_{\beta}\right)$ associated with the cyclic rotation
$\rho:\Lambda_{\beta}\rightarrow \Lambda_{r(\beta)}$ maps degree 0 Reeb chords as follows:
\[
\Phi_\rho^*\left(b_k\right)=b_{k-1} \quad \forall 1< k\leq l \quad \text{and}\quad \Phi_\rho^*\left(b_1\right)=M_{i_l+1,i_l}^{(i_l)}t_{i_l}.
\]
The dga homomorphism $\Phi_{\rho^{-1}}^*:\cA\left(\Lambda_{l(\beta)}\right)\rightarrow \cA\left(\Lambda_{\beta}\right)$ associated with the inverse cyclic rotation $\rho^{-1}:\Lambda_{\beta}\rightarrow \Lambda_{l(\beta)}$ maps degree 0 Reeb chords as follows:
\[
\Phi_{\rho^{-1}}^*\left(b_k\right)=b_{k+1} \quad \forall 1\leq k<l \quad \text{and} \quad \Phi_{\rho^{-1}}^*\left(b_l\right)=t_{i_1}M_{i_1,i_1+1}^{(i_1)}.
\]
\end{prop}
\begin{proof} We only prove the formula for $\Phi_\rho^*$. The proof  for $\Phi_{\rho^{-1}}^*$ is similar.

We break the cyclic rotation $\rho$ into steps according to Figure \ref{cyclic rotation figure} and use the bordered dga method \cite{Siv10} to compute the functorial homomorphism for each step. First, by considering the bordered dga on the complement of the right cusps region, we see that
\[
\Phi_\beta^*\circ \Phi_\alpha^*\left(b_k\right)=\left\{\begin{array}{ll}
    c &  \text{if $k=1$}, \\
    b_{k-1} & \text{otherwise},
\end{array}\right.
\]
where $c$ is depicted in Figure \ref{cyclic rotation figure}.

The Reidemeister II move is performed away from the crossing $c$ and the braid region. Therefore  $\Phi_\gamma^*(c)=c$ and $\Phi_\gamma^*\left(b_k\right)=b_k$. For the same reason,  $\Phi_\delta^*\left(b_k\right)=b_k$, which, combined with the formulas of $\Phi_\beta^*\circ \Phi_\alpha^*$ and $\Phi_\gamma^*$, implies that $\Phi_\rho^*\left(b_k\right)=b_{k-1}$ for $1<k\leq l$.

It remains to compute $\Phi_\delta^*(c)$. Define $i:=i_l$. Let us consider the bordered dga of the region on the right of the braid region (including the crossing $c$). The differentials of the degree 1 Reeb chords of the bordered dga before $\Phi_\delta^*$ are
\[
\begin{array}{l} \partial a_i=t_i^{-1}+x_{24}+x_{23}b_l, \\ [5pt]
\partial a_{i+1}=t_{i+1}^{-1}+x_{13}+cx_{23}, \\ [5pt]
\partial d = x_{14}+x_{13}b_l+cx_{24}+cx_{23}b_l.
\end{array}
\quad \quad \quad \quad \quad 
\begin{tikzpicture}[baseline=17, scale=0.7]
    \draw (0,0) node [left] {$4$} to [out=0,in=180] (1,0.5) to [out=0,in=180] (2,1.375) to [out=180,in=0] (1,2) to [out=180,in=0] (0,1.5) node [left] {$2$};
    \draw (0,0.5) node [left] {$3$} to [out=0,in=180] (1,0) to [out=0,in=180] (2,0.625) to [out=180,in=0] (1,1.5) to [out=180,in=0] (0,2) node [left] {$1$};
    \node at (0.5,0.25) [below] {$b_l$};
    \node at (0.5,1.75) [above] {$c$};
    \node at (1.5,1) [right] {$d$};
    \node at (2,1.375) [right] {$a_i$};
    \node at (2,0.625) [right] {$a_{i+1}$};
\end{tikzpicture}
\]
where $t_i$ denotes the marked point near the Reeb chord $a_i$.
On the other hand, we know that the differentials of the degree 1 Reeb chords of the bordered dga after $\Phi_\delta^*$ are
\[
\partial a_i=t_i^{-1}+x_{24}+x_{23}b_l, \quad \quad \quad \partial a_{i+1} =t_{i+1}^{-1}+x_{13}+x_{14}t_ix_{23}+x_{13}b_lt_ix_{23}+x_{12}a_it_ix_{23}. 
\]
By comparison, we see that $\Phi_\delta^*\left(a_i\right)=a_i$, $\Phi_\delta^*\left(a_{i+1}\right)=a_{i+1}$, and most importantly,
\[
\Phi_\delta^*(c)=x_{14}t_i+x_{13}b_lt_i+x_{12}a_it_i.
\]
Now if we include the bordered dga of the remaining part of $\Lambda_{\left(i_1,\dots, i_l\right)}$, we see that $x_{12}=0$ and $x_{14}+x_{13}b_l=M_{i+1,i}^{(i)}$. Therefore we  conclude that
\[
\Phi_\rho^*\left(b_1\right)=\Phi_\delta^*\circ \Phi_\gamma^*\circ \Phi_\beta^*\circ \Phi_\alpha^*\left(b_1\right)=\Phi_\delta^*(c)=M_{i+1,i}^{(i)}t_i. \qedhere
\]
\end{proof}

\medskip 

\noindent\textbf{(III) Braid Move.} A braid move $B$ is a Reidemeister III move within the braid region, which naturally gives rise to an invertible exact Lagrangian concordance $B$. We decorate $B$ with oriented marked curves that are the traces of the marked points on either end.

Consider the braid move $B:\Lambda_{(\dots s_is_js_i\dots)}\rightarrow \Lambda_{(\dots s_js_is_j \dots)}$, where $|i-j|=1$. Let $b_1,b_2,b_3$ be the Reeb chords corresponding to the crossings involved in $B$. By  \cite[\S 8.2, 8.3]{Che02}, 
\begin{equation}\label{braid move}
\Phi_B^*\left(b_1\right)=b_3, \quad \quad \Phi_B^*\left(b_2\right)=b_2+b_3b_1, \quad \quad \Phi_B^*\left(b_3\right)=b_1.
\end{equation}
The move $B$ is local. Therefore the rest Reeb chords are invariant under $\Phi_B^*$.

\medskip

\noindent\textbf{(IV) Minimum Cobordism.} Let $O$ denote the Legendrian unknot whose Thurston-Bennequin number is $-1$. Without loss of generality, we assume that $O$ has only one Reeb chord $a$. Then $|a|=1$ in the CE dga $\cA(O)$. By \cite{EP}, the unknot $O$ has a unique exact Lagrangian filling $M$, called the \emph{minimum cobordism}, which topologically is a hemisphere capping off $O$. 
\begin{figure}[H]
\begin{tikzpicture}[scale=0.6]
\draw (135:1.414) arc (45:315:1.414);
\draw (-45:1.414) arc (-135:135:1.414);
\draw (1,-1) -- (-1,1);
\draw [decoration={markings,mark=at position 0.75 with {\arrow{>}}},postaction={decorate}] (-1,-1) -- (1,1);
\node at (0,0) [below]{$a$};
\tikzset{shift={(-2,0)}}
\foreach \i in {1,...,5}
    {
    \node at (45*\i:1.414) [] {$\ast$};
    \node at (45*\i:1.75) [] {$t_\i$};
    }
    \node at (-45:1.414) [] {$\ast$};
    \node at (-45:1.75) [] {$t_k$};
    \node at (-90:1.75) [] {$\cdots$};
\tikzset{shift={(4,0)}}
\foreach \i in {1,...,5}
    {
    \node at (180-45*\i:1.414) [] {$\ast$};
    \node at (180-45*\i:1.85) [] {$t_{k+\i}$};
    }
    \node at (-135:1.414) [] {$\ast$};
    \node at (-135:1.75) [] {$t_m$};
    \node at (-90:1.75) [] {$\cdots$};
\end{tikzpicture}\hspace{1cm}
\begin{tikzpicture}[scale=.8]
\draw [decoration={markings,mark=at position 0.25 with {\arrow{>}}},postaction={decorate}] (0,0) ellipse (2 and 0.5);
\node (4) at (1,0.414) [] {$\ast$};
\node (5) at (-1,0.414) [] {$\ast$};
\node (3) at (1,-0.414) [] {$\ast$};
\node (1) at (-1,-0.414) [] {$\ast$};
\node at (2,0) [] {$\ast$};
\node at (-2,0) [] {$\ast$};
\node (2) at (0,-0.5) [] {$\ast$};
\node at (-1,0.414) [above] {$t_m$};
\node at (-2,0) [left] {$t_1$};
\node at (-1,-0.414) [above left] {$t_2$};
\node at (0,-0.5) [above] {$t_3$};
\node at (1,-0.414) [above right] {$t_4$};
\node at (2,0) [right] {$t_5$};
\node at (1,0.414) [above] {$t_6$};
\node at (0,0.5) [above] {$\cdots$};
\node at (0,-2) [below] {$\tau$};
\draw [decoration={markings,mark=at position 0.5 with {\arrow{>}}},postaction={decorate}] (-2,0) arc (-180:-90:2);
\draw [decoration={markings,mark=at position 0.5 with {\arrow{>}}},postaction={decorate}] (2,0) arc (0:-90:2);
\draw [decoration={markings,mark=at position 0.5 with {\arrow{>}}},postaction={decorate}] (1) to [out=-90,in=150] (0,-2);
\draw [decoration={markings,mark=at position 0.5 with {\arrow{>}}},postaction={decorate}] (3) to [out=-90,in=30] (0,-2);
\draw [decoration={markings,mark=at position 0.5 with {\arrow{>}}},postaction={decorate}] (2) to (0,-2);
\draw [dashed, decoration={markings,mark=at position 0.5 with {\arrow{>}}},postaction={decorate}] (5) to [out=-85,in=135] (0,-2);
\draw [dashed, decoration={markings,mark=at position 0.5 with {\arrow{>}}},postaction={decorate}] (4) to [out=-95,in=45] (0,-2);
\end{tikzpicture}
    \caption{Legendrian unkont and its minimum cobordism}\label{minimum cobordism picture}
\end{figure}

Now suppose $O$ is decorated with $m$ marked points. 
Let $\cP$ be the decoration of $M$ with oriented marked curves that flow from the marked points on $O$ to the unique $t$-minimum $\tau$ on $M$. Abusing notation, let us label the the marked curves with the same symbol $t_i$ as the marked points they originate from. 
Then \[\cA\left(\emptyset, \cP\right)=\frac{\mathbb{Z}_2\left\langle t_1^{\pm 1},t_2^{\pm 1},\dots, t_m^{\pm 1}\right\rangle}{t_1t_2\cdots t_m = 1}.\] 
The functorial homomorphism $\Phi_M^*: \cA\left(O,\cP\right)\rightarrow \cA\left(\emptyset,\cP\right)$ maps
\begin{equation}\label{minimum cobordism}
\Phi_M^*\left(t_i\right)=t_i \quad \text{and} \quad \Phi_M^*(a)=0.
\end{equation}

\smallskip

\subsection{Cluster Charts from Admissible Fillings}\label{Sec 4.1} By \eqref{alphaLP}, every decorated admissible filling $(L,\cP)$ of $\Lambda_\beta$ gives rise to a morphism of algebraic varieties
\[
\alpha_L=\phi_+\circ \Phi_L:~\Aug\left(\emptyset,\cP\right)\xrightarrow{\Phi_L} \Aug\left(\Lambda_\beta,\cP\right)\xrightarrow{\phi_+}\Aug\left(\Lambda_\beta\right).
\]
In this section, we show that the morphism $\alpha_L$ is an open embedding of an algebraic torus and its image is a cluster chart of  $\Aug\left(\Lambda_\beta\right)$.

Among the basic exact Lagrangian cobordisms defining admissible cobordisms, the decorated saddle cobordism $S_k$ is the only one that creates new oriented marked curves (and hence marked points) in the decoration. Thus, on any admissible filling $L$, we have exactly $l+n$ many oriented marked curves in the decoration $\cP$ on $L$, where $l$ is the length of the braid word $\beta$ and $n$ is the number of strands. There are $n$ many $t$-minima on $L$, one for each strand in $\beta$, and hence we have $n$ relations \eqref{monodromy condition} among the formal variables associated with the marked curves. Moreover, since the original cuspidal marked points $t_i$ of $\Lambda_\beta$ end at distinct $t$-minima on $L$, we can use these $n$ relations to eliminate $t_i$, leaving the formal variables $p_i$ formal variables. This proves the following lemma.

\begin{lem} For any admissible filling $L$ of  $\Lambda_\beta$ with decoration $\cP$, $\Aug\left(\emptyset ,\cP\right)\cong \left(\mathbb{F}^\times\right)^l$.
\end{lem}

For any admissible filling $L$ of $\Lambda_\beta$ with decoration $\mathcal{P}$, every component of the complement $L-\mathcal{P}$ is simply connected. Thus, we can think of the numerical values of the formal variables $p_i$ and $t_i$ as a recording the transition functions of a trivialization of some rank 1 local system on $L$ (in the normal direction determined by the orientation of the marked curve). From this perspective, the condition \eqref{monodromy condition} at $t$-minima can be viewed as a compatibility condition for the transition functions. As a consequence, we get the following Lemma. 

\begin{lem} 
\label{triv=augp}
For any admissible filling $L$ of $\Lambda_\beta$ with decoration $\cP$, there is a natural isomorphism $\Aug\left(\emptyset, \cP\right)\cong \Triv_1(L,\cP)$, where $\Triv_1(L,\cP)$ denotes the moduli space of trivializations of rank 1 local systems on $L$ with respect to the family of oriented marked curves $\cP$.
\end{lem}

Next, using the isomorphism between $\Aug\left(\Lambda_\beta\right)$ and $\conf^e_\beta\left(\mathcal{C}\right)$ in Theorem \ref{3.6}, we prove that, for the first three types of the basic admissible cobordisms, the corresponding functorial morphism of augmentation varieties is intertwined with a certain quasi-cluster morphism between double Bott-Samelson cells. For the minimum cobordism, we show that the corresponding functorial morphism is an isomorphism of algebraic tori.

\medskip

\noindent\textbf{(I) Saddle Cobordism.} Let $\beta=s_{i_1}\cdots s_{i_l}$ be a positive braid word. Consider the decorated saddle cobordism $\left(S_k,\cP\right)$ that resolves the crossing $i_k$ into a new pair of marked points $p_k^{\pm 1}$. 
Let $\Lambda_-$ denote the obtained  decorated Legendrian link. The underlying undecorated Legendrian link of $\Lambda_-$ is $\Lambda_{\beta_{\hat{k}}}$, where $\beta_{\hat{k}}:=s_{i_1}\cdots s_{i_{k-1}}s_{i_{k+1}}\cdots s_{i_l}$. Let us move the marked points $p_k^{\pm 1}$ to the right and absorb them into two of the cuspidal marked points $t_i$. This procedure induces an isomorphism 
\[
\Aug\left(\Lambda_-,\cP\right)\xrightarrow{\cong} \Aug\left(\Lambda_{\beta_{\hat{k}}}\right)\times \left(\mathbb{F}^\times\right)_{p_k}.
\]
Denote by $\tau$  the composition 
\[
\Aug\left(\Lambda_{\beta_{\hat{k}}}\right)\times \left(\mathbb{F}^\times\right)_{p_k}\xrightarrow{\cong} \Aug\left(\Lambda_-,\cP\right)\xrightarrow{\phi_+\circ \Phi_{S_k}}\Aug\left(\Lambda_\beta\right).
\]

\begin{prop}\label{commutative diagram for pinch} The following diagram commutes:
\[
\xymatrix{\Aug\left(\Lambda_{\beta_{\hat{k}}}\right)\times \left(\mathbb{F}^\times\right)_{p_k} \ar[r]^{\gamma\times \id}_\cong \ar[d]_\tau & \conf^e_{\beta_{\hat{k}}}\left(\mathcal{C}\right)\times \left(\mathbb{F}^\times\right)_{q_k} \ar[d]^{l^{-1} \circ \psi\circ l} \\
\Aug\left(\Lambda_\beta\right)\ar[r]_{\gamma}^\cong & \conf^e_\beta(\mathcal{C})}
\]
where $l$ is the sequence of left reflections on double Bott-Samelson cells that reflects the first $k-1$ flags from the bottom to the top, and $\psi$ is the open embedding  in \eqref{def open embedding}.
\end{prop}
\begin{proof}
The left  map $\tau$ is composed of an isomorphism corresponding to the migration of the new pair of marked points $p_k^{\pm 1}$ and the functorial morphism $\phi_+\circ \Phi_{S_k}$. We  show that the open embedding
\[\psi:\conf^{s_{i_{k-1}}\cdots s_{i_1}}_{s_{i_{k+1}}\cdots s_{i_l}}(\mathcal{C})\times \left(\mathbb{F}^\times\right)_{q_k}\longrightarrow \conf^{s_{i_{k-1}}\cdots s_{i_1}}_{s_{i_k}\cdots s_{i_l}}(\mathcal{C})\] admits a similar factorization, and prove that the two factorizations coincide under $\gamma$.

Following the proof of Theorem \ref{3.6},  $\gamma$ is defined by setting
$
b_s=q_s$ for  $1\leq s\leq l$,
where $b_s$ are the Reeb coordinates on $\Aug\left(\Lambda_\beta\right)$ and $q_s$ are the affine coordinates on $\conf^e_\beta\left(\mathcal{C}\right)$.  Together with $\Phi_{S_k}^*\left(p_k\right)=b_k$,  we have
\begin{equation}\label{p_k=q_k}
p_k=q_k.
\end{equation}

Now consider the standard representative of a point in the image of $\psi$:
\begin{equation} \label{standard representative for pinch}
\begin{tikzpicture}[baseline=10]
\node (u0) at (0,1.5) [] {$\U_-$};
\node (u1) at (3,1.5) [] {$\U_-Z_{i_{k-1}}\left(q_{k-1}\right)$};
\node (u2) at (6,1.5) [] {$\cdots$};
\node (u3) at (10,1.5) [] {$\U_-Z_{i_1}\left(q_1\right)\dots Z_{i_{k-1}}\left(q_{k-1}\right)$};
\node (d0) at (0,0) [] {$\U_+$};
\node (d1) at (3,0) [] {$Z_{i_k}\left(q_k\right)\U_+$};
\node (d2) at (6,0) [] {$\cdots$};
\node (d3) at (10,0) [] {$Z_{i_k}\left(q_k\right)\dots Z_{i_l}\left(q_l\right)\U_+$};
\draw [->] (u0) -- node [above] {$s_{i_{k-1}}$} (u1);
\draw [->] (u1) -- node [above] {$s_{i_{k-2}}$} (u2);
\draw [->] (u2) -- node [above] {$s_{i_1}$} (u3);
\draw [->] (d0) -- node [below] {$s_{i_k}$} (d1);
\draw [->] (d1) -- node [below] {$s_{i_{k+1}}$} (d2);
\draw [->] (d2) -- node [below] {$s_{i_l}$} (d3);
\draw (u0) -- (d0);
\draw (u3) -- (d3);
\draw [dashed] (-0.5,-0.6) -- (12.25,-0.6) -- (12.25,0.5) -- (0.5,0.5) -- (0.5,1) -- (12.5,1) -- (12.5,2.1) -- (-0.5,2.1) -- cycle;
\end{tikzpicture},
\end{equation}
where  $q_k\neq 0$. We have 
\[
Z_{i_k}\left(q_k\right)=e_{-i_k}\left(q_k^{-1}\right)q_k^{\alpha_{i_k}^\vee}e_{i_k}\left(q_k^{-1}\right).
\]
Let us delete the $\U_+$ at the lower left corner of \eqref{standard representative for pinch}, and act  by $\left(e_{-i_k}\left(q_k^{-1}\right)\right)^{-1}$ on the rest decorated flags. It gives rise to the 
first factor of the preimage of \eqref{standard representative for pinch} under 
 $\psi$:
\[
\begin{tikzpicture}[baseline=10, scale=0.8]
\node (u0) at (-3,1.5) [] {$\U_-$};
\node (u1) at (1,1.5) [] {$\U_-Z_{i_{k-1}}\left(q_{k-1}\right)e_{-i_k}\left(q_k^{-1}\right)$};
\node (u2) at (5,1.5) [] {$\cdots$};
\node (u3) at (10,1.5) [] {$\U_-Z_{i_1}\left(q_1\right)\dots Z_{i_{k-1}}\left(q_{k-1}\right)e_{-i_k}\left(q_k^{-1}\right)$};
\node (d1) at (0,0) [] {$q_k^{\alpha_{i_k}^\vee}e_{i_k}\left(q_k^{-1}\right)\B_+$};
\node (d2) at (4,0) [] {$\cdots$};
\node (d3) at (10,0) [] {$q_k^{\alpha_{i_k}^\vee}e_{i_k}\left(q_k^{-1}\right)Z_{i_{k+1}}\left(q_{k+1}\right)\dots Z_{i_l}\left(q_l\right)\B_+$};
\draw [->] (u0) -- node [above] {$s_{i_{k-1}}$} (u1);
\draw [->] (u1) -- node [above] {$s_{i_{k-2}}$} (u2);
\draw [->] (u2) -- node [above] {$s_{i_1}$} (u3);
\draw [->] (d1) -- node [below] {$s_{i_{k+1}}$} (d2);
\draw [->] (d2) -- node [below] {$s_{i_l}$} (d3);
\draw (u0) -- (d1);
\draw (u3) -- (d3);
\draw [dashed] (-3.5,1) -- (13.9,1) -- (13.9,2.2) -- (-3.5,2.2) -- cycle;
\end{tikzpicture}.
\]
The following procedure  transforms the above configuration to a standard representative.
\begin{enumerate}
    \item[(a)] Move the unipotent factor $e_{-i_k}\left(q_k^{-1}\right)$ inside each decorated flag in the top row all the way to the left so that it can be absorbed into $\U_-$.
    \item[(b)] Move the unipotent factor $e_{i_k}\left(q_k^{-1}\right)$ inside each decorated flag in the bottom row all the way to the right so that it can be absorbed into $\B_+$.
    \item[(c)] Move the torus factor $q_k^{\alpha_{i_k}^\vee}$ inside each decorated flag in the bottom row all the way to the right so that it can be absorbed into $\B_+$.
    \item[(d)] Replace every $\B_+$ in the bottom row by $\U_+$ to obtain a standard representative.
\end{enumerate}
Among these operations, we claim that (a) and (b) correspond to the matrix scanning algorithms for $\Phi_{S_k}^*$, and (c) corresponds to moving the new pair of marked points to the far right after the saddle cobordism $S_k$.

Let us start with (a). Let $q'$ be a collection of $\mathbb{F}$-valued parameters such that
\[
\U_-Z_{i_s}\left(q_s\right)\cdots Z_{i_{k-1}}\left(q_{k-1}\right)e_{-i_k}\left(q_k^{-1}\right)=\U_-Z_{i_s}\left(q'_s\right)\cdots Z_{i_{k-1}}\left(q'_{k-1}\right).
\]
Note that these $q'$ parameters are part of the affine coordinates for $\conf^e_{\beta_{\hat{k}}}(\cC)$. Since we would like to compare the pull-back of the Reeb coordinates on the augmentation varieties versus the affine coordinates on the double Bott-Samelson cells, we need to express the $q$ parameters in terms of the $q'$ parameters. To do so, let us multiply the above equation by $e_{i_k}\left(q_k^{-1}\right)$ on both sides, which yields
\[
\U_-Z_{i_s}\left(q_s\right)\cdots Z_{i_{k-1}}\left(q_{k-1}\right)=\U_-Z_{i_s}\left(q'_s\right)\cdots Z_{i_{k-1}}\left(q'_{k-1}\right)e_{-i_k}\left(q_k^{-1}\right).
\]

We then observe that in order for $Z_i(q)lZ_i(q')$ to hold for $l\in \U_-$, we need 
\begin{equation}\label{l identity}
q=q'+l_{i+1,i}.
\end{equation}
Set $l^{(k-1)}=e_{-i_k}\left(q_k^{-1}\right)$. 
For $s<k$, using \eqref{l identity} recursively, we obtain 
\begin{equation}\label{left q'}
q_s=q'_s+l_{i_s+1,i_s}^{(s)}, \quad \quad \quad \quad l^{(s-1)}=Z_{i_s}\left(q'_s\right)l^{(s)}Z_{i_s}\left(q_s\right)^{-1}\in \U_-.
\end{equation}

Now we turn to (b). The identity we need is 
\[
e_{i_k}\left(q_k^{-1}\right)Z_{i_{k+1}}\left(q_{k+1}\right)\cdots Z_{i_s}\left(q_s\right)\B_+=Z_{i_{k+1}}\left(q'_{k+1}\right)\cdots Z_{i_s}\left(q'_s\right)\B_+,
\]
which is equivalent to
\[
Z_{i_{k+1}}\left(q_{k+1}\right)\cdots Z_{i_s}\left(q_s\right)\B_+=e_{i_k}\left(q_k^{-1}\right)Z_{i_{k+1}}\left(q'_{k+1}\right)\cdots Z_{i_s}\left(q'_s\right)\B_+.
\]
To express the $q$ parameters in terms of the $q'$ parameters, we can first set $u^{(k+1)}=e_{i_k}\left(q_k^{-1}\right)$, and then recursively, we have
\begin{equation}\label{right q'}
q_s=q'_s+u_{i_s,i_s+1}^{(s)} \in \mathbb{F}, \quad \quad \quad \quad u^{(s+1)}:=Z_{i_s}\left(q_s\right)^{-1}u^{(s)}Z_{i_s}\left(q'_s\right) \in \U_+.
\end{equation}

Let $b'$ denote the Reeb coordinates on the Legendrian link after the saddle cobordism but before moving the pair of the newly created marked points to the right. 
We want to show that, under the assumption 
\begin{equation}\label{b'_s=q'_s}
b'_s=q'_s,
\end{equation}
we have 
\begin{equation}\label{b_s=q_s}
b_s=q_s.
\end{equation}

By comparing \eqref{left q'} with \eqref{left-scanning} and  \eqref{right q'} with \eqref{right-scanning},  it suffices to show 
\begin{align*}
\vec{1}+L^{(s)}=&l^{(s)} \quad \text{for $s<k$, and}  \\
\vec{1}+U^{(s)}=&u^{(s)} \quad \text{for $s>k$,}
\end{align*}
where $\vec{1}$ denotes the identity matrix of the appropriate size.

Let us do a backward induction on $s$ to prove the $s<k$ case; the $s>k$ case is similar. The base case $s=k-1$ is clear. By a calculation similar to \eqref{l identity}, for any square matrix $M$ over $\mathbb{F}$ (of characteristic $2$) and any element $x\in \mathbb{F}$, we have
\[
Z_i\left(x\right)\left(\vec{1}+M^-\right)Z_i\left(x+M_{i+1,i}\right)^{-1}=\vec{1}+\left(Z_i\left(x\right)M^-Z_i\left(x+M_{i+1,i}\right)^{-1}\right)^-.
\]
Using this identity, we see that
\begin{align*}
l^{(s-1)}=&Z_{i_s}(q'_s)l^{(s)}Z_{i_s}(q_s)^{-1}\\
=&Z_{i_s}\left(b'_s\right)l^{(s)}Z_{i_s}(b_s)^{-1}\\
=&Z_{i_s}\left(b'_s\right)\left(\vec{1}+L^{(s)}\right)Z_{i_s}\left(b'_s+L^{(s)}_{i_s+1,i_s}\right)^{-1}\\
=&\vec{1}+\left(Z_{i_s}\left(b'_s\right)L^{(s)}Z_{i_s}\left(b'_s+L^{(s)}_{i_s+1,i_s}\right)^{-1}\right)^-\\
=&\vec{1}+L^{(s-1)}.
\end{align*}

For step (c), we claim that moving the torus factor $q_k^{\alpha_{i_k}^\vee}$ through the product
\[
Z_{i_{k+1}}\left(q'_{k+1}\right)Z_{i_{k+2}}\left(q'_{k+2}\right)\cdots Z_{i_s}\left(q'_s\right)
\]
corresponds to moving the new marked points $p_k^{\pm 1}$ to the right through the crossings $s_{i_{k+1}}\cdots s_{i_s}$. 
To include marked points in the braid region, we modify the algorithm to compute the CE dga in Proposition \ref{partial a} by interpolating diagonal matrices from marked points. Observe that moving the new marked points through the crossings $s_{i_{k+1}}\cdots s_{i_s}$ changes the Reeb coordinates $b'_s$ of $\Lambda_-$ to $b''_s$, which are determined by the identity
\[
p_k^{\alpha_{i_k}^\vee}Z_{i_{k+1}}\left(b'_{k+1}\right)Z_{i_{k+2}}\left(b'_{k+2}\right)\cdots Z_{i_s}\left(b'_s\right)=Z_{i_{k+1}}\left(b''_{k+1}\right)Z_{i_{k+2}}\left(b''_{k+2}\right)\cdots Z_{i_s}\left(b''_s\right)D,
\]
where $D$ is a diagonal matrix recording the strand level  of the marked points $p_k^{\pm 1}$. Correspondingly, let $q''_s$ be uniquely chosen such that for all $s>k$,
\[
q_k^{\alpha_{i_k}^\vee}Z_{i_{k+1}}\left(q'_{k+1}\right)Z_{i_{k+2}}\left(q'_{k+2}\right)\cdots Z_{i_s}\left(q'_s\right)\B_+=Z_{i_{k+1}}\left(q''_{k+1}\right)Z_{i_{k+2}}\left(q''_{k+2}\right)\cdots Z_{i_s}\left(q''_s\right)\B_+.
\]
By \eqref{p_k=q_k} and \eqref{b'_s=q'_s}, we deduce that 
\begin{equation}\label{b''_s=q''_s}
b''_s=q''_s, \quad \forall s>k. 
\end{equation}

Note that $\left(b'_1,\dots, b'_{k-1},b''_{k+1}, \dots b''_l\right)$ are the Reeb coordinates on $\Aug\left(\Lambda_{\beta_{\hat{k}}}\right)$, and $\left(q'_1,\dots, q'_{k-1},q''_{k+1},\dots, q''_l\right)$ are the affine coordinates on $\conf^e_{\beta_{\hat{k}}}(\mathcal{C})$. Therefore, \eqref{p_k=q_k}, \eqref{b'_s=q'_s}, and \eqref{b''_s=q''_s} imply the commutativity of the diagram in the proposition.
\end{proof}

\begin{cor}\label{3.30} Let $\left(S_k,\cP\right):\Lambda_-\rightarrow \Lambda_\beta$ be the decorated saddle cobordism that resolves the crossing $i_k$ into a pair of marked points $p_k^{\pm 1}$. Then the functorial morphism $\phi_+\circ \Phi_{S_k}:\Aug\left(\Lambda_-,\cP\right)\rightarrow \Aug\left(\Lambda_\beta\right)$ is an open embedding.
\end{cor}
\begin{proof} Note that in the commutative diagram in Proposition \ref{commutative diagram for pinch}, the top map and the bottom map are both isomorphisms, whereas the map on the right is an open embedding. Therefore the map on the left is also an open embedding.
\end{proof}

\medskip

\noindent\textbf{(II) Cyclic Rotation.} Our next proposition shows that the cyclic rotation morphism between augmentation varieties is equivalent to the composition of a pair of reflections between the double Bott-Samelson varieties. Following \cite{SWflag}, the reflections on double Bott-Samelson varieties are quasi-cluster isomorphisms.

\begin{prop}\label{3.31} The following two diagrams commute:
\begin{equation}\label{cyc rot comm diag}
\xymatrix{\Aug\left(\Lambda_{\beta s_i}\right) \ar[r]^{\gamma}_\cong \ar[d]_{\Phi_\rho}^\cong & \conf^e_{\beta s_i}\left(\mathcal{C}\right) \ar[d]^{l\circ r}_\cong \\
\Aug\left(\Lambda_{s_i\beta}\right) \ar[r]_{\gamma}^\cong & \conf^e_{s_i\beta}(\mathcal{C}),} \quad \quad \quad 
\xymatrix{\Aug\left(\Lambda_{s_i\beta}\right) \ar[r]^{\gamma}_\cong \ar[d]_{\Phi_{\rho^{-1}}}^\cong & \conf^e_{s_i\beta }\left(\mathcal{C}\right) \ar[d]^{r^{-1}\circ l^{-1}}_\cong \\
\Aug\left(\Lambda_{\beta s_i}\right) \ar[r]_{\gamma}^\cong & \conf^e_{\beta s_i}(\mathcal{C}),}
\end{equation}
where $r:\conf^e_{\beta s_i}(\mathcal{C})\stackrel{\cong}{\rightarrow} \conf^{s_i}_\beta(\mathcal{C})$ is the right reflection isomorphism and $l:\conf^{s_i}_\beta(\mathcal{C})\stackrel{\cong}{\rightarrow} \conf^e_{s_i\beta}(\mathcal{C})$ is the left reflection isomorphism.
\end{prop}
\begin{proof} Due to symmetry, it suffices to prove the first commutative diagram.

Suppose that $\beta=s_{i_1}\ldots s_{i_{l-1}}$. The right reflection $r$ maps 
\[
\begin{tikzpicture}[baseline=10,scale=0.67]
\node (d0) at (0,0) [] {$\U_+$};
\node (d1) at (2,0) [] {$x_1\U_+$};
\node (d2) at (3,0) [] {$\cdots$};
\node (d3) at (5,0) [] {$x_{l-1}\U_+$};
\node (d4)  at (8.5,0) [] {$x_{l-1}Z_i\left(q\right)\U_+$};
\node (u) at (4.75,1.5) [] {$\U_-$};
\draw (d0) -- (u) -- (d4);
\draw [->] (d0) -- node [below] {$s_{i_1}$} (d1);
\draw [->] (d2) -- node [below] {$s_{i_{l-1}}$} (d3);
\draw [->] (d3) -- node [below] {$s_i$} (d4);
\end{tikzpicture} \quad \longmapsto \quad 
\begin{tikzpicture}[baseline=10,scale=0.67]
\node (d0) at (0,0) [] {$\U_+$};
\node (d1) at (2,0) [] {$x_1\U_+$};
\node (d2) at (4,0) [] {$\cdots$};
\node (d3) at (6.5,0) [] {$x_{l-1}\U_+$};
\node (u1)  at (5,1.5) [] {$\U_-Z_i\left(q'\right)$};
\node (u0) at (2,1.5) [] {$\U_-$};
\draw (d0) -- (u0);
\draw (u1) -- (d3);
\draw [->] (d0) -- node [below] {$s_{i_1}$} (d1);
\draw [->] (d1) -- node [below] {$s_{i_2}$} (d2);
\draw [->] (d2) -- node [below] {$s_{i_l}$} (d3);
\draw [->] (u0) -- node [above] {$s_i$} (u1);
\end{tikzpicture}. 
\]
Following the definition of $r$, all affine coordinates on $\conf^{s_i}_\beta(\mathcal{C})$ are pulled back to the corresponding affine coordinates except for $q'$. 

Now we compute the pull back $r^*\left(q'\right)$. From the assumption $\xymatrix{\U_-\ar@{-}[r] & x_{l-1}Z_i(q)\U_+}$, we know that $z:=x_{l-1}Z_i(q)$ is Gaussian decomposable, i.e., there exists unique matrices $[z]_\pm \in \U_\pm$ and $[z]_0\in \T$ such that $z=[z]_-[z]_0[z]_+$. We act on the left configuration by $[z]_-^{-1}$, turning it into the picture on the left below. 
\[
\begin{tikzpicture}[scale=0.7,baseline=10]
\node (0) at (0,0) [] {$\cdots$};
\node (1) at (3,0) [] {$[z]_0\U_+$};
\node (2) at (0,1.5) [] {$\U_-$};
\draw [->] (0) -- node [below] {$s_i$} (1);
\draw (1) -- (2);
\end{tikzpicture}\quad \quad \quad \quad
\begin{tikzpicture}[scale=0.7,baseline=10]
\node (0) at (0,1.5) [] {$\U_-$};
\node (1) at (0,0) [] {$\cdots$};
\node (2) at (3,1.5) [] {$\U_-\overline{s}_i$};
\node (3) at (3,0) [] {$[z]_0\U_+$};
\draw [->] (0) -- node [above] {$s_i$} (2);
\draw (3) -- node [right] {$s_i$} (2);
\draw [->] (1) -- node [below] {$s_i$} (3);
\end{tikzpicture}
\]
According to the definition of $r$, the new flag in the top row is the unique flag that is of Tits distance $s_i$ from $\U_-$ and of Tits codistance $s_i$ from $[z]_0\U_+$. It is not hard to see that this flag must be $\U_-\overline{s}_i$. To restore to the standard representative for the preimage, we need to act again by $[z]_-$, which implies that
\[
\U_-Z_i\left(q'\right)=\U_-\overline{s}_i[z]_-^{-1}.
\]
Following the Gaussian elimination process, one can see that
\[
[z]_-=\begin{pmatrix} 1 & 0  & \cdots & 0 & 0\\
\frac{\Delta_1\left(\overline{s}_1^{-1}z\right)}{\Delta_1(z)} & 1  & \cdots & 0 & 0 \\
\vdots & \vdots  & \ddots & \vdots & \vdots \\
* & *  & \cdots & 1 & 0 \\
* & *  & \cdots & \frac{\Delta_{n-1}\left(\overline{s}_{n-1}^{-1}z\right)}{\Delta_{n-1}(z)} & 1\end{pmatrix}.
\]
This implies that $Z_i\left(\frac{\Delta_i\left(\overline{s}_i^{-1}z\right)}{\Delta_i(z)}\right)[z]_-\overline{s}_i^{-1}$ is still a lower triangular unipotent matrix. Therefore, we get
\[
r^*\left(q'\right)=\frac{\Delta_i\left(\overline{s}_i^{-1}z\right)}{\Delta_i(z)}.
\]

On the other hand, since the left reflection map $l$ only moves decorated flags within the compatible region, it follows that the pull-back map $l^*$ is the identity map on affine coordinates. Composing $r^*$ and $l^*$ we get that
\begin{equation}\label{lcr}
(l\circ r)^*\left(q_k\right)=q_{k-1} \quad \forall 1<k\leq l \quad \text{and} \quad 
(l\circ r)^*\left(q_1\right)=\frac{\Delta_i\left(\overline{s}_i^{-1}z\right)}{\Delta_i(z)},
\end{equation}
where
\[
z=Z_{i_1}\left(q_1\right)Z_{i_2}\left(q_2\right)\cdots Z_{i_{l-1}}\left(q_{l-1}\right)Z_i\left(q_{l}\right).
\]

Let us now make use of the natural isomorphism $\gamma$. We first observe that $\gamma^*\left(q_i\right)=b_i$ and $\gamma^*(z)=M:=Z_{i_1}\left(b_1\right)\cdots Z_{i_{l-1}}\left(b_{l-1}\right)Z_i\left(b_l\right)$. Therefore we have
\begin{equation}\label{gamma^*1}
\gamma^*\circ(l\circ r)^*\left(q_1\right)=\frac{\Delta_i\left(\overline{s}_i^{-1}M\right)}{\Delta_i(M)}=\frac{\Delta^{\{1,\dots,i\}}_{\{1,\dots, i-1,i+1\}}(M)}{\Delta_i(M)},
\end{equation}
where $\Delta_I^J$ denotes the determinant of the submatrix formed by the rows in the set $I$ and the columns in the set $J$. By Propositions \ref{cnjona} and \ref{non-vanishing} we can further deduce that
\begin{equation}\label{Delta_i(M)}
\Delta_i(M)=\prod_{k=1}^it_k^{-1} \quad \text{and} \quad \Delta^{\{1,\dots, i\}}_{\{1,\dots, i-1,i+1\}}(M)=M_{i+1,i}^{(i)}\prod_{k=1}^{i-1}t_k^{-1}.
\end{equation}
Combining \eqref{gamma^*1}, \eqref{Delta_i(M)}, and Proposition \ref{3.15}, we obtain the following pull-back image for the affine coordinate $q_1$:
\[
\gamma^*\circ(l\circ r)^*\left(q_1\right)=t_iM_{i+1,i}^{(i)}=\Phi_R^*\circ \gamma^*\left(q_1\right).
\]
For all other affine coordinates $q_k$ with $1<k\leq l$, we can deduce from \eqref{lcr} and Proposition \ref{3.15} that $\gamma^*\circ(l\circ r)^*\left(q_k\right)=\Phi_\rho^*\circ \gamma^*\left(q_k\right)$.
\end{proof}

\medskip

\noindent\textbf{(III) Braid Move.} Suppose $|i-j|=1$ and suppose $\beta'$ and $\beta$ are two braid words that only differ at three consecutive crossings by replacing $(i,j,i)$ with $(j,i,j)$. From the matrix identity
\[
Z_i\left(q_1\right)Z_j\left(q_2\right)Z_i\left(q_3\right)=Z_j\left(q_3\right)Z_i\left(q_2+q_1q_3\right)Z_j\left(q_1\right)
\]
and \eqref{braid move} we deduce that the following diagram commutes
\begin{equation}\label{braid move comm diagram}
\xymatrix{ \Aug\left(\Lambda_{\beta'}\right) \ar[d]_{\Phi_B} \ar[r]^\gamma_\cong &  \conf^e_{\beta'}(\mathcal{C}) \ar@{=}[d]\\
 \Aug\left(\Lambda_{\beta}\right) \ar[r]_\gamma^\cong
& \conf^e_\beta(\mathcal{C}) }
\end{equation}
where $\Phi_B$ is the functorial morphism induced from the braid move Legendrian isotopy $B:\Lambda_{\beta'}\rightarrow \Lambda_{\beta}$. Note that the two $\gamma$ maps are not identical because the top one is defined by the braid word $\beta'$ and the bottom one is defined by the braid word $\beta$.

\medskip 

\noindent\textbf{(IV) Minimum Cobordism.} Consider a decorated Legendrian unknot $O$ with $\tb=-1$ as drawn in Figure \ref{minimum cobordism picture}. The differential of the unique degree 1 Reeb chord $a$ is
\[
\partial a= t_1t_2\cdots t_k+t_m^{-1}t_{m-1}^{-1}\cdots t_{k+1}^{-1}.
\]
Therefore $\Aug(O)$ is the vanishing locus of $t_1t_2\cdots t_k+t_m^{-1}t_{m-1}^{-1}\cdots t_{k+1}^{-1}$ in $\left(\mathbb{F}^\times\right)^m_{t_1,\dots, t_m}$. 

Let $(M,\cP):\emptyset\rightarrow O$ be the decorated minimum cobordism that fills $O$. By definition, $\Aug\left(\emptyset, \cP\right)$ is defined to be the subtorus of $\left(\mathbb{F}^\times\right)^m_{t_1,\dots, t_m}$ satisfying $\prod_i t_i=1$. In characteristic 2, the equation $t_1t_2\cdots t_k+t_m^{-1}t_{m-1}^{-1}\cdots t_{k+1}^{-1}=0$ is equivalent to the equation $\prod_i t_i=1$. Moreover, recall from \eqref{minimum cobordism} that $\Phi_M^*\left(t_i\right)=t_i$ for all $t_i$. Therefore we can conclude the following Lemma.

\begin{lem}\label{3.35} The functorial morphism $\Phi_M:\Aug\left(\emptyset,\cP\right)\rightarrow \Aug(O,\cP)$ is an isomorphism of algebraic tori.
\end{lem}

We are now ready to prove the main theorem of this section.

\begin{thm}\label{thm correspondence} For any admissible filling $L$ of  $\Lambda_\beta$ with decoration $\cP$, the functorial morphism $\phi_+\circ \Phi_L:\Aug\left(\emptyset, \cP\right)\rightarrow \Aug\left(\Lambda_\beta\right)$ is an open embedding of an algebraic torus, and its image  is a cluster chart on $\Aug\left(\Lambda_\beta\right)$.
\end{thm}
\begin{proof} Among the four types of building blocks, we know that cyclic rotations and braid moves are Legendrian isotopies, which are invertible exact Lagrangian concordance. This implies that their induced functorial morphisms between Augmentation varieties are always isomorphisms. Moreover, commutative diagrams \eqref{cyc rot comm diag} and \eqref{braid move comm diagram} yield that $\Phi_{\rho^{\pm 1}}$ and $\Phi_B$ are both quasi-cluster isomorphisms, which map cluster charts to cluster charts. Therefore it suffices to prove the theorem for admissible fillings $L:\emptyset \rightarrow \Lambda_\beta$ that are of the form $S_{k_1}\circ S_{k_2}\circ \cdots S_{k_l}\circ \left(\bigsqcup_n M\right)$, where $l$ is the length of $\beta$ and $n$ is the number of strands in $\beta$.

First we observe that $\Phi_{\bigsqcup_n M}=\prod_n \Phi_M$. Let $\bigsqcup_n O$ be the split union of $n$ decorated Legendrian unknots right before the final minimum cobordisms. Then by Lemma \ref{3.35}, we know that $\Phi_{\bigsqcup_n M}:\Aug(\emptyset,\cP)\rightarrow \Aug\left(\bigsqcup_n O,\cP\right)$ is an isomorphism between algebraic tori. Therefore it remains to show that $\phi_+\circ \Phi_{S_{k_1}}\circ \cdots \circ \Phi_{S_{k_l}}:\Aug\left(\bigsqcup_n O,\cP\right)\longrightarrow \Aug\left(\Lambda_\beta\right)$ is an open embedding from an algebraic torus onto a cluster chart. 

Let us do an induction on the length $l$ of $\beta$. For the base case with $l=1$, the statement follows from Proposition \ref{commutative diagram for pinch} and Corollary \ref{3.30}. For $l>1$, we consider the following commutative diagram:
\[
\xymatrix{
\Aug\left(\bigsqcup_nO,\cP'\right)\times \left(\mathbb{F}^\times\right)_{p_{k_1}} \ar[rrr]^{\left(\phi_+\circ \Phi_{S_{k_2}}\circ \cdots \circ \Phi_{S_{k_l}}\right)\times \id} \ar[d]_\cong & & & \Aug\left(\Lambda_{\beta_{\hat{k}_1}}\right)\times \left(\mathbb{F}^\times\right)_{p_{k_1}} \ar[d]^\cong \ar[dr]^\tau & \\
\Aug\left(\bigsqcup_n O,\cP\right) \ar@{=}[d] \ar[rrr]^{\phi_+\circ \Phi_{S_{k_2}}\circ \cdots \circ \Phi_{S_{k_l}}} & & & \Aug\left(\Lambda_{\beta_{\hat{k}_1}},\cP_{k_1}\right) \ar[r]^(0.6){\phi_+\circ \Phi_{S_{k_1}}} & \Aug\left(\Lambda_\beta\right) \ar@{=}[d] \\
\Aug\left(\bigsqcup_n O,\cP\right) \ar[rrrr]_{\phi_+\circ \Phi_{S_{k_1}}\circ \cdots \circ \Phi_{S_{k_l}}} & & & & \Aug\left(\Lambda_\beta\right)
}
\]
where $\cP'=\cP\setminus \left\{p_{k_1}\right\}$, and $\cP_{k_1}$ denotes the decoration on the saddle cobordism $S_{k_1}$. By the inductive hypothesis, we know that the top morphism is an open embedding onto a cluster chart. On the other hand, Proposition \ref{commutative diagram for pinch} and Corollary \ref{3.30} implies that $\phi_+\circ \Phi_{S_{k_1}}$ is an open embedding and a quasi-cluster morphism. Therefore, it follows from the commutative that the bottom morphism is also an open embedding onto a cluster chart. This finishes the proof of the theorem.
\end{proof}

\begin{cor}\label{all Reeb chords are cluster variables} Every degree 0 Reeb chord $b_k$ of $\Lambda_\beta$ is a mutable cluster variable of the cluster structure on $\Aug(\Lambda_\beta)$.
\end{cor}
\begin{proof} From the proof of Proposition \ref{commutative diagram for pinch} and Theorem \ref{thm correspondence}, we see that $b_k$ is a mutable cluster coordinate on the cluster chart corresponding to the admissible filling $L\circ S_k$ where $L$ is any admissible filling of $\Lambda_{\beta_{\hat{k}}}$.
\end{proof}

\begin{cor}\label{cor 3.40} Suppose $L$ and $L'$ are Hamiltonian isotopic admissible fillings of $\Lambda_\beta$, then they give rise to the same cluster seed. 
\end{cor}
\begin{proof} 
By construction, any admissible filling $(L,\cP)$ has sufficient $\cP$. By Corollary \ref{sameimage}, the cluster charts corresponding to $L$ and $L'$ are equal as open subvarieties. By Proposition \ref{distinguishing charts}, we know that $L$ and $L'$ must correspond to the same cluster seed.
\end{proof}

 The theory of cluster algebras gives rise to a computable numerical invariant for each admissible filling.
 Let $\alpha_0$ be the cluster seed associated to the admissible filling 
\[
L_0:=\left(\bigsqcup_n M\right)\circ S_l\circ S_{l-1}\circ \cdots \circ S_2\circ S_1.
\]
We set $\alpha_0$ as the initial cluster seed. Fomin and Zelevinsky  \cite[(6.4)]{FZIV} constructed an integer matrix $G_\alpha$, called the \emph{$g$-matrix}, for every cluster seed $\alpha$. Following \cite{GHKK}, each $\alpha$ corresponds to a cluster chamber $\mathcal{C}_\alpha$ in the scattering diagram associated with the cluster algebra, and the column vectors of $G_\alpha$ are the primitive vectors  spanning $\mathcal{C}_\alpha$. Thus, the sums of the column vectors of the $g$-matrices are a complete invariant for the cluster seeds. We conclude the following corollary.

\begin{cor} 
\label{cor:g-vector-inv}
For each admissible filling $L$, let $\alpha_L$ be its corresponding cluster seed and let $G_L$ be the $g$-matrix of $\alpha_L$ with respect to the initial seed $\alpha_0$. Let $g_L$ be the sum of column vectors of $G_L$. If $L$ and $L'$ are Hamiltonian isotopic, then $g_L=g_{L'}$.
\end{cor}

\smallskip

\subsection{Computing Cluster Seeds associated with Admissible Fillings}\label{sec algorithm}

In this section, we present an explicit algorithm to compute the cluster seeds (including their cluster coordinates and quivers) associated with admissible fillings. Throughout this section, we fix an $n$-stranded braid word $\beta=s_{i_1}\dots s_{i_l}$. 

\medskip

\noindent\textbf{(0) Initial Seed.} Let us first consider the cluster chart  that is the image of the functorial morphism $\phi_+\circ \Phi_{L_0}$ with
$L_0=\left(\bigsqcup_n M\right)\circ S_l\circ S_{l-1}\circ \cdots \circ S_2\circ S_1$. The following statement is a direct consequence of Proposition \ref{commutative diagram for pinch} and Theorem \ref{thm correspondence}.

\begin{prop}\label{prop initial seed} Under the  isomorphism $\gamma:\Aug(\Lambda_\beta)\rightarrow \conf^e_\beta(\mathcal{C})$, the cluster seed $\alpha_0$ is identified with the unique triangulation defined by the braid word $\beta$ on $\conf^e_\beta(\mathcal{C})$.
\end{prop} 

The cluster coordinates on $\alpha_0$ are
\[
A_k=\Delta_{i_k}\left(Z_{i_1}\left(b_1\right)Z_{i_2}\left(b_2\right)\cdots Z_{i_k}\left(b_k\right)\right), \quad \quad \forall 1\leq k\leq l.
\]
Comparing  $Q_\beta$  with the quivers associated with triangulations for $\conf^e_\beta(\mathcal{C})$ in Appendix \ref{Sec 3}, we see that $Q_\beta$ is precisely the quiver for the initial cluster seed $\alpha_0$. Note that the cluster coordinate $A_k$ is associated with the region (quiver vertex) to the immediate right of the $k$th crossing, and the cluster coordinates that are on the furthest right on each horizontal level are automatically frozen. We call  $\alpha_0=\left(\left\{A_k\right\}_{1\leq k\leq l},Q_\beta\right)$ the \emph{initial seed} and $Q_\beta$ the \emph{initial quiver} associated with the braid word $\beta$. Other cluster seeds can be obtained from the initial seed via a sequence of cluster mutations, and we will describe an explicit cluster mutation sequence for each building block of admissible fillings.

\medskip

\noindent\textbf{(I) Saddle Cobordism.} We make use of Proposition \ref{commutative diagram for pinch} to derive the mutation sequence for saddle cobordisms. From this proposition we know that a saddle cobordism $S_k:\Lambda_{\beta_{\hat{k}}}\rightarrow \Lambda_\beta$ corresponds an open quasi-cluster morphism. In order to get the image, which is an open cluster subvariety, we need to 
\begin{enumerate}
    \item apply a sequence of left reflection maps $l$ that reflects first $k-1$ flags from the bottom to the top;
    \item perform the open embedding $\psi$ described in Appendix \ref{sec 3.3};
    \item apply the inverse sequence of left reflection maps $l^{-1}$.
\end{enumerate}
Our goal is to produce the initial quiver $Q_{\beta_{\hat{k}}}$ for the positive braid Legendrian link $\Lambda_{\beta_{\hat{k}}}$ (without marked points in the braid region). The mutation sequence to turn $Q_\beta$ to $Q_{\beta_{\hat{k}}}$ will be a composition of mutation sequences that correspond to the three steps above.

Since (2) involves setting aside a quiver vertex that will no longer be considered as part of the quiver for $\Lambda_{\beta_{\hat{k}}}$, we introduce a new concept called active vertices for the quivers associated with admissible fillings. 

\begin{defn} Let $L:\Lambda_-\rightarrow \Lambda_\beta$ be an admissible cobordism. An unfrozen quiver vertex is said to be \emph{active} if it is still considered as part of the quiver for the positive braid Legendrian link $\Lambda_-$ after disregarding all the marked point in the braid region. A quiver vertex is said to be \emph{inactive} if it is not active. 
\end{defn}

Note that in the initial quiver $Q_\beta$, all unfrozen vertices are active.

Let us now describe the mutation sequences for each of the three steps involved in locating the open cluster subvariety.

\begin{enumerate}[leftmargin=*]
    \item In terms of the triangulation description of cluster seeds in double Bott-Samelson cells, each left reflection in $l$ reflects a flag from the bottom left hand corner to the top left hand corner by turning the left most triangle upside down. But then in order to prepare for the next left reflection, we should move this newly turned triangle to the right of the triangle with base $\xymatrix{\B^k\ar[r]^{s_{i_k}} & \B^{k+1}}$ using cluster mutations. 
    
    Let us denote the active quiver vertices on the $i$th level as 
    \[
    \textstyle\binom{i}{1}, \binom{i}{2},\dots, \binom{i}{m_i}
    \]
    from left to right. For each level $i$ and two integers $a,b$ satisfying $1\leq a\leq b\leq m_i$, we define a mutation sequence 
    \begin{equation}\label{def eta}
    \eta(i,a,b):= \mu_{\binom{i}{b}}\circ \mu_{\binom{i}{b-1}}\circ \cdots \circ \mu_{\binom{i}{a}}.
    \end{equation}
   For each crossing $i_j$ in the braid $\beta$ with $1\leq j<k$, we define 
    \[
    t_j:= \#\left\{s \ \middle| \  j<s\leq k, i_s=i_j\right\}.
    \]
    The sequence of left reflections $l$ corresponds to the sequence of mutations:
    \begin{equation}\label{E_l}
    E_l:=\eta\left(i_{k-1},1,t_{k-1}\right)\circ \cdots \circ\eta\left(i_2,1,t_2\right)\circ\eta\left(i_1,1,t_1\right).
    \end{equation}
    
    \item In this step, we need to remove the left most triangle, which has base $\xymatrix{\B^k\ar[r]^{s_{i_k}} & \B^{k+1}}$, from the triangulation. This corresponds to deactivating the left most active vertex on the $i_k$th level. Due to this deactivation, there will be one fewer active vertex on the $i_k$th level. To avoid confusion, let us denote the new braid by $\beta'$ and denote the active quiver vertices on the $i$th level as $\binom{i}{1}', \binom{i}{2}',\dots, \binom{i}{m_i'}'$. Note that
    \[
    \textstyle\binom{i}{a}' = \left\{\begin{array}{ll}
        \binom{i}{a} & \text{if $i\neq i_k$},  \\
        \binom{i}{a+1} & \text{if $i=i_k$}. 
    \end{array}\right.
    \]
    
    \item Note that $\beta'= \beta_{\hat{k}} = s_{i_1}\cdots s_{i_{k-1}}s_{i_{k+1}}\cdots s_{i_l}$. For each $i_j$ with $1\leq j<k$, we define
    \[
    t'_j:=\#\left\{s\ \middle| \ j<s<k, i_s=i_j\right\}.
    \]
    Define $\eta'(i,a,b)$ similar to \eqref{def eta} with each mutation indexed by primed indices $\binom{i}{s}'$. The muatation sequence corresponding to the sequence of left reflections $l^{-1}$ is
    \begin{equation}\label{E_l^-1}
    E_{l^{-1}}:=\eta'\left(i_1,1,t'_1\right)^{-1}\circ \eta'\left(i_2,1,t'_2\right)^{-1}\circ \cdots \circ \eta'\left(i_{k-1},1,t'_{k-1}\right)^{-1}.
    \end{equation}
\end{enumerate}

Combining the three steps, the total mutation sequence for a saddle cobordism $S_k$ is 
\begin{equation}\label{saddle mut seq}
E_{S_k}:=E_{l^{-1}}\circ E_l,
\end{equation}
where $E_l$ is defined in \eqref{E_l} and $E_{l^{-1}}$ is defined in \eqref{E_l^-1}.

\medskip

\noindent\textbf{(II) Cyclic Rotation.} Let $\rho_i:\Lambda_{\beta s_i}\rightarrow \Lambda_{s_i\beta}$ be a cyclic rotation cobordism. According to Proposition \ref{3.31}, the functorial morphism $\Phi_{\rho_i}:\Aug\left(\Lambda_{\beta s_i}\right)\rightarrow \Aug\left(\Lambda_{s_i\beta}\right)$ corresponds to the composition 
\[
\conf^e_{\beta s_i}(\mathcal{C})\xrightarrow{r} \conf^{s_i}_\beta(\mathcal{C})\xrightarrow{l} \conf_{s_i\beta}(\mathcal{C}).
\]
The change of initial quiver associated with this composition of reflection maps can be realized via a mutation sequence that mutates every active quiver vertex on the $i$th level. When we left-compose $\rho_i$ onto an admissible cobordism, we are changing from the initial quiver of $\Aug\left(\Lambda_{s_i\beta}\right)$ to the initial quiver for $\Aug\left(\Lambda_{\beta s_i}\right)$. Therefore the corresponding mutation sequence is
\begin{equation} \label{R mut seq}
E_{\rho_i}:=\eta\left(i,1,m_i\right),
\end{equation}
where $\eta$ is defined in \eqref{def eta}. Consequently, 
\begin{equation}\label{R-1 mut seq}
E_{\rho_i^{-1}}:=\eta\left(i,1,m_i\right)^{-1}.
\end{equation}

\medskip

\noindent\textbf{(III) Braid Move.} From the commutative diagram \eqref{braid move comm diagram} we know that a braid move cobordism $B:\Lambda_\beta'\rightarrow \Lambda_\beta$ corresponds to a braid move on the bases of the corresponding double Bott-Samelson cell triangulation. It is known that the latter is a single mutation that takes place at a unique quiver vertex. In terms of the initial quiver $Q_\beta$, this unique quiver vertex is the unique vertex that is associated with the region completely enclosed by the three strands involved in the braid move. Therefore we conclude that
\begin{equation}\label{braid move mut seq}
E_B:=\mu_c.
\end{equation}
\begin{figure}[H]
    \centering
\begin{tikzpicture}[baseline=10]
\draw (0,0) -- (1,0) to [out=0,in=180] (2,0.5) to [out=0,in=180] (3,1);
\draw (0,0.5) to [out=0,in=180] (1,1) -- (2,1) to [out=0,in=180] (3,0.5);
\draw (0,1) to[out=0,in=180] (1,0.5) to [out=0,in=180] (2,0) -- (3,0);
\node (c) at (1.5,0.75) [] {$c$};
\draw [->] (c) -- (3,0.75);
\draw [->] (0,0.75) -- (c);
\draw [->] (c) -- (0.5,0.25);
\draw [->] (2.5,0.25) -- (c);
\draw [->] (1,0.25) -- (2,0.25);
\end{tikzpicture} \quad \quad $\overset{\mu_c}{\longleftrightarrow}$\quad \quad
\begin{tikzpicture}[baseline=10]
\draw (0,0)  to [out=0,in=180] (1,0.5) to [out=0,in=180] (2,1)  -- (3,1);
\draw (0,0.5) to [out=0,in=180] (1,0) -- (2,0) to [out=0,in=180] (3,0.5);
\draw (0,1)  -- (1,1) to [out=0,in=180] (2,0.5) to [out=0,in=180] (3,0);
\node (c) at (1.5,0.25) [] {$c$};
\draw [->] (c) -- (3,0.25);
\draw [->] (0,0.25) -- (c);
\draw [->] (c) -- (0.5,0.75);
\draw [->] (2.5,0.75) -- (c);
\draw [->] (1,0.75) -- (2,0.75);
\end{tikzpicture}
\caption{Braid Move}
    \label{}
\end{figure}
Note that after a braid move, the active vertex $c$ needs to move to the adjacent level, as depicted in the picture above.

\medskip

\noindent\textbf{(IV) Minimum Cobordism.} A minimum cobordism $M$ induces an isomosphism $\Phi_M$ between algebraic tori. Therefore it corresponds to the empty mutation sequence, i.e., \begin{equation}\label{minimum mut seq}
    E_M=\emptyset.
\end{equation}

\medskip

\noindent\textbf{(V) Summary.} For any admissible filling $L$ of $\Lambda_\beta$, the corresponding cluster seed $\alpha_L$ can be computed as follows. First we compute the initial seed $\alpha_0$ associated with the braid word $\beta$; then we write $L$ as a composition of elementary building blocks $L=L_m\circ \cdots \circ L_2\circ L_1$, and mutate the initial seed $\alpha_0$ accordingly, yielding
\[
\alpha_L:=E_{L_m}\circ \cdots \circ E_{L_2}\circ E_{L_1}\left(\alpha_0\right).
\]
Each mutation subsequence $E_{L_i}$ is given by one of \eqref{saddle mut seq}, \eqref{R mut seq}, \eqref{R-1 mut seq}, \eqref{braid move mut seq}, and \eqref{minimum mut seq}. We have implemented a characteristic 0 version of this algorithm in a javascript program.\footnote{See \url{https://users.math.msu.edu/users/wengdap1/filling_to_cluster.html}.} For any admissible filling $L$, this program computes 
    \begin{itemize}[leftmargin=2.1em]
        \item the functorial homomorphism images $\Phi_L^*\left(b_i\right)$ for all degree 0 Reeb chords;
        \item the mutation sequence from the initial cluster $\alpha_0$ to the cluster $\alpha_L$;
        \item the cluster seed of $\alpha_L$, including both the cluster variables and the associated quiver;
        \item the seed invariant vector $g_L$ (Corollary \ref{cor:g-vector-inv}).
    \end{itemize}

\bigskip

\section{Infinitely Many Fillings}

In this section, we solve the infinite-filling problem for positive braid Legendrian links. 
One key ingredient in our proof is the cluster Donaldson-Thomas transformations. 
Throughout this section, all mentions of the quiver $Q_\beta$ refer to its unfrozen part.

\subsection{Full Cyclic Rotation and Donaldson-Thomas Transformation}

\begin{defn} For a positive braid word $\beta$ of length $l$, the \emph{full cyclic rotation} $\R$ is the exact Lagrangian concordance $\rho^l: \Lambda_\beta \rightarrow \Lambda_\beta$, where $\rho$ is the cyclic rotation \eqref{def rho}. 
\end{defn}

The cluster DT transformation is a unique central element of the cluster modular group acting on the associated cluster varieties (Definition \ref{def DT}). Combinatorially, the cluster DT transformation can be manifested as a maximal green sequence, or more generally, a reddening sequence of quiver mutations \cite{KelDT}.

\begin{lem}\label{6.7} For any positive braid word $\beta$, we have  $\Phi_\R=\DT^{2}$ on $\Aug\left(\Lambda_\beta\right)$.
\end{lem}

\begin{proof} Suppose $\beta=s_{i_1}\dotsb s_{i_l}$. By \cite{SWflag}, the DT transformation on  $\conf^e_\beta(\cC)$ is 
\[
\DT=t\circ \left(r_{i_1}\circ r_{i_2} \circ \dots \circ r_{i_l}\right),
\]
where $t$ is a biregular isomorphism induced by the transposition action on $\G=\SL_n$ and $r_i$ are right reflection maps. Let us denote the left reflection of $s_i$ by $l^i$. Then
\begin{align*}
\DT^{2}
=&t\circ \left(r_{i_1}\circ r_{i_2} \circ \dots \circ r_{i_l}\right)\circ t\circ\left(r_{i_1}\circ r_{i_2} \circ \dots \circ r_{i_l}\right)\\
=&t\circ t\circ \left(l^{i_1}\circ l^{i_2}\circ \cdots\circ l^{i_l}\right)\circ \left(r_{i_1}\circ r_{i_2} \circ \dots \circ r_{i_l}\right)\\
=&\left(l^{i_1}\circ l^{i_2}\circ \cdots\circ l^{i_l}\right)\circ \left(r_{i_1}\circ r_{i_2} \circ \dots \circ r_{i_l}\right) \\
=&\left(l^{i_1}\circ r_{i_1}\right)\circ \left(l^{i_2}\circ r_{i_2}\right)\circ \cdots \circ \left(l^{i_l} \circ r_{i_l}\right).
\end{align*}
The first commutative diagram in Proposition \ref{3.31} asserts that $l^{i_k}\circ r_{i_k}=\Phi_\rho$. Therefore $\DT^2=\Phi_{\rho^l}=\Phi_\R$.
\end{proof}

\begin{thm}\label{5.11} For any positive braid word $\beta$, if the $\DT$ transformation on $\Aug\left(\Lambda_\beta\right)$ is aperiodic, then $\Lambda_\beta$ admits infinitely many admissible fillings.
\end{thm}
\begin{proof} Let $L_0$ be the admissible filling that pinches the crossings in $\beta$ from left to right and then fills the resulted unlinks with minimum cobordisms. Let $L_m = \R^{m}\circ L_0$. We claim that $L_m$ is not Hamiltonian isotopic to $L_k$ for $m \neq  k$. To see this, note that by Lemma \ref{6.7}, the cluster seeds of $L_m$ can be computed by mutating the initial seed according to $\DT^{2m}$; the aperiodicity of $\DT$ implies that the cluster seeds of $L_m$ and $L_k$ are distinct for $m\neq k$. The statement follows from Corollary \ref{cor 3.40}.
\end{proof}

\begin{rmk}\label{periodicDT}

The full cyclic rotation was observed by K\'alm\'an \cite{Kalman}. For torus links $\Lambda_{(n,m)}$, where $\beta = (s_1s_2 \dotsb s_{n-1})^m$, \cite{Kalman} further defined another Legendrian loop $\K=\rho^{n-1}$, with the property $\R = \K^{m}$. K\'alm\'an showed that $\Phi_\K$ has finite order.

The quivers associated to $\Aug\left(\Lambda_{(n,m)}\right)$ and those associated to the Grassmannian $\Gr_{n,n+m}$ share the same unfrozen parts. Hence, their DT transformations have the same order. The $\DT$ on $\Gr_{n,n+m}$ has finite order because it is related to the periodic Zamolodchikov operator by $\DT^2=\mathrm{Za}^{m}$ \cite{Kelperiod,weng, SWflag}. In fact, K\'alm\'an's loop induces the Zamolodchikov operator. Summarizing,
\[
\Phi_\R = \Phi_{\K}^m = \DT^2  = \mathrm{Za}^{m}.
\]
\end{rmk}{}

\begin{thm}\label{6.13} Let $Q$ be an acyclic quiver. Its associated $\DT$ transformation is of finite order if and only if $Q$ is of finite type. 
\end{thm}
\begin{proof} Combinatorially, the DT transformation arises from a maximal green sequence of quiver mutations \cite{KelDT}. When $Q$ is acyclic, one may label the vertices of $Q$ by $1,\ldots, l$ such that $i<j$ if there is an arrow from $i$ to $j$. The mutation sequence $\mu_n\circ\cdots \circ \mu_1$ is  maximal green and therefore gives rise to the DT transformation associated with $Q$.

The DT transformation acts on the cluster variety $\mathscr{A}_Q$ associated with the quiver $Q$. Following \cite{lee2018frieze}, the {\it frieze variety} $X(Q)$ is defined to be the Zariski closure of the DT-orbit containing the point $P=(1, \ldots, 1)\in \mathscr{A}_Q$. 
Theorem 1.1 of {\it loc.cit.} states that
\begin{enumerate}
    \item If $Q$ is representation finite (i.e., the underlying graph of $Q$ is a Dynkin diagram of type ADE), then the frieze variety $X(Q)$ is of dimension $0$.
    \item If $Q$ is tame then the frieze variety $X(Q)$ is of dimension $1$.
    \item If $Q$ is wild then the frieze variety $X(Q)$ is of dimension at least $2$.
\end{enumerate}{}
As a direct consequence, if $Q$ is not of finite type, then the DT-orbit of $P$ contains infinitely many points, and therefore DT is not periodic. 
If $Q$ is of finite type, then its cluster variety is of finite type, and therefore its DT transformation is periodic.
\end{proof}

\begin{rmk} Keller pointed out to us that the aperiodicity of $\DT$ for acyclic quiver $Q$ of infinite type follows from the aperiodicity of the Auslander-Reiten translation functor on the derived category of representations of $Q$.
\end{rmk}

\begin{cor}\label{rainbowinfinite}
For any positive braid word $\beta$, if $Q_\beta$ is acyclic and of infinite type, then $\Lambda_\beta$ admits infinitely many admissible fillings.
\end{cor}
\begin{proof} It follows from Theorem \ref{5.11} and Theorem \ref{6.13}.
\end{proof}

\vskip 1mm 

\subsection{Infinitely Many Fillings for Infinite Type}\label{sec inifnite fillings}

This section is devoted to the proof of the following result.

\begin{thm}\label{MainStep1}
If $[\beta]$ is a positive braid of infinite type, then the positive braid Legendrian link $\Lambda_\beta$ admits infinitely many non-Hamiltonian isotopic exact Lagrangian fillings.  
\end{thm}

\begin{defn} Given two positive braid words $\beta$ and $\gamma$, we say $\beta$ \emph{dominates} $\gamma$ if there is an admissible cobordism from $\Lambda_\gamma$ to $\Lambda_\beta$. Dominance is a partial order on braid words.
\end{defn}

Recall that a quiver is \emph{connected} if its underlying graph is connected. Connectedness of quivers is invariant under  mutations. Under the connectedness assumption, Theorem \ref{MainStep1} is a consequence of Corollary \ref{rainbowinfinite} and the following Propositions.

\begin{prop}\label{3.2} Suppose $\beta$ dominates $\gamma$. If $\Lambda_\gamma$ admits infinitely many  admissible fillings, then so does $\Lambda_\beta$. 
\end{prop}
\begin{proof} Recall from Corollary \ref{cor 3.40} that the cluster seeds can be used to distinguish admissible fillings. Since the functorial morphism between augmentation varieties induced by any admissible cobordism is a cluster morphisms, it must map distinct cluster seeds to distinct cluster seeds.
\end{proof}

\begin{prop}\label{Mainquiver}
For any braid word $\beta$ with connected $Q_\beta$, either one of the following two scenarios happens:
\begin{enumerate}[label*=\emph{(\arabic*)},leftmargin = *]
    \item there is an admissible concordance from $\Lambda_\gamma$ to $\Lambda_\beta$ and $Q_\gamma$ is a quiver of finite type.
    \item $\beta$ dominates a braid word $\gamma$ and $Q_\gamma$ is acyclic and of infinite type. 
\end{enumerate} 
\end{prop}

\begin{prop}\label{3.4} If Proposition \ref{Mainquiver} (1) happens, then $[\beta]$ is of finite type.

If Proposition \ref{Mainquiver} (2) happens, then $[\beta]$ is of infinite type.
\end{prop}
\begin{proof} Admissible concordances give rise to sequences of mutations (Section \ref{sec algorithm}). 
If Proposition \ref{Mainquiver} (1) happens, then $Q_\beta$ is mutation equivalent to $Q_\gamma$. The latter is of finite type. Therefore $[\beta]$ is of finite type. 

If Proposition \ref{Mainquiver} (2) happens, then by Theorem \ref{thm correspondence}, $Q_\beta$ is mutation equivalent to a quiver which contains $Q_\gamma$ as a full subquiver. Suppose that $[\beta]$ is of finite type. Then $Q_\gamma$ is mutation equivalent to finite type quiver, which contradicts with the assumption that $Q_\gamma$ is acyclic and of infinite type. Therefore $[\beta]$ is of infinite type. 
\end{proof}

Proposition \ref{3.4} implies the exclusiveness of the two scenarios of Proposition \ref{Mainquiver}. To conclude the proof of Proposition \ref{Mainquiver}, it remains to prove that the two scenarios cover all braid words with connected quivers. The strategy of our proof is as follows.
\begin{itemize}
    \item Suppose there is an admissible concordance $\Lambda_\gamma\rightarrow \Lambda_\beta$  such that $Q_\gamma$ is acyclic. If $Q_\gamma$ is of finite type,  then $\beta$ satisfies (1); otherwise, $\beta$ satisfies (2). 
    \item Otherwise, we prove that $\beta$ satisfies (2).
\end{itemize}

\medskip

\noindent\textbf{(I) Preparation.} We adopt the following notations for  operations on braid words.
\begin{enumerate}
    \item[1.] $\overset{\text{R1}}{=}$ denotes the positive Markov destabilization, which deletes the $s_1$ (resp. $s_{n-1}$) if it only occurs once in $\beta$.
   
    \item[2.] $\overset{\text{R3}}{=}$ denotes the braid move R3, which switches $s_is_{i+1}s_i$ and $s_{i+1}s_is_{i+1}$.
    
    \item[3.] $\overset{\rho}{=}$ denotes the cyclic rotation, which turns $\beta s_i$ into $s_i\beta$ or vice versa.
    
    \item[4.] $\overset{c}{=}$ denotes the commutation which turns $s_is_j$ into $s_js_i$ whenever $|i-j|>1$.

    \item[5.] $\succ$ denotes deleting letters; $\beta\succ \gamma$ means that $\gamma$ can be obtained by deleting letters in $\beta$. In particular, when $\beta\succeq \gamma$, we say that $\gamma$ is a \emph{subword} of $\beta$.
    
    \item[6.] $\overset{\textrm{oppo}}{\rightsquigarrow}$ denotes taking the opposite word ${\beta}^{\mathrm{op}}$. The quiver $Q_{\beta^{\mathrm{op}}}$ alters the orientation of every arrow in $Q_\beta$. 
\end{enumerate}
Operations 1 - 4 induce Legendrian isotopies between corresponding positive braid Legendrian links, which are building blocks for admissible concordance. Operations  5 induces pinch cobordisms between Legendrian links. Operation 6 is a symmetry that can be used to reduce the number of cases considered in the proof.

\begin{lem}\label{basicbraidsinf}
The quivers for the following braids are acyclic and of infinite type:
\begin{enumerate}[label*=\emph{(\arabic*)}]
    \item $s_1^2s_2^2s_1^2s_2^2$, or more generally, $s_i^2s_{i+1}^2s_i^2s_{i+1}^2$;
    \item $s_1s_3s_2^2s_1s_3s_2^2$.
\end{enumerate}
\end{lem}
\begin{proof}
The quivers for (1) and (2) are $\tilde{\mathrm{D}}_5$ and $\tilde{\mathrm{D}}_4$ respectively. 
\[
\begin{tikzpicture}
	\draw [thick,->](0.15,0) -- (0.85, 0);
	\filldraw (0,0) circle (2pt);
	\filldraw (1,0) circle (2pt);
	\filldraw (-1,0) circle (2pt);
	\filldraw (0,-1) circle (2pt);
	\filldraw (1,-1) circle (2pt);
	\filldraw (-1,-1) circle (2pt);
    \draw [thick,->] (0.15,-1) -- (0.85, -1);
	\draw [thick,<-](0,-0.15) -- (0, -0.85);
	\draw [thick,<-](-0.15,0) -- (-0.85, 0);
	\draw [thick,->](-0.85,-1) -- (-0.15,-1);
\end{tikzpicture}
\qquad\qquad\qquad
\begin{tikzpicture}[baseline = -15]
	\draw [thick,->](0.15,0) -- (0.85, 0);
	\filldraw (0,0) circle (2pt);
	\filldraw (1,0) circle (2pt);
	\filldraw (-1,0) circle (2pt);
	\filldraw (0,1) circle (2pt);
	\filldraw (0,-1) circle (2pt);
	\draw [thick,<-](0,0.15) -- (0, 0.85);
	\draw [thick,<-](0,-0.15) -- (0, -0.85);
	\draw [thick,<-](-0.15,0) -- (-0.85, 0);
\end{tikzpicture} \qedhere
\]
\end{proof}

\begin{lem}  \label{ws4lemma}
Suppose $w_1,w_2,w_3, w_4\in \left\{ s_1s_3, s_1^2, s_3^2\right\}$. Then $w_1s_2w_2s_2w_3s_2w_4s_2$ dominates a braid with an acyclic quiver of infinite type.
\end{lem}
\begin{proof} Note that $\beta = w_1 s_2 {\color{red}{w_2}} s_2 w_3 s_2 {\color{red}{w_4}} s_2 \succ w_1s_2^2 w_3s_2^2.$
If $w_1=w_3$, then the Lemma follows from Lemma \ref{basicbraidsinf}. The same argument applies to  $w_2=w_4$. 
In the rest of the proof, we assume that $w_1 \neq w_3$ and  $w_2 \neq w_4$.

Let $k$ be the size of the set $\{i \mid w_i = s_1s_3\}$. Here $k\leq 2$; otherwise,  $w_1 = w_3$ or  $w_2 = w_4$. Using the symmetry between $s_1$ and $s_3$, we further assume that there are more $s_1^2$ than $s^2_3$ in $\{w_1, w_2, w_3, w_4\}$. 
We shall exhaust all the possibilities of $k$.

\vskip 2mm

\paragraph{{\it Case 1: k=2}} After taking necessary cyclic rotations and/or the opposite word, we have $w_1 = w_2 = s_1s_3$, and the values of $w_3, w_4$ split into two subcases.

If $w_3 =s_1^2$ and $w_4 = s_3^2$, then $\Lambda_\beta$ is admissibly concordant to the standard $\mathrm{E}_9$ link:
                \begin{align*}
                    \beta
                    &= s_1s_3s_2s_3s_1s_2s_1s_1 {\color{teal}{s_2s_3s_3s_2}}
                    \stackrel{\rho}{=}
                    s_2s_3s_3s_2s_1 {\color{blue}{s_3s_2s_3}} s_1s_2s_1s_1 \\
                    &\stackrel{\textrm{R3}}{=}
                    s_2s_3s_3 {\color{blue}{s_2s_1s_2}} s_3s_2 s_1s_2s_1s_1
                    \stackrel{\textrm{R3}}{=}
                    s_2{\color{blue}{s_3s_3s_1}}s_2 {\color{blue}{s_1s_3}}s_2 s_1s_2s_1s_1 \\
                    &
                    \stackrel{c}{=}
                    s_2s_1{\color{blue}{s_3s_3s_2s_3}}s_1s_2 s_1s_2s_1s_1 
                    \stackrel{\textrm{R3}}{=}
                    s_2s_1s_2{\color{teal}{s_3}}s_2s_2s_1s_2 s_1s_2s_1s_1    \\
                    &\stackrel{\textrm{R1}}{=} {\color{blue}{s_2s_1s_2}}s_2s_2s_1{\color{blue}{s_2 s_1s_2}}s_1s_1 
                    \stackrel{\textrm{R3}}{=}
                    s_1{\color{blue}{s_2s_1s_2s_2}}s_1s_1 s_2s_1s_1s_1    \\
                    &\stackrel{\textrm{R3}}{=}
                    s_1s_1s_1s_2s_1s_1s_1 s_2s_1s_1s_1 
                    =
                    s_1^3 s_2s_1^3 s_2 {\color{teal}{s_1^3}} 
                    \stackrel{\rho}{=} s_1^6 s_2s_1^3 s_2
                \end{align*}
 
 If ${w_3 =w_4 =s_1^2}$, we make the following moves and then apply Lemma \ref{basicbraidsinf} (1):
            \begin{align*}
                \beta 
                &= s_1 {\color{blue}{s_3 s_2 s_3}} s_1 s_2 s_1^2 s_2 s_1^2 s_2 
                \stackrel{\textrm{R3}}{=}
                s_1s_2 {\color{teal}{s_3}} s_2 s_1s_2 s_1^2 s_2 s_1^2 s_2 
                \stackrel{\textrm{R1}}{=}
                s_1s_2^2 s_1s_2 s_1^2 s_2 s_1^2 s_2 \\
                &= s_1s_2^2 s_1s_2 s_1 s_1 s_2 s_1 {\color{teal}{s_1 s_2}} 
                \stackrel{\rho}{=} 
                s_1 s_2 s_1s_2^2 s_1 s_2 s_1 {\color{blue}{s_1 s_2 s_1}}  
                \stackrel{\textrm{R3}}{=} 
                s_1 {\color{red}{s_2}} s_1s_2^2 s_1 {\color{red}{s_2}} s_1 s_2 {\color{red}{s_1}} s_2
                {\succ}
                s_1^2s_2^2s_1^2s_2^2.
            \end{align*}
            
    \vskip 2mm
\paragraph{{\it Case 2: k=1}} We assume that $w_1 =s_1s_3$ after a necessary cyclic rotation. Then $w_2,w_3,w_4$ are either $s_1^2$ or $s_3^2$. Note that $w_2\neq w_4$.  By the symmetry between $s_1$ and $s_3$, and taking rotations and the opposite word if necessary, it suffices to consider $w_2=w_3=s_1^2$ and $w_4=s_3^2$. The $\Lambda_\beta$ is admissible concordance to the standard $\mathrm{E}_9$ link:
        \begin{align*}
            \beta 
            &=
            s_3s_1s_2s_1^2s_2s_1^2{\color{teal}{s_2s_3^2s_2}}
            \stackrel{\rho}{=}
            s_2{\color{blue}{s_3^2s_2s_3}}s_1s_2s_1^2s_2s_1^2
            \stackrel{\textrm{R3}}{=}
            s_2s_2^2{\color{teal}{s_3}}s_2s_1s_2s_1^2s_2s_1^2 \\
            &\stackrel{\textrm{R1}}{=}
            s_2s_2^2s_2s_1s_2s_1^2s_2s_1^2
            =
            {\color{blue}{s_2^4s_1s_2}}s_1^2s_2s_1^2
            \stackrel{\textrm{R3}}{=}
            s_1^4s_2s_1s_1^2s_2{\color{teal}{s_1^2}}
            \stackrel{\rho}{=}
            s_1^6s_2s_3s_2.
        \end{align*}
  \vskip 2mm      
 \paragraph{{\it Case 3: k=0}} Assume that $w_1=w_2 =s_1^2$ and $w_3=w_4 = s_3^2$. Then $Q_\beta$ is of type $\tilde{\mathrm{D}}_8$:
    \[
    \begin{tikzpicture}
	\draw [thick,->](0.15,0) -- (0.85, 0);
	\filldraw (0,0) circle (2pt);
	\filldraw (1,0) circle (2pt);
	\filldraw (-1,0) circle (2pt);
	\filldraw (0,-1) circle (2pt);
	\filldraw (1,-1) circle (2pt);
	\filldraw (2,-1) circle (2pt);
		\filldraw (3,-2) circle (2pt);
	\filldraw (2,-2) circle (2pt);
	\filldraw (1,-2) circle (2pt);
    \draw [thick,->] (0.15,-1) -- (0.85, -1);
	\draw [thick,<-](0,-0.15) -- (0, -0.85);
	\draw [thick,<-](-0.15,0) -- (-0.85, 0);
	\draw [thick,->](1.15,-1) -- (1.85,-1);
	\draw [thick,->](1.15,-2) -- (1.85,-2);
	\draw [thick,->](2.15,-2) -- (2.85,-2);
	\draw [thick,<-](2,-1.15) -- (2, -1.85);
	\end{tikzpicture}. \qedhere
    \]
\end{proof}

\begin{defn}\label{defn:subword}
Let $\beta$ be a braid word of $n$ strands. For $1\leq i < j\leq n-1$, we define
$$\beta(i,j) := \textrm{the sub-word of }{\beta} \textrm{ that contains } s_i, s_{i+1},\dotsb, s_j.$$
For example, if $\beta = s_1 s_2s_3 s_1^2 s_2 s_5s_2s_3s_4$, then $\beta(2,3) = s_2s_3s_2^2s_3$.
\end{defn}

\begin{lem}\label{quiverdeg} Let $\beta$ be a braid word of $n$ strands.
\begin{enumerate}[label=\emph{(\arabic*)}]
    \item If $s_i^2$ is not a subword of $\beta$, then $Q_\beta=Q_{\beta(1,i-1)}\sqcup Q_{\beta(i+1,n-1)}$.

    \item If $\beta(i,i+1)$ does not contain a sub-word of intertwining pairs, namely neither $s_{i}s_{i+1}s_{i}s_{i+1}$ nor $s_{i+1}s_{i}s_{i+1}s_{i}$, then $Q_\beta=Q_{\beta(1,i)}\sqcup Q_{\beta(i+1,n-1)}$.
\end{enumerate}
\end{lem}
\begin{proof} The brick diagram has an empty level $i$ in case (1) and does not have arrows between level $i$ and level $i+1$ in case (2).
\end{proof}

\begin{lem} Let $n\geq 3$ and let $\beta$ be an $n$-strand braid word such that $Q_\beta$ is connected.  If $\beta\succ s_1^2$ and  $\beta \succ s_{n-1}^2$, then $Q_\beta$ is acyclic if and only if for $1\leq i\leq n-2$, we have \[\beta(i,i+1)=s_i^{a_1}s_{i+1}^{b_1}s_i^{a_2}s_{i+1}^{b_2} ~\mbox{or}~ s_{i+1}^{b_1}s_i^{a_1}s_{i+1}^{b_2}s_i^{a_2}, \quad \quad \mbox{where } a_1, a_2, b_1, b_2\geq 1\]
\end{lem}
\begin{proof} The \emph{if} direction is obvious. To see the \emph{only if} direction, let us assume without loss of generality that $\beta(i,i+1)$ begins with $s_i$. If $\beta(i,i+1)$ does not end after $s_i^{a_1}s_{i+1}^{b_1}s_i^{a_2}s_{i+1}^{b_2}$, then there is at least one $s_i$ after $s_{i+1}^{b_2}$, giving $Q_\beta$ an $a_2$-cycle between levels $i$ and $i+1$.
\end{proof}

\begin{assumption}\label{mainassumption}
Note that $2$-strand braids correspond to type A quivers. It suffices to consider  braid words $\beta$ of at least 3 strands. Let us single out the generator $s_2$. After necessary rotations, we assume that $\beta$ does not start with $s_2$ but ends with $s_2$, that is, 
$$\beta = w_1s_2^{b_1} w_2 s_2^{b_2}\dotsb w_m s_2^{b_m},$$
where each $w_i$ is a word of $s_1, s_3, s_4, \dotsb, s_{n-1}$. 

We assume that every $w_i$ contains at least one $s_1$ or $s_3$; otherwise, we can move the whole $w_i$ across the $s_2$'s at either end and merge it with $w_{i-1}$ or $w_{i+1}$. We further assume that $\sum b_i$ achieves minimum. Under this assumption, the length of every $w_i(1,3)$ is at least  $2$. Otherwise, with the letters $s_4\dotsb, s_{n-1}$ migrated away, we have $s_2w_is_2=s_2s_1s_2$ or $s_2s_3s_2$, and we can use R3 to reduce $\sum b_i$. 

We assume that $m\geq 2$; otherwise, $Q_\beta$ is disconnected by Lemma \ref{quiverdeg}. Meanwhile, if $m\geq 4$, then after necessarily deleting letters, we land on the case of Lemma \ref{ws4lemma}, and the braid $\beta$ dominates a braid with an acyclic quiver of infinite type.

In the rest of this section, without loss of generality, we assume that 
\begin{equation}
\label{beta.exp}
\beta = w_1s_2^{b_1}w_2s_2^{b_2}\dotsb w_ms_2^{b_m},
\end{equation}where  $b_i \geq 1$, $m=2$ or $3$, and $w_i\succeq s_1^2, s_3^2$ or $s_1s_3$.
\end{assumption}
We prove Proposition \ref{Mainquiver} by induction on the number of strands of $\beta$.

\medskip

\noindent\textbf{(II) Proof of Proposition \ref{Mainquiver} for 3-strand braids{}}
 If $m=2$ in \eqref{beta.exp}, then $Q_\beta$ is acyclic and therefore the proposition follows. It remains to consider $m=3$.  Suppose that at least one of the $b_i$'s, say $b_3$ after necessary cyclic rotations, is greater than 1.  The proposition follows since 
  \[\beta \succeq w_1s_2{\color{red}{w_2}}s_2w_3s_2^2
    {\succ} 
    w_1s_2^2w_3s_2^2\succeq s_1^2s_2^2s_1^2s_2^2.\]
It remains to consider $b_1=b_2=b_3=1$, i.e., 
\[\beta = s_1^{a_1}s_2s_1^{a_2}s_2s_1^{a_3}s_2.\] 

If two of $a_i$'s, say $a_1$ and $a_2$ after necessary rotations, are equal to 2, then
        \begin{align*}
            \beta &= s_1{\color{blue}{s_1s_2s_1}}s_1s_2s_1^{a_3}s_2 
                \stackrel{\text{R3}}{=} s_1s_2s_1{\color{blue}{s_2s_1s_2}}s_1^{a_3}s_2 
                \stackrel{\text{R3}}{=} s_1s_2s_1s_1s_2s_1s_1^{a_3}{\color{teal}{s_2}} \\
                &\stackrel{\rho}{=} {\color{blue}{s_2s_1s_2}}s_1s_1s_2s_1s_1^{a_3} 
                \stackrel{{\text{R3}}}{=} {\color{teal}{s_1}}s_2s_1s_1s_1s_2s_1s_1^{a_3}
                \stackrel{{\rho}}{=} s_2s_1s_1s_1s_2s_1s_1^{a_3}s_1 = s_2s_1^3s_2s_1^{a_3+2}        
        \end{align*}
The quiver for the last word is  acyclic. The proposition is proved.

Otherwise, at least two of the $a_i$'s, say $a_1$ and $a_2$ after necessary rotations, are greater than 2. The proposition follows since 
        \begin{align*}
        \beta 
            &\succ s_1^{3}s_2s_1^{3}s_2s_1^{2}s_2 
            = s_1^{2}{\color{blue}{s_1s_2s_1}}s_1^{2}s_2s_1^{2}s_2 
            \stackrel{\text{R3}}{=} 
                s_1^{2}s_2{\color{red}{s_1}}s_2s_1^{2}s_2{\color{red}{s_1^{2}}}s_2
            {\succ} s_1^2s_2^2s_1^2s_2^2.
        \end{align*}

\medskip

\noindent\textbf{(III) Proof of Proposition \ref{Mainquiver} for braids of at least 4 strands.}
Assume that $\beta$ is expressed as in \eqref{beta.exp}. Note that $s_1$ commutes with all other generators in $w_i$. Therefore we further assume that 
\begin{itemize}
    \item $w_i = s_1^{a_i}v_i = v_i s_1^{a_i}$, where $v_i$ is a word of $s_3, \dotsb, s_{n-1}$.
\end{itemize}
We shall start with the proof of the following two lemmas. 

\begin{lem}\label{3.13} Suppose $\beta \succ s_1$ and $\beta \succ s_2$. If $Q_{\beta(1,2)}$ is of Dynkin type $\mathrm{A}$, then there exists an admissible concordance $\Lambda_\gamma\rightarrow \Lambda_\beta$ such that $\gamma$ has fewer strands than $\beta$.
\end{lem}
\begin{proof} Since $Q_{\beta(1,2)}$ is of Dynkin type $\mathrm{A}$, $\beta(1,2)$ must be of the form $s_1^{a_1}s_2^{b_1}s_1^{a_2}s_2^{b_2}$ with $\min\{a_1, a_2\}=\min\{b_1,b_2\}=1$. After necessary cyclic rotations and/or taking the opposite word, we assume $a_1=b_1=1$. Then
\[
\beta = v_1 {\color{blue}s_1 s_2 s_1^{a_2}} v_2 s_2^{b_2} \stackrel{\textrm{R3}}{=} v_1 s_2^{a_2} {\color{teal}s_1} s_2 v_2 s_2^{b_2}   \stackrel{\textrm{R1}}{=}v_1 s_2^{a_2} s_2 v_2 s_2^{b_2}.
\]
The braid reduces to the case of one fewer strand. 
\end{proof}

\begin{lem} \label{A sub-tree} Suppose $Q_\beta$ is connected. If $Q_{\beta(1,3)}$ is acyclic and $Q_{\beta(1,2)}$ is not of type $\mathrm{A}$, then Proposition \ref{Mainquiver} is true for $[\beta]$.
\end{lem}
\begin{proof} Define 
\[k: = \max\{i ~|~  Q_{\beta(1,i)} \mbox{ is acyclic}\}.
\]
If $k=n$, then $Q_\beta$ is acyclic and the Lemma is proved. If $Q_{\beta(1,k)}$ is of infinite type, then the Lemma follows since  $\beta \succ \beta(1,k)$. Now we assume $k<n$ and $Q_{\beta(1,k)}$ is of finite type.

Note that $Q_{\beta(1,3)}$ is a subquiver of $Q_{\beta(1,k)}$. By assumption, $Q_{\beta(1,2)}$ is not of type $\mathrm{A}$. Therefore $Q_{\beta(1,k)}$ must be of type $\mathrm{D}$ or $\mathrm{E}$. Hence, $Q_{\beta(i,j)}$ is of type $\mathrm{A}$ for $1<i<j\leq k$.
In particular, $Q_{\beta(k-1,k)}$ is of type $\mathrm{A}$ and $Q_\beta$ is connected. Therefore we have
$$\beta(k-1,k) = s_{k-1}^{e_1} s_{k}^{f_1} s_{k-1}^{e_2} s_{k}^{f_2},\quad \text{or} \quad \beta(k-1,k) = s_{k}^{f_1} s_{k-1}^{e_1} s_{k}^{f_2}s_{k-1}^{e_2},$$
where
$$\min\{e_1,e_2\} = \min\{f_1,f_2\}=1.$$

Below we consider the first case $\beta(k-1,k) = s_{k-1}^{e_1} s_{k}^{f_1} s_{k-1}^{e_2} s_{k}^{f_2}$. The second case follows by taking the opposite word of $\beta$. The letters $s_1,\dotsb, s_{k-2}$ commute with  $s_{k+1},\dotsb, s_n$. After necessary communications of the letters in $\beta$, we can write
$$\beta = \gamma_1\delta_1\gamma_2\delta_2,$$
where $\gamma_i~ (i=1,2)$ is a word of  $s_1\dotsb, s_{k-1}$ with $e_i$ many of $s_{k-1}$, and $\delta_i$ is a word of $s_k,\dotsb, s_n$ with $f_i$ many of $s_{k}$. We remark that we have \emph{not} performed cyclic rotations yet and will only do it carefully, so that the quiver for $\beta(1,k-1)=\gamma_1\gamma_2$ is not distorted. 

Recall that $\min\{f_1,f_2\} = 1$. We consider the case $f_1=1$. The argument for $f_2=1$ is a similar repetition. Let us write 
$\delta_1 = xs_k y,$
where $x,y$ are words of $s_{k+1}\dotsb, s_n$ and they commute with $\gamma_1,\gamma_2$. We pass $y$ through $\gamma_2$, and we pass $x$ through $\gamma_1$ and rotation, obtaining $\gamma_1(xs_ky)\gamma_2\delta_2 \rightsquigarrow \gamma_1 s_k \gamma_2 (y\delta_2x)$. This move does not change the quiver for $\beta (1,k)$, and is a Legendrian isotopy. Consequently, we can assume $\delta_1 = s_k$  and write $\beta=\gamma_1s_k\gamma_2\delta_2$.

Now we consider 
\[\delta_2(k,k+1)=s_k^{g_1}s_{k+1}^{h_1}s_k^{g_2}s_{k+1}^{h_2}\cdots s_k^{g_l},\] 
where $g_1, g_l\geq 0$ and all other powers $\geq 1$.
The Lemma holds for the following two cases. 
\begin{enumerate}
    \item If $\delta_2\succ s_{k+1}s_k^2s_{k+1}$, then $\beta \succ \gamma_1 s_k \gamma_2s_{k+1}s_k^2s_{k+1}:=\beta_1$. 
    
    \item If $\delta_2 \succ s_{k+1}^2s_ks_{k+1}^2$, then $\beta \succ \gamma_1 s_k \gamma_2s_{k+1}^2s_ks_{k+1}^2:=\beta_2$. 
\end{enumerate}
The quivers for $\beta_1$ and $\beta_2$ are acyclic and of infinite type, as depicted below: 
$$
\begin{tikzpicture}
    \filldraw (0,-1) circle (2pt);
    \filldraw (-1,0) circle (2pt);
    \node at (-1.5,0) [] {$\dotsb$};
    \filldraw (0,0) circle (2pt);
    \filldraw (1,0) circle (2pt);
    \node at (1.5,0) [] {$\dotsb$};
    \draw [thick,->](0.15,0) -- (0.85, 0);
    \draw [thick,<-](-0.15,0) -- (-0.85, 0);
    \draw [thick,<-](0,-0.15) -- (0, -0.85);
    \node at (0, -1.3) [] {\vdots};
    \filldraw (0,-1.8) circle (2pt);
    \filldraw (1,-1.8) circle (2pt);
    \filldraw (0,-2.8) circle (2pt);
    \draw [thick,->](0.15,-1.8) -- (0.85, -1.8);
    \draw [thick,->](0,-2.65) -- (0, -1.95);
    \node at (-2, -2.3) [] {$(1)$};
\end{tikzpicture}
\qquad\quad
\begin{tikzpicture}
    \filldraw (0,-1) circle (2pt);
    \filldraw (-1,0) circle (2pt);
    \node at (-1.5,0) [] {$\dotsb$};
    \filldraw (0,0) circle (2pt);
    \filldraw (1,0) circle (2pt);
    \node at (1.5,0) [] {$\dotsb$};
    \draw [thick,->](0.15,0) -- (0.85, 0);
    \draw [thick,<-](-0.15,0) -- (-0.85, 0);
    \draw [thick,<-](0,-0.15) -- (0, -0.85);
    \node at (0, -1.3) [] {\vdots};
    \filldraw (0,-1.8) circle (2pt);
    \filldraw (1,-2.8) circle (2pt);
    \filldraw (-1,-2.8) circle (2pt);
    \filldraw (0,-2.8) circle (2pt);
    \draw [thick,->](0.15,-2.8) -- (0.85, -2.8);
    \draw [thick,->](-0.85,-2.8) -- (-0.15, -2.8);
    \draw [thick,->](0,-2.65) -- (0, -1.95);
    \node at (-2, -2.3) [] {$(2)$};
\end{tikzpicture}
$$

By the definition of $k$, the quiver for $\beta(k,k+1)$ is not acyclic. Therefore we have $\delta_2\succeq s_{k+1}s_ks_{k+1}s_k$. We assume that $\delta_2$ does not satisfy the above $(1)$ or $(2)$. Then
 \[\delta_2(k,k+1)=s_k^{g_1}s_{k+1}^{h_1}s_ks_{k+1}^{h_2}s_k^{g_3},\]  where $g_1\geq 0$, $g_3\geq 1$, and $\min\{h_1,h_2\} =1$. Depending on whether $h_1=1$ or $h_2=1$, we have the following two cases: 
\[
\beta(1,k+1)
=
\gamma_1s_k\gamma_2{\color{blue}s_k^{g_1}s_{k+1}s_k}s_{k+1}^{h_2}s_k^{g_3}
\overset{\text{R3}}{=}
\gamma_1s_k\gamma_2s_{k+1}{\color{teal}s_ks_{k+1}^{g_1+h_2}s_k^{g_3}}
\overset{\rho}{=}
s_ks_{k+1}^{g_1+h_2}s_k^{g_3}\gamma_1s_k\gamma_2s_{k+1},
\]
\[
\beta(1,k+1)
=
\gamma_1s_k\gamma_2s_k^{g_1}s_{k+1}^{h_1}{\color{blue}s_ks_{k+1}s_k^{g_3}}
\overset{\text{R3}}{=}
\gamma_1s_k\gamma_2s_k^{g_1}s_{k+1}^{h_1+g_3}s_ks_{k+1}.
\]
In both cases, the only R3 move is $s_ks_{k+1}s_k \rightsquigarrow s_{k+1}s_ks_{k+1}$ (only from left to right) performed in $\delta_2$, hence the move can be extended from $\beta(1,k+1)$ to $\beta$. The cyclic rotations can also be extended to $\beta$ without changing the quiver for $\beta(1,k)$. In the end, we performed a Legendrian isotopy and get a new braid word $\beta'$ with acyclic $Q_{\beta'(1,k+1)}$. We repeat the above argument for $\beta'$ and $k+1$.
This completes the case $f_1=1$. 
\end{proof}

Now we prove the proposition. If $m=2$ in \eqref{beta.exp}, then $Q_{\beta(1,3)}$ is acyclic. If $Q_{\beta(1,2)}$ is not of type $\mathrm{A}$, then the proposition follows directly from Lemma \ref{A sub-tree}. Otherwise, we apply Lemma \ref{3.13}. It remains to consider the case $m=3$, in which we have
\begin{equation}\beta =w_1s_2^{b_1}w_2s_2^{b_2}w_3s_2^{b_3}= v_1s_1^{a_1}s_2^{b_1}v_2s_1^{a_2}s_2^{b_2}v_3s_1^{a_3}s_2^{b_3}.
\end{equation}

Let us set
\[
p= \# \{i~|~a_i\neq 0 \}, \hskip 14mm q=\#\{i~|~ v_i \succeq s_3\}.
\] Here $p,q \in \{2, 3\}$. We consider cases by $(p,q)$.

\medskip
\paragraph{\it Case 1: $(p,q) = (2,2)$} 
After suitable cyclic rotation, we assume $a_3=0$. Then $Q_{\beta(1,3)}$ is acyclic. The rest goes through the same line as the above proof for the case $m=2$.

 \medskip   
 \paragraph{\it Case 2: $(p,q) = (2,3)$} After suitable cyclic rotation, we assume $a_3=0$. If $b_1 \geq 2$, then 
 \[\beta= v_1s_1^{a_1}s_2^{b_1}v_2s_1^{a_2}s_2^{b_2}v_3s_2^{b_3} \succ v_1s_1s_2^{b_1}v_2s_1s_2^{b_2+b_3} \succeq s_3s_1s_2^2 s_3s_1s_2^2.
 \]The proposition follows. So we assume $b_2=1$. 
 
 Now if $a_2 =1$, then using  $s_1^{a_1}s_2s_1 = s_2s_1s_2^{a_1}$, we can reduce the number of strands. The same argument works for $a_1=1$.
It remain to consider $a_1 \geq 2$ and $a_2\geq 2$. Then 
    $$\beta \succeq 
    s_1^2 {\color{red}{v_1}} s_2 {\color{red}{v_2}} s_1^2 s_2 {\color{red}{v_3}} s_2 
    {\succ} 
    s_1^2 {\color{blue}s_3 s_2 s_3} s_1^2s_2 s_2 
    \stackrel{\textrm{R3}}{=}
    s_1^2 s_2 {\color{red}{s_3}} s_2 s_1^2s_2 s_2
    {\succ} 
    s_1^2 s_2^2  s_1^2s_2^2.
    $$
 
 \medskip   
 \paragraph{\it Case 3: $(p,q) = (3,3)$} We have 
 \[w_i=v_is_i^{a_i}\succeq s_3s_1,  \quad \quad \forall i=1,2,3.\] 
 If there is a $b_i$, say $b_1$ after necessary cyclic rotation, greater than $1$, then 
 \[\beta \succ s_3s_1 s_2^{b_1}s_3s_1s_2^{b_2+b_3}\succeq s_1s_3s_2^2s_1s_3s_2^2.
 \]
The proposition follows. It remains to consider  $b_1= b_2=b_3 =1$.

If there is some $w_i$ with $w_i(1,3) = s_1s_3$, after suitable cyclic rotation we can assume  $w_2(1,3) = s_1s_3$. Note that the rest letters of $w_2$ are $s_4, \ldots, s_n$. They commute with the $s_2$ at either end and can be merged into $w_1$ or $w_3$. Therefore, we may assume $w_2=s_1s_3$ and  use  identity 
   \begin{equation}\label{3strandstocancel1}
        s_1^{a_1}s_2s_1s_3s_2s_1^{a_3} = s_3^{a_3}s_2s_1s_3s_2s_3^{a_1}
    \end{equation}
to reduce the number of strands.
    
    If none of the $w_i$'s has $w_i(1,3) = s_1s_3$. Then $w_i(1,3) \succeq s_1^2s_3$ or $w_i(1,3) \succeq s_1s_3^3$ for $i=1,2,3$. Two of them must be the same kind, and they have adjacent indices after cyclic rotation. For example, if $w_1,w_2$ are of the same type $s_1s_1s_3$, then
    $$\beta \succeq 
    s_1s_3s_1 s_2 s_1s_3s_1 s_2 {\color{red}{w_3}} s_2
    {\succ} 
    s_1s_3 {\color{blue}{s_1s_2s_1}} s_3s_1s_2s_2
    \stackrel{\textrm{R3}}{=}
    s_1s_3 s_2 {\color{red}{s_1}} s_2 s_3s_1s_2s_2
    {\succ} 
    s_1s_3 s_2^2 s_3s_1s_2^2.
    $$
    Other combinations of $w_i$ are similar.

\medskip
\paragraph{\it Case 4: $(p,q) = (3,2)$} If $\beta$ is a $4$-strand word, then by the symmetry between $s_1$ and $s_3$, it reduces to the case $(p,q)=(2,3)$. Below we assume $\beta$ is at least of $5$ strands.

After cyclic rotations, we assume that $v_3$ does not contain  $s_3$. Then $v_3$ commutes with $s_2$ and can be merged into $w_2$. By \eqref{beta.exp}, we assume that  
\[w_1\succeq s_1s_3, \quad \quad w_2\succeq s_1s_3, \quad \quad w_3 = s_1^{a_3} \mbox{ with } a_3\geq 2.\] 
If $b_1\geq 2$, then the proposition follows since
\[\beta\succeq s_1s_3s_2^2s_1s_3s_2 {\color{red}{w_3}} s_2 
{\succ} 
s_1s_3s_2^2s_1s_3s_2^2.\]  Below we consider $b_1=1$.

If $w_1\succeq s_1s_3^2$ and $w_2\succeq s_1s_3^2$, then the proposition follows since
\[\beta\succeq 
s_1s_3 {\color{blue}{s_3s_2s_3}} s_3s_1s_2^2
\stackrel{\textrm{R3}}{\succ}
s_1s_3s_2{\color{red}{s_3}}s_2s_3s_1s_2^2
{\succ} 
s_1s_3s_2^2s_1s_3s_2^2.\] 
 Hence, we assume that one of $w_1, w_2$ contains a single $s_3$. After suitable cyclic rotations and taking the opposite word if necessary, we assume that $w_2$ contains a single $s_3$. Moreover, all the letters $s_4,\dotsb, s_{n-1}$ in $w_2$ can be merged to $w_1$ by moving them in two directions and taking necessary cyclic rotations.
To summarize, it remains to consider 
\[\beta = v_1 s_1^{a_1} s_2 s_1^{a_2} s_3 s_2^{b_2} s_1^{a_3} s_2^{b_3},\quad \mbox{where }v_1\succeq s_3,~a_3\geq 2,~ \mbox{and }a_1, a_2, b_2, b_3\geq 1.\]
 We split our proof into two cases based on the value of $a_2$.

\medskip

A. If $a_2 =1$, then $b_2\geq 2$. Otherwise, $\beta =v_1 {\color{purple}{s_1^{a_1} s_2 s_1 s_3 s_2 s_1^{a_3}}} s_2^{b_3}$, and we can apply Identity \eqref{3strandstocancel1} to the purple part to reduce the number of strands. We further assume $b_3 =1$; otherwise, $b_3\geq 2$, and together with $a_3,b_2\geq 2$, we have
    $$\beta \succeq s_1^{a_1}{\color{red}{s_2}}s_1^{a_2}s_2^{b_2}s_1^{a_3}s_2^{b_3}
    {\succ} 
    s_1^{a_1+a_2}s_2^{b_2}s_1^{a_3}s_2^{b_3}
    {\succeq}
    s_1^{2}s_2^{2}s_1^{2}s_2^{2}.
    $$
    To recollect, we have
    \[\beta = v_1 s_1^{a_1} s_2 s_1^{a_2} s_3 s_2^{b_2\geq 2} s_1^{a_3\geq 2} s_2.\] 
    If $w_1(1,3)=s_1s_3$, then $w_1 =x s_1s_3 y$. After rotating $s_1^{a_3} s_2$, and moving $x,y$, we have
    $$
    \beta = 
    x s_1s_3 y s_2 s_1^{a_2} s_3 s_2^{b_2} {\color{teal}{s_1^{a_3} s_2}}
    \stackrel{\rho}{=}
    s_1^{a_3} s_2 {\color{teal}{x}} s_1s_3 {\color{teal}{y}} s_2 s_1^{a_2} s_3 s_2^{b_2}
    \stackrel{\textrm{c},\rho}{=}
    {\color{purple}{s_1^{a_3} s_2 s_1s_3 s_2 s_1^{a_2}}} ys_3 s_2^{b_2}x.
    $$
    We apply identity \eqref{3strandstocancel1} to the purple part to reduce the number of strands. 
    Therefore we can assume $w_1(1,3)\succeq s_1s_3^2$ or $s_1^2s_3$. 
    
    Now we focus on $w_1(1,4)$. The connectedness of $Q_\beta$ implies that $w_1(1,4)$ has at least two copies of $s_4$, with at least one $s_3$ sandwiched in between. Hence there are four possibilities:
    $$
    w_1(1,4)\succeq \; 
    (a)\; s_1^2s_4s_3s_4, \;
    (b)\; s_1s_4s_3s_3s_4, \;
    (c)\; s_1s_4s_3s_4s_3, \;
    (d)\; s_1s_3s_4s_3s_4. 
    $$
The proposition follows via direct calculations:
    \begin{enumerate}
        \item[($a$)] $w_1(1,4)\succeq s_1^2s_4s_3s_4 =  s_1^2s_3s_4s_3\succ s_1^2s_3^2$. Then
        $$
        \quad \quad \quad \beta \succ s_1^2{\color{blue}{s_3^2s_2s_3}}s_1s_2^2s_1^2s_2
        \stackrel{\textrm{R3}^2}{=}
        s_1^2s_2s_3s_2{\color{blue}{s_2s_1s_2}}s_2s_1^2s_2
        \stackrel{\textrm{R3}}{=}
        s_1^2s_2{\color{red}{s_3}}s_2s_1{\color{red}{s_2}}s_1s_2{\color{red}{s_1^2}}s_2
        {\succ} 
        s_1^2s_2^2s_1^2s_2^2.
        $$
        
        \item[($b$)] $w_1(1,4)\succeq s_1s_4s_3s_3s_4$. Then
        \begin{align*}
        \beta 
        & \succeq s_1{\color{teal}{s_4}}s_3^2{\color{teal}{s_4}}s_2s_1s_3s_2^2s_1^2s_2
        \stackrel{\textrm{c},\rho}{=}
        s_1s_3^2s_2s_1{\color{blue}{s_4s_3s_4}}s_2^2s_1^2s_2
        \stackrel{\textrm{R3}}{=}
        s_1s_3^2s_2s_1s_3{\color{red}{s_4}}s_3s_2^2s_1^2s_2 \\
        &{\succ} 
        s_1s_3^2s_2s_1{\color{teal}{s_3}}s_3s_2^2s_1^2s_2
        \stackrel{\textrm{c}}{=}
        s_1s_3{\color{blue}{s_3s_2s_3}}s_1s_3s_2^2s_1^2s_2
        \stackrel{\textrm{R3}}{=}
        s_1s_3s_2{\color{red}{s_3}}s_2s_1s_3s_2^2{\color{red}{s_1^2s_2}}\\
        &{\succ}  s_1s_3s_2^2s_1s_3s_2^2.
        \end{align*}
        
        \item[($c$)]  $w_1(1,4)\succeq s_1s_4s_3s_4s_3 = s_1s_3s_4s_3s_3 \succ s_1s_3^3$. Then
        \begin{align*}
        \beta 
        & \succeq s_1{\color{blue}{s_3^3s_2s_3}}s_1s_2^2s_1^2s_2
        \stackrel{\textrm{R3}}{=}
        s_1s_2{\color{teal}{s_3}}s_2^3s_1s_2^2s_1^2s_2
        \stackrel{\textrm{R1}}{=}
        s_1s_2^4s_1s_2^2{\color{teal}{s_1^2s_2}} \\
        &\stackrel{\rho}{=}
        {\color{blue}{s_1^2}}s_2s_1s_2^4s_1s_2^2
        \stackrel{\textrm{R3}}{=}
        {\color{teal}{s_2}}s_1s_2^2s_2^4s_1s_2^2
        \stackrel{\rho}{=}
        s_1s_2^6s_1s_2^3.
        \end{align*}
        We end up with an $\mathrm{E}_9$ quiver, which is acyclic and of infinite type.
        
        \item[($d$)] $w_1(1,4)\succeq s_1s_3s_4s_3s_4 = s_1s_3s_3s_4s_3 \succ s_1s_3^3$. The rest follows from the same calculation as in $(c)$.
    \end{enumerate}

\medskip    

B. If $a_2\geq 2$, then we look at $w_1(1,3)$. 

If $w_1(1,3) = s_1s_3$, then $v_1=xs_3y$, where $x,y$ are words of $s_4, \ldots, s_{n-1}$.  Let $\tilde{x}$ and $\tilde{y}$ be the opposite word of $x$ and $y$ respectively. Then 
    \[
    \beta = {\color{teal}{x}}s_3y s_1 s_2 s_1^{a_2} s_3 {\color{teal}s_2^{b_2} s_1^{a_3} s_2^{b_3}} \stackrel{\rho, c}{=} s_2^{b_2} s_1^{a_3} s_2^{b_3}  s_3  s_1 s_2 s_1^{a_2} ys_3x \stackrel{\textrm{oppo}}{\rightsquigarrow} \tilde{x}s_3\tilde{y}s_1^{a_2}s_2s_1s_3s_2^{b_3}s_1^{a_3}s_2^{b_2}.
    \]
 It goes back to Case A. 
 
 If $w_1(1,3) \succeq s_1^2s_3$, then 
    $$
    \beta
    \succeq s_3s_1{\color{blue}{s_1s_2s_1}}s_1s_3s_2^{b_2}s_1^{a_3}s_2^{b_3}
    \stackrel{\textrm{R3}}{=}
    s_3s_1s_2{\color{red}{s_1}}s_2s_1s_3{\color{red}{s_2^{b_2}s_1^{a_3}s_2^{b_3}}}
    {\succ} s_1s_3s_2^2s_1s_3s_2^2.
    $$
    
It remains to consider $w_1(1,3)\succeq s_1s_3^2$.  There are three possibilities for $w_1(1,4)$:
    $$
    w_1(1,4)\succeq \; 
    (e)\; s_1s_4s_3s_4s_3, \;
    (f)\; s_1s_3s_4s_3s_4, \;
    (g)\; s_1s_4s_3s_3s_4. 
    $$
For both ($e$) and ($f$), after $s_4s_3s_4 = s_3s_4s_3$, we have $w_1(1,4)\succ s_1s_3^3$. Then 
        \begin{align*}
         \beta 
         & \succ s_1{\color{blue}{s_3^3s_2s_3}}s_1^2s_2s_1^2s_2
        \stackrel{\textrm{R3}}{=}
        s_1s_2{\color{teal}{s_3}}s_2^3s_1^2s_2s_1^2s_2
        \stackrel{\textrm{R1}}{=}
        s_1s_2^4s_1^2{\color{teal}{s_2s_1^2s_2}}  \\
        & \stackrel{\rho}{=}
        s_2s_1^2s_2s_1s_2^4s_1^2
        \stackrel{\textrm{R3}}{=}
        s_2s_1^2s_1^4s_2s_1s_1^2
        =
        s_2s_1^6s_2s_1^3
        \end{align*}
        This is again the $\mathrm{E}_9$ quiver.
For ($g$), we have $w_1(1,4)\succeq s_1s_4s_3s_3s_4$. Then
        \begin{align*}
        \beta 
        &\succeq
        s_1{\color{teal}{s_4}}s_3^2{\color{teal}{s_4}}s_2s_1^2s_3s_2s_1^2s_2
        \stackrel{\textrm{c},\rho}{=}
        s_1s_3^2s_2s_1^2{\color{blue}{s_4s_3s_4}}s_2s_1^2s_2
        \stackrel{\textrm{R3}}{=}
        s_1s_3^2s_2s_1^2s_3{\color{red}{s_4}}s_3s_2{\color{red}{s_1^2}}s_2 \\
        &{\succ} 
        s_1s_3^2s_2s_1^2{\color{teal}{s_3}}s_3s_2s_2
        \stackrel{\textrm{c}}{=}
        s_1s_3{\color{blue}{s_3s_2s_3}}s_1^2s_3s_2s_2
        \stackrel{\textrm{R3}}{=}
        s_1s_3s_2{\color{red}{s_3}}s_2{\color{red}{s_1}}s_1s_3s_2s_2 \\
        &{\succ} 
        s_1s_3s_2^2s_1s_3s_2^2.
        \end{align*}

This completes the proof of Proposition \ref{Mainquiver}.

\begin{cor} For positive braids $[\beta]$ with connected $Q_\beta$, the two cases in Proposition \ref{Mainquiver} coincides with the dichotomy between finite and infinite types for positive braids.
\end{cor}
\begin{proof} It follows from Proposition \ref{Mainquiver} and Proposition \ref{3.4}.
\end{proof}

\noindent \textit{Proof of Theorem \ref{MainStep1} for disconnected $Q_\beta$.} Suppose $Q_\beta$ has two components. Because vertices on the same level are connected, there exists a unique $1\leq i<n$ such that no arrow appears between level $i$ and $i+1$. We consider $\beta(1,i)$ and $\beta(i+1,n-1)$. Since we can pinch some crossings of $\beta$ to obtain $\beta(1,i)$ and $\beta(i+1,n-1)$, if one of them has infinitely many 
admissible fillings, so does $\beta$ by Proposition \ref{3.2}. Otherwise by Propositions \ref{Mainquiver} (1) and \ref{3.4}, both $Q_{\beta(1,i)}$ and $Q_{\beta(i+1,n-1)}$ are mutation equivalent to finite type quivers, and hence $[\beta]$ is of finite type. In general, we can induct on the number of components in the quiver of the braid. \qed

\bigskip

\subsection{Finite Type Classification}

In this section, we focus on positive braid Legendrian links  of finite type.

\begin{thm}\label{Dynkin quiver} Let $\beta$ be a braid word such that $Q_\beta$ is mutation equivalent to a Dynkin quiver and $\Lambda_\beta$ does not contain a split union of knots. Then $\Lambda_\beta$ is Legendrian isotopic to a standard link in Definition \ref{stad.links}.
\end{thm}
\begin{proof} By Proposition \ref{Mainquiver} (1), it suffices to assume that $Q_\beta$ is a  Dynkin quiver. 

If $Q_\beta$ is of type $\mathrm{A}$, we repeated utilize Lemma \ref{3.13} to reduce the number of strands of $\Lambda_\beta$ until it becomes $2$-strand link, which is a standard link of type A.

If $Q_\beta$ is of type $\mathrm{D}$ or $\mathrm{E}$, then it contains  a unique trivalent vertex.  If $n\geq 4$, we can apply Lemma \ref{3.13} to $\beta(1,2)$ or $\beta(n-2,n-1)$, whichever does not contain the trivalent vertex, to reduce $n$ until $n=3$. Note that $\beta$ can be written as \eqref{beta.exp}. Since $[\beta]$ is of finite type, following the discussion in (II) of Section \ref{sec inifnite fillings}, we may assume $m=2$ in \eqref{beta.exp}. After necessary rotation, we get
\[
\beta= s_1^{a_1}s_2^{b_1}s_1^{a_2}s_2^{b_2}, \quad \quad \mbox{where } a_1\geq 2, ~a_2 \geq 2, ~ \min\{b_1, b_2\}=1.
\]

 The trivalent vertex in a Dynkin $\mathrm{DE}$ quiver has three legs, at least one of which is of length $1$. For $Q_\beta$, two legs lie in level $1$ and one leg stretches to level $2$. We show that $b_1=b_2= 1$ after suitable Legendrian isotopy. Otherwise, one of the level $1$ legs is of length $1$. Then up to cyclic rotations, we get $a_2=2$. Depending on $b_1=1$ or $b_2=1$, we have the following Legendrian isotopies:
\begin{align*}
    \beta &= 
    s_1^{a_1}s_2s_1^2s_2^{b_2}
    =
    s_1^{{a_1}-1}{\color{blue}{s_1s_2s_1}}s_1s_2^{b_2}
    \stackrel{\textrm{R3}}{=}
    s_1^{a_1-1}s_2s_1{\color{blue}{s_2s_1s_2^{b_2}}}
    \stackrel{\textrm{R3}}{=}
    s_1^{a_1-1}s_2s_1^{b_2+1}s_2{\color{teal}{s_1}}
    \stackrel{\rho}{=}
    s_1^{a_1}s_2s_1^{b_2+1}s_2,    \\
    \beta &= 
    s_1^{a_1}s_2^{b_1}s_1^2s_2
    =
    {\color{teal}{s_1}}s_1^{{a_1}-1}s_2^{b_1}s_1^2s_2
    \stackrel{\rho}{=}
    s_1^{a_1-1}s_2^{b_1}s_1{\color{blue}{s_1s_2s_1}}
    \stackrel{\textrm{R3}}{=}
    s_1^{a_1-1}{\color{blue}{s_2^{b_1}s_1s_2}}s_1s_2
    \stackrel{\textrm{R3}}{=}
    s_1^{a_1}s_2s_1^{b_1+1}s_2.
\end{align*}
Eventually, after necessary cyclic notations, we get the standard links. 
\end{proof}

\begin{defn} Let $\beta$ be an $n$-strand braid word and let $\gamma$ be an $m$-strand braid word. Denote by $\gamma^{\#_j}$ the word obtained from $\gamma$ via $s_i\mapsto s_{i+j}$. 

The \emph{connect sum} of $\beta$ and $\gamma$ is the braid word $\beta \# \gamma:=\beta\gamma^{\#_{n-1}}$.

The  \emph{split union} of $\beta$ and $\gamma$ is the braid word $\beta \sqcup \gamma:=\beta\gamma^{\#_{n}}$.

Note that $\left[\beta \# \gamma\right] \in \Br_{n+m-1}^+$ and $[\beta\sqcup\gamma] \in \Br_{n+m}^+$.
\end{defn}

The connect sum of two positive braid links is again a positive braid link. By \cite{EV}, positive braid links attain a unique maximum tb Legendrian representative. The connect sum of two links is well-defined once specifying which components to attach the $1$-handle. Once well-defined, the connect sum is associative and commutative.

\begin{rmk} 
Here is  a list of the numbers of components for the standard $\mathrm{ADE}$ links:
\begin{center}
    \begin{tabular}{|c|c|c|} \hline
        \rule{0pt}{2.3ex} knots & 2-component links & 3-component links  \\[0.045cm] \hline
        \rule{0pt}{2.3ex} $\mathrm{A}_\text{even}$, $\mathrm{E}_6$, $\mathrm{E}_8$ & $\mathrm{A}_\text{odd}$, $\mathrm{D}_\text{odd}$, $\mathrm{E}_7$ & $\mathrm{D}_\text{even}$\\[0.045cm] \hline
    \end{tabular}.
\end{center}
\end{rmk}

\begin{thm}\label{4.6} If $[\beta]$ is of finite type, then $\Lambda_\beta$ is Legendrian isotopic to a split union of unknots and connect sums of standard links of ADE types.
\end{thm}

\begin{proof} For an $n$-strand positive braid $\beta$, the vertices of $Q_\beta$ are separated into $n-1$ many levels, each of which forms a type $\mathrm{A}$ quiver.
If $Q_\beta$ is disconnected, then 
\begin{enumerate}
\item two adjacent levels of $Q_\beta$ have vertices but no arrows in between; and/or  
\item a level of $Q_\beta$ has no vertex.
\end{enumerate}

For (1), after necessary rotation, we get  $\beta(i, i+1) = s_i^a s_{i+1}^b$ for some $i$. We may further commute $s_1,\dotsb, s_{i-1}$ with $s_{i+2},\dotsb, s_{n-1}$, obtaining \[\beta = \beta(1,i) \beta(i+1,n-1).\] Hence, $\beta$ is a connect sum of two braid words. 

For (2), we get $\beta(i,i) = s_i$ or empty  for some $i$. If it is empty, then 
\[\beta = \beta(1,i-1) \beta(i+1,n),\] which is a split union of two braid words. If  $\beta(i,i) = s_i$, then the braid is a connect sum via the following Legendrian isotopy:
\[
\begin{tikzpicture}[baseline=20,scale=0.9]
\draw [teal](0,-0.1) rectangle (0.5,0.35);
\draw (0,0.4) rectangle (0.5,0.85);
\draw (1,-0.1) rectangle (1.5,0.35);
\draw [teal](1,0.4) rectangle (1.5,0.85);
\draw (0.5,0.25) -- (1,0.5);
\draw (0.5,0.5) -- (1,0.25);
\draw (0.5,0.75) -- (1,0.75);
\draw (0.5,0) -- (1,0);
\draw (1.5,0.75) to [out=0,in=180] (1.75,0.875) to [out=180,in=0] (1.5,1)--(0,1) to [out=180,in=0] (-0.25,0.875) to [out=0,in=180] (0,0.75);
\draw (1.5,0.5) to [out=0,in=180] (2,0.875) to [out=180,in=0] (1.5,1.25) -- (0,1.25) to [out=180,in=0] (-0.5,0.875) to [out=0,in=180] (0,0.5);
\draw (1.5,0.25) to [out=0,in=180] (2.25,0.875) to [out=180,in=0] (1.5,1.5) -- (0,1.5) to [out=180,in=0] (-0.75,0.875) to [out=0,in=180] (0,0.25);
\draw (1.5,0) to [out=0,in=180] (2.5,0.875) to [out=180,in=0] (1.5,1.75) -- (0,1.75) to [out=180,in=0] (-1,0.875) to [out=0,in=180] (0,0);
\end{tikzpicture} \quad \overset{\rho}{\rightsquigarrow} \quad 
\begin{tikzpicture}[baseline=20,scale=0.9]
\draw [teal] (2,-0.1) -- (2.5,-0.1) -- (2.5,0.35) -- (2,0.35);
\draw (2,-0.1) -- (1.5,-0.1) -- (1.5,0.35) -- (2,0.35);
\draw (0.5,0.4) -- (1,0.4) -- (1,0.85) -- (0.5,0.85);
\draw [teal] (0.5,0.4) -- (0,0.4) -- (0,0.85) --  (0.5,0.85);
\draw (1,0.75) -- (2.5,0.75) to [out=0,in=180] (2.75,0.875) to [out=180,in=0] (2.5,1)--(0,1) to [out=180,in=0] (-0.25,0.875) to [out=0,in=180] (0,0.75);
\draw (0,0.25) -- (1,0.25) to [out=0,in=180] (1.5,0.5) -- (2.5,0.5) to [out=0,in=180] (3,0.875) to [out=180,in=0] (2.5,1.25) -- (0,1.25) to [out=180,in=0] (-0.5,0.875) to [out=0,in=180] (0,0.5);
\draw (1,0.5) to [out=0,in=180] (1.5,0.25);
\draw (2.5,0.25) to [out=0,in=180] (3.25,0.875) to [out=180,in=0] (2.5,1.5) -- (0,1.5) to [out=180,in=0] (-0.75,0.875) to [out=0,in=180]  (0,0.25);
\draw (2.5,0) to [out=0,in=180] (3.5,0.875) to [out=180,in=0] (2.5,1.75) -- (0,1.75) to [out=180,in=0] (-1,0.875) to [out=0,in=180] (0,0) -- (1.5,0);
\end{tikzpicture} \quad \rightsquigarrow\quad 
\begin{tikzpicture}[baseline=20,scale=0.9]
\draw (0,0.4) rectangle (0.5,0.85);
\draw (1,-0.1) rectangle (1.5,0.35);
\draw (0.5,0.75) to [out=0,in=180] (0.75,0.875) to [out=180,in=0] (0.5,1) -- (0,1) to [out=180,in=0] (-0.25,0.875) to [out=0,in=180] (0,0.75);
\draw (1.5,0) to [out=0,in=180] (2.5,0.875) to [out=180,in=0] (1.5,1.75) -- (-1,1.75) to [out=180,in=0] (-2,0.875) to [out=0,in=180] (-1,0) --  (1,0);
\draw [red] (0.5,0.25) to [out=0,in=180] (1,0.5) -- (1.5,0.5) to [out=0,in=180] (2,0.875) to [out=180,in=0] (1.5,1.25) -- (-1,1.25) to [out=180,in=0] (-1.5,0.875) to [out=0,in=180] (-1,0.5) -- (0,0.5);
\draw (0.5,0.5) to [out=0,in=180] (1,0.25);
\draw (1.5,0.25) to [out=0,in=180] (2.25,0.875) to [out=180,in=0] (1.5,1.5) -- (-1,1.5) to [out=180,in=0] (-1.75,0.875) to [out=0,in=180] (-1,0.25) -- (0.5,0.25);
\end{tikzpicture}
\]
\[
\overset{\text{R2,R3}}{\rightsquigarrow}\quad
\begin{tikzpicture}[baseline=20,scale=0.9]
\draw (0,0.4) rectangle (0.5,0.85);
\draw (1,-0.1) rectangle (1.5,0.35);
\draw (0.5,0.75) to [out=0,in=180] (0.75,0.875) to [out=180,in=0] (0.5,1) -- (0,1) to [out=180,in=0] (-0.25,0.875) to [out=0,in=180] (0,0.75);
\draw (1.5,0) to [out=0,in=180] (2.5,0.875) to [out=180,in=0] (1.5,1.75) -- (-1,1.75) to [out=180,in=0] (-2,0.875) to [out=0,in=180] (-1,0) -- (1,0);
\draw [red] (-1,0.25) to [out=0,in=180] (-0.25,0.625) to [out=180,in=0] (-0.625,0.75) to [out=180,in=0] (-1,0.625) to[out=0,in=180] (-0.5,0.5) -- (0,0.5);
\draw (0.5,0.5) to [out=0,in=180] (1,0.25);
\draw (1.5,0.25) to [out=0,in=180] (2.25,0.875) to [out=180,in=0] (1.5,1.5) -- (-1,1.5) to [out=180,in=0] (-1.75,0.875) to [out=0,in=180] (-1,0.25);
\end{tikzpicture}\quad 
\overset{\text{R1}}{\rightsquigarrow}\quad 
\begin{tikzpicture}[baseline=20,scale=0.9]
\draw (0,0.4) rectangle (0.5,0.85);
\draw (1,-0.1) rectangle (1.5,0.35);
\draw (0.5,0.5) to [out=0,in=180] (1,0.25);
\draw (0.5,0.75)  to [out=0,in=180] (0.75,0.875) to [out=180,in=0] (0.5,1) -- (0,1) to [out=180,in=0] (-0.25,0.875) to [out=0,in=180] (0,0.75);
\draw (1.5,0.25) to [out=0,in=180] (2,0.875) to [out=180,in=0] (1.5,1.25) -- (0,1.25) to [out=180,in=0] (-0.5,0.875) to [out=0,in=180] (0,0.5);
\draw (1.5,0) to [out=0,in=180] (2.25,0.875) to [out=180,in=0] (1.5,1.5) -- (0,1.5) to [out=180,in=0] (-0.75,0.875) to [out=0,in=180] (0,0) -- (1,0);
\end{tikzpicture}
\quad \rightsquigarrow\quad
\begin{tikzpicture}[baseline=20,scale=0.9]
\draw (0,0.4) rectangle (0.5,0.85);
\draw (1,0.15) rectangle (1.5,0.6);
\draw (0.5,0.5) -- (1,0.5);
\draw (0.5,0.75) -- (1.5,0.75) to [out=0,in=180] (1.75,0.875) to [out=180,in=0] (1.5,1) -- (0,1) to [out=180,in=0] (-0.25,0.875) to [out=0,in=180] (0,0.75);
\draw (1.5,0.5) to [out=0,in=180] (2,0.875) to [out=180,in=0] (1.5,1.25) -- (0,1.25) to [out=180,in=0] (-0.5,0.875) to [out=0,in=180] (0,0.5);
\draw (1.5,0.25) to [out=0,in=180] (2.25,0.875) to [out=180,in=0] (1.5,1.5) -- (0,1.5) to [out=180,in=0] (-0.75,0.875) to [out=0,in=180] (0,0.25)--(1,0.25);
\end{tikzpicture}
\]
Each quiver component is Legendrian isotopic to the standard $\mathrm{ADE}$ links. There could also be a split union of unknot for every pair of consecutive levels $\beta(i,i)$ and $\beta(i+1,i+1)$ that are both empty. This completes the proof.
\end{proof}

\bigskip

\setcounter{section}{0}

\renewcommand{\thesection}{\Alph{section}}

\section{Appendix I: Cluster varieties} 
\label{app.A}
We provide a rapid review on cluster varieties in the skew-symmetric cases. Below we set $[n]_+ := \max\{0,n\}$ for $n\in \mathbb{R}$.

\subsection{Definitions}
\label{append.A.1}
A {\it quiver} is a triple $Q=\left(I,I^\uf, \varepsilon\right)$, where  $I$ is a finite set, $I^\uf$ is a subset of  $I$, and $\varepsilon$ is an $I\times I$ skew-symmetric matrix whose entries $\varepsilon_{ij}$ are integers when $i\in I$ and $j\in I^\uf.$ 

Let $k\in I^\uf$. The \emph{mutation} in the direction $k$ produces a new  quiver $\mu_kQ = \left(I',I'^\uf,\varepsilon'\right)$ where $I'=I$, $I'^\uf=I^\uf$, and 
\[
\varepsilon'_{ij}=\left\{\begin{array}{ll}
    -\varepsilon_{ij} & \text{if $k\in \{i,j\}$},  \\
    \varepsilon_{ij} & \text{if $\varepsilon_{ik} \varepsilon_{kj}<0$, and $k\notin\{i,j\}$},\\
        \varepsilon_{ij}+ |\varepsilon_{ik}| \varepsilon_{kj}  & \text{if $\varepsilon_{ik} \varepsilon_{kj}\geq 0$, and $k\notin\{i,j\}$}.\\
\end{array}\right.
\]
Two quivers are  \emph{mutation equivalent} if they are related by a  sequence of mutations. Denote by $|Q|$ the class of quivers that are mutation equivalent to $Q$.

Each  $Q$ induces a directed graph with vertex set $I$. For $i, j\in I$, the number of arrows from $i$ to $j$ is $[\varepsilon_{ij}]_+$. 
Vertices parametrized by $i\in I- I^\uf$ are called \emph{frozen vertices}. In this paper, arrows among frozen vertices are allowed to be of half weight and will be illustrated by  dashed arrows. 

The \emph{unfrozen part} of $Q$ is the full subquiver $Q^\uf$ containing the unfrozen vertices.

A quiver $Q$ is said to be \emph{acyclic} if there is no directed cycle inside $Q^\uf$. 

A quiver $Q$ is said to be \emph{connected} if the underlying graph of $Q^\uf$ is connected.

A quiver $Q$ is said to have \emph{full-rank} if the submatrix $\varepsilon|_{I^\uf\times I}$ is of full-rank. 

Connectedness and being full-rank are invariant under mutations and therefore  descend to properties of mutation equivalence classes of quivers.

\begin{defn} \label{clusterK2}
A \emph{cluster $\mathrm{K}_2$ variety} $\mathscr{A}$ is an affine variety together with a collection $\mathcal{C}$ of triples $\alpha= (Q_\alpha, {\rm T}_\alpha, \vec{A}_\alpha),$
where
\begin{itemize}
    \item $Q_\alpha=\left(I, I^\uf, \varepsilon\right)$ is a quiver;
    \item ${\rm T}_\alpha$ is a split algebraic torus of rank $|I|$ inside $\mathscr{A}$;
    \item $\vec{A}_\alpha=\{A_{i;\alpha}\}_{i\in I}$ is a coordinate system of ${\rm T}_\alpha$.
\end{itemize}
We require that
\begin{itemize}

    \item For any unfrozen vertex $k$ of the quiver $Q_\alpha$, there is an  $\alpha'=(Q_{\alpha'}, {\rm T}_{\alpha'}, \vec{A}_{\alpha'}) \in \mathcal{C}$, where  $Q_{\alpha'}=\mu_kQ_{\alpha}$, and the transition map between $\vec{A}_{\alpha'}$ and  $\vec{A}_{\alpha}$ is
\[
    A_{i; \alpha'}=\left\{\begin{array}{ll}
        A_{k; \alpha}^{-1}\left(\prod_jA_{j;\alpha}^{[-\varepsilon_{kj}]_+}+\prod_jA_{j;\alpha}^{[\varepsilon_{kj}]_+}\right)
        & \text{if $i=k$},  \\
        A_{i;\alpha} & \text{if $i\neq k$}. 
    \end{array}\right.
\]
    We say that $\alpha'$ is  a \emph{mutation} of $\alpha$ in the direction  $k$ and write $\alpha'=\mu_k\alpha$. 
    \item Every pair $\alpha, \alpha' \in \mathcal{C}$ are related by a finite sequence of mutations. 
    \item The complement of the union of ${\rm T}_\alpha$ for all $\alpha$ is of codimension 2 in $\mathscr{A}$.
\end{itemize}
Each $\alpha$ is called a \emph{cluster seed}, ${\rm T}_\alpha$ is called a \emph{cluster chart}, $\vec{A}_\alpha$ is called a \emph{cluster}, and  $A_{i;\alpha}$ is called a \emph{cluster $\mathrm{K}_2$ coordinate} or a \emph{cluster variable}. 
Each $A_{i;\alpha}$ for $i\in I- I^\uf$ is invariant under mutations and is called a {\it frozen variable.}
We will suppress the subscript $\alpha$ when the context is clear. 
\end{defn}

\begin{rmk}
 The coordinate ring of a cluster chart ${\rm T}_\alpha$ is a Laurent polynomial ring $\mathbb{L}_{\alpha}$ in the variables $A_{i; \alpha}$. The intersection
$\bigcap_{\alpha \in \mathcal{C}} \mathbb{L}_\alpha$
is an {\it upper cluster algebra} of  \cite{BFZ}. A cluster ${\rm K}_2$ variety $\mathscr{A}$ is an affine variety whose coordinate ring is an upper cluster algebra.
It is worth mentioning that our cluster ${\rm K}_2$ varieties are different from the cluster $\mathcal{A}$ varieties in \cite{FGensemble}. The latter is defined to be the union of the tori ${\rm T}_\alpha$ for $\alpha\in \mathcal{C}$, and is not affine in general. 

 Each cluster seed $\alpha$ of $\mathscr{A}$ encodes a 2-form on ${\rm T}_\alpha$:
\begin{equation}
\label{2.form.cluster.A}
    \Omega_\alpha:= \sum \varepsilon_{ij} \frac{{\rm d} A_{i;\alpha}}{A_{i;\alpha}} \wedge \frac{{\rm d} A_{j;\alpha}}{A_{j;\alpha}}.
\end{equation}
By Corollary 6.5 of \cite{FGensemble}, this 2-form does not depend on the choice of cluster seeds and therefore defines a global 2-form $\Omega$ on $\mathscr{A}$.

 Borrowing ideas from mirror symmetry, Gross, Hacking, Keel, and Kontsevich  interpreted the cluster structures in terms of wall-crossing structures called scattering diagrams \cite{GHKK}. In detail, associated to a quiver $Q$ is a scattering diagram $\mathfrak{D}$. Inside $\mathfrak{D}$ is a simplicial fan consisting of cones called cluster chambers. The paper \cite{GHKK} establishes a one-to-one correspondence between the cluster seeds of $\mathscr{A}$ and the cluster chambers of $\mathfrak{D}$. The mutation from $\alpha$ to $\mu_k \alpha$ corresponds to crossing the sharing facet (a.k.a the {\it wall}) of their corresponding cluster chambers.

\end{rmk}

The following proposition is crucial for this paper. 

\begin{prop}\label{distinguishing charts} Let $Q$ be a quiver of full rank and let $\mathscr{A}$ be its associated cluster $\mathrm{K}_2$ variety over an algebraically closed field (of any characteristic). The cluster charts of distinct cluster seeds of $\mathscr{A}$ do not coincide. 
\end{prop}
\begin{rmk}
Proposition \ref{distinguishing charts} may not hold when $Q$ is not of full rank, e.g., if $Q$ contains one vertex and no arrows, then its cluster variety has two cluster seeds but only one cluster chart.
\end{rmk}
\begin{proof}
Let $\mathscr{A}$ be defined over an algebraically closed field of characteristic $p$. The characteristic 0 case follows by the same argument.
Let $\alpha$ and $\alpha'$ be two distinct cluster seeds of $\mathscr{A}$. By Corollary 6.3 of \cite{FZIV}, the transition map between $\vec{A}_{\alpha'}=\{A_{i}'\}$ and $\vec{A}_\alpha=\{A_i\}$ takes the form
\begin{equation}
\label{seperation formula}
    A_{i}'= \left(F_i{\Big|}_{X_k=\prod_lA_l^{\varepsilon_{kl}}}\right)\prod_{j\in I} A_{j}^{g_{ij}}, 
\end{equation}
where $g_{ij}$ are integers, and each $F_i$ is a polynomial in the variables $X_{k}$ for $k \in I^\uf$. 
The matrix $G=(g_{ij})$ is called a $g$-matrix. The polynomials $F_i$ are called F-polynomials.

By \cite{GHKK}, each $F_i$ is a generating function that counts broken lines in the scattering diagram associated to $Q$. For distinct $\alpha$ and $\alpha'$, there is at least one wall between their corresponding chambers. 
In particular, there is an $i\in I^\uf$ such that $F_i\neq 1$. 
By \cite{LS} and \cite{GHKK}, we  have
\begin{itemize}
    \item all coefficients of $F_i$ are positive integers;
    \item the constant term of $F_i$ is 1;
    \item the coefficient of the highest term of $F_i$ is 1. 
\end{itemize}
Here the highest term of $F_i$ is the monomial $\prod_j X_j^{a_j}$ such that for any other term $\prod_j X_j^{b_j}$ in $F_i$, we have $a_j\geq b_j$ for all $j$.
The above last two properties are equivalent due to \cite[Prop.5.3]{FZIV}.

By the above discussion, there exists an $i\in I^\uf$ such that the polynomial $F_i$ has at least two terms even after reducing to a polynomial with coefficients in the finite field $\mathbb{F}_p$.
The quiver $Q$ is of full rank. The substitution $X_k=\prod_lA_l^{\varepsilon_{kl}}$ gives rise to an injective homomorphism from  the polynomial ring $\mathbb{F}_p[X_i]_{i\in I^\uf}$ to the Laurent polynomial ring $\mathbb{F}_p[A_j^{\pm 1}]_{j\in I}$.  Therefore  $A'_i$ is not a Laurent monomial of $A_j$ for $j\in I$. On the other hand, biregular isomorphisms between algebraic tori over an algebraically closed field are of monomial coordinate transformations. Thus ${\rm T}_\alpha \neq {\rm T}_{\alpha'}$.
\end{proof}

\subsection{Cluster ensembles}

Following \cite{FGensemble}, {\it cluster Poisson varieties} are the cluster dual of cluster ${\rm K}_2$ varieties. Each cluster Poisson variety $\mathscr{X}$ is equipped with a collection of torus charts with  coordinate systems $\vec{X}_{\alpha} =\{{X}_{i;\alpha}^{\pm 1}\}_{i\in I}$.
The transition map between $\vec{X}_{\alpha'}=\vec{X}_{\mu_k\alpha}$ and $\vec{X}_\alpha$ is
\[
X_{i;\alpha'}=\left\{\begin{array}{ll}
    X_{k;\alpha}^{-1} & \text{if $i=k$,}  \\
    X_{i; \alpha} X_{k; \alpha}^{\left[\varepsilon_{ik}\right]_+}\left(1+X_{k;\alpha}\right)^{-\varepsilon_{ik}} & \text{if $i\neq k$.}
\end{array}\right.
\]
The coordinates $X_{i;\alpha}$ are called \emph{cluster Poisson coordinates}.

Let $\mathscr{A}$ and $\mathscr{X}$ be a pair of cluster varieties associated to a mutation equivalence class $|Q|$. There is a natural one-to-one correspondence between the cluster seeds of $\mathscr{A}$ and the cluster seeds of $\mathscr{X}$. Each pair of corresponding cluster seeds is called  a \emph{cluster seed} of  $(\mathscr{A}, \mathscr{X})$. Following \cite[\S 1.2]{FGensemble}, there is a canonical map $p:\mathscr{A}\rightarrow \mathscr{X}$ such that\footnote{Since $\varepsilon_{ij}$ may not be integers when $i,j$ are frozen, the map $p$ is not necessarily  algebraic. In Section \ref{acndjoaa}, we consider the unfrozen quotient $\mathscr{X}^\uf$ of $\mathscr{X}$. The induced map $p:\mathscr{A}\rightarrow\mathscr{X}^\uf$ is  algebraic.}
\[
p^*\left(X_{i;\alpha}\right)=\prod_j A_{j;\alpha}^{\varepsilon_{ij}}
\]
for every cluster seed of $(\mathscr{A}, \mathscr{X})$.
The triple $(\mathscr{A},\mathscr{X},p)$ is called a \emph{cluster ensemble}.

\begin{defn} Suppose $\sigma:I'\rightarrow I$ is an injective map such that
\begin{enumerate}
    \item $\sigma|_{I'^\uf}:I'^\uf\rightarrow I^\uf$ is a bijection,
    \item $\varepsilon'_{ij}=\varepsilon_{\sigma(i)\sigma(j)}$ for all $i,j\in I'$.
\end{enumerate}
Then $\sigma$ induces a morphism of algebraic tori $\sigma:\alpha'\rightarrow \alpha$ and $\sigma:\chi\rightarrow \chi'$, which are extended to morphisms of cluster varieties $\sigma:\mathscr{A}'\rightarrow \mathscr{A}$ and $\sigma:\mathscr{X}\rightarrow \mathscr{X}'$, called \emph{cluster morphisms}. If  $\sigma$ is bijective, then the induced cluster morphisms are called \emph{cluster isomorphisms}. 
\end{defn}

\begin{exmp} Consider the inclusion of the unfrozen part $Q^\uf=\left(I^\uf,I^\uf,\varepsilon|_{I^\uf\times I^\uf}\right)$ into $Q$. This inclusion induces cluster morphisms $\mathscr{A}^\uf\rightarrow \mathscr{A}$ and $\mathscr{X}\rightarrow \mathscr{X}^\uf$.
More properties about these cluster morphisms can be found in \cite[\S 3]{Sh}.
\end{exmp}

\begin{defn} A \emph{cluster automorphism} is a cluster isomorphism from a cluster variety to itself. Cluster automorphisms form a group $\mathcal{G}$ called the \emph{cluster modular group}.
\end{defn}

Fix an initial cluster seed of $(\mathscr{A}, \mathscr{X})$. Every cluster automorphism maps the initial seed to another cluster seed. We can express the obtained new cluster coordinates in terms of the initial ones as in \eqref{seperation formula}.
 As a consequence, one may assign the $c$-matrix, $g$-matrix, and $F$-polynomials of \cite{FZIV} to each cluster automorphism with respect to a fixed initial seed.

\begin{prop} A cluster automorphism $\sigma$ is the identity map  on $\mathscr{A}$ if and only if it is the identity map on $\mathscr{X}$.
\end{prop}
\begin{proof} The separation formula of Fomin-Zelevinsky \cite{FZIV} implies that $\sigma$ is the identity map on $\mathscr{A}$ (resp. $\mathscr{X}$) if and only if its $g$-matrix $G$ (resp. $c$-matrix $C$) is the identity matrix with respect to one (equivalently any) initial seed. The proposition then follows from the tropical duality theorem \cite[Theorem 1.2]{NZ}, which says that $C^{-1}=G^t$.
\end{proof}

\begin{defn}[\cite{GS2}]\label{def DT} A \emph{cluster Donaldson-Thomas transformation} on a cluster variety is a cluster automorphism whose $c$-matrix  is equal to minus identity. 
\end{defn}
For any cluster ensemble, its cluster Donaldson-Thomas transformation, if exists, is a unique central element in the cluster modular group.

\subsection{Quasi-cluster morphisms} \label{acndjoaa}
 Define $N:=\bigoplus_{i\in I}\mathbb{Z}e_i$ for a quiver $Q=\left(I,I^\uf,\varepsilon\right)$, and let $N^\uf$ its the sub-lattice spanned by $e_i$ for $i\in I^\uf$. The exchange matrix $\varepsilon$ equips $N$ with a  skew-symmetric form $\{\cdot, \cdot\}:N\times N\rightarrow \mathbb{Q}$ such that $\left\{e_i,e_j\right\}=\varepsilon_{ij}$. Let $M$ be the dual lattice of $N$. 
 
 One should think of $N$ as the character lattice of a cluster chart $\chi$ and think of $M$ as the character lattice of the cluster chart $\alpha$ dual to $\chi$. 
 For $n\in N$ and $m\in M$ we denote the corresponding character functions as $X^n$ and $A^m$. In particular, $X^{e_i}$ are precisely the cluster Poisson coordinates $X_i$, and the map $p:\mathscr{A}\rightarrow\mathscr{X}$ is induced by the linear map $p^*:N\rightarrow M, n\mapsto \{n,\cdot\}$. We will use this set-up to define quasi-cluster morphisms. More detailed discussions can be found in \cite{Fra,GS3,SWflag}.

\begin{defn}\label{def for quasi} Let $N$ and $N'$ be the lattices associated to $Q$ and $Q'$ respectively. Suppose $\sigma:N'\rightarrow N$ is an injective linear map such that
\begin{enumerate}
    \item $\sigma|_{N'^\uf}$ is an isomoprhism onto $N^\uf$;
    \item for any $i\in I'^\uf$, we have $\sigma\left(e'_i\right)=e_j$ for some $j\in I^\uf$,
    \item $\sigma$ preserves the skew-symmetric forms.
\end{enumerate}
Then $\sigma$ induces a morphism of algebraic tori $\sigma:\chi\rightarrow \chi'$ which extends to a morphism $\sigma:\mathscr{X}\rightarrow \mathscr{X}'$. On the dual side, $\sigma$ induces a linear map $\sigma:M\rightarrow M'$, which defines a morphism of algebraic tori $\sigma:\alpha'\rightarrow \alpha$ and extends to a morphism $\sigma:\mathscr{A}'\rightarrow \mathscr{A}$. 
We call the induced morphisms $\sigma:\mathscr{X}\rightarrow \mathscr{X}'$ and $\sigma:\mathscr{A}'\rightarrow \mathscr{A}$ \emph{quasi-cluster morphisms}.

A \emph{quasi-cluster isomorphism} is a quasi-cluster morphism where $\sigma: N'\rightarrow N$ is an isomorphism. 
A \emph{quasi-cluster automorphism} is a quasi-cluster isomorphism from a cluster variety to itself. 
Quasi-cluster automorphisms form a group $\mathcal{QG}$ called the \emph{quasi-cluster modular group}.
\end{defn}

The cluster modular group $\mathcal{G}$ is a subgroup of the quasi-cluster modular group $\mathcal{QG}$. There is a natural map $\mathcal{QG}\rightarrow \mathcal{G}^\uf$ where $\mathcal{G}^\uf$ denotes the cluster modular group for the unfrozen part.

\begin{rmk} Quasi-cluster automorphisms are also known as \emph{(quasi-)cluster transformations}.
\end{rmk}

The restriction of quasi-cluster morphisms to cluster charts commute with cluster mutations. Consequently, we have the following theorem.

\begin{thm} Let $\mathcal{V}$ and $\mathcal{W}$ be two cluster varieties of the same type (either $\mathrm{K}_2$ or Poisson). If $\sigma:\mathcal{V}\rightarrow \mathcal{W}$ is a quasi-cluster morphism, then there is a one-to-one correspondence between their cluster seeds, and $\sigma$ restricts to a morphism between the corresponding cluster charts. 
\end{thm}

Below we construct two types of quasi-cluster morphisms that are crucial for us.

\vskip 2mm

\paragraph{\bf Changing a Frozen Vertex.}
Recall the lattice $N$ associated with a quiver $Q=\left(I,I^\uf,\varepsilon\right)$.
 Let $k$ be a frozen vertex of $Q$. Let $\left(\delta_j\right)_{j\in I}$ is an $|I|$-tuple of integers.  We consider a lattice $N'=\bigoplus_{i\in I} \bbZ e_i'$ and define a linear map
 $\sigma: N'\rightarrow N$ such that 
\[
\sigma\left(e'_i\right):=\left\{\begin{array}{ll}
    e_i &  \text{if $i\neq k$}, \\
    \sum_{j\in I} \delta_je_j & \text{if $i=k$}.
\end{array}\right.
\]

 The exchange matrix $\varepsilon$ equips $N$ with a skew-symmetric form $\{\cdot, \cdot\}$, whose pull-back through $\sigma$ induces a skew-symmetric form $\{\cdot, \cdot\}'$ on $N'$. Let $\varepsilon'$ be an $I\times I$ matrix such that
\[
 \varepsilon'_{ij}=\left\{ e'_i, e'_j\right\}':=\left\{ \sigma(e'_i), \sigma(e'_j)\right\}.
\]
 Let $\mathscr{A}'$ and $\mathscr{X}'$ be the cluster varieties associated with the quiver $Q'=\left(I, I^\uf, \varepsilon'\right)$. Note that $\sigma$ satisfies the conditions (1) and (2) of Definition \ref{def for quasi}. Therefore it defines quasi-cluster morphisms 
\[
\sigma:\mathscr{A}'\rightarrow \mathscr{A} \quad \quad \text{and} \quad \quad \sigma:\mathscr{X}\rightarrow \mathscr{X}'.
\]
 Let $\alpha$ (resp. $\alpha'$) be the $\mathrm{K}_2$ cluster chart associated with the quiver $Q$ (resp. $Q'$). Let $\chi$ (resp. $\chi'$) be the Poisson cluster chart associated with the quiver $Q$ (resp. $Q'$). Then the pull-back maps of $\sigma$ can be written in terms of these cluster charts as
\begin{equation} \label{quasi open embedding}
\sigma^*\left(A_i\right)=\left\{\begin{array}{ll}
    A'_iA'^{\delta_i}_k & \text{if $i\neq k$}, \\
    A'_k & \text{if $i=k$}.
\end{array}\right. \quad \quad \text{and} \quad \quad 
\sigma^*\left(X'_i\right)=\left\{\begin{array}{ll}
    X_i & \text{if $i\neq k$}, \\
    \prod_j X_j^{\delta_j} & \text{if $i=k$}. 
\end{array}\right.
\end{equation}

\begin{prop}\label{changing frozen vertex} With the above set-up, the following statements are true.
\begin{enumerate}
    \item If $\delta_k=1$, then the quasi-cluster morphisms $\sigma$ are quasi-cluster isomorphisms.
    \item If $\sum_j \varepsilon_{ij}\delta_j=0$ for every $i\in I^\uf$, then there is no arrow between the vertex $k$ and the unfrozen part of $Q'$.
\end{enumerate}  
\end{prop} 
\begin{proof} (1) is obvious. For (2), it suffices to note that for $i \in I^\uf$, 
\[
\varepsilon_{ik}'=\left\{e'_i,e'_k\right\}'=\left\{\sigma\left(e'_i\right),\sigma\left(e'_k\right)\right\}
=
\bigg\{e_i,\sum_j \delta_je_j\bigg\}
=
\sum_j \delta_j\varepsilon_{ij}
=
0.
\]
Hence, there is no arrow between the vertex $k$ and the unfrozen part of $Q'$.
\end{proof}

\paragraph{\bf Merging Frozen Vertices.}

Let $t_1$ and $t_2$ be frozen vertices in a quiver $Q=\left(I,I^\uf,\varepsilon\right)$. Define the quiver $Q'=\left(I',I'^\uf,\varepsilon'\right)$, where $I':=(I\setminus \left\{t_1,t_2\right\})\sqcup\{t\}$, $I'^\uf:=I^\uf$, and 
\[
\varepsilon'_{ij}:=\left\{\begin{array}{ll}
    \varepsilon_{ij} & \text{if $i,j\neq t$}, \\
    \varepsilon_{t_1j}+\varepsilon_{t_2j} & \text{if $i=t$},\\
    \varepsilon_{it_1}+\varepsilon_{it_2} & \text{if $j=t$}.
\end{array}\right.
\]
We say that $Q'$ is obtained from $Q$ by \emph{merging} the frozen vertices $t_1$ and $t_2$ into a single frozen vertex $t$. Let $N$ and $N'$ be the lattices associated with the quivers $Q$ and $Q'$ respectively. There is an injective linear map
\begin{align*}
    \sigma:N'&\rightarrow N\\
    e'_i & \mapsto e_i \quad \text{for $i\neq t$},\\
    e'_t & \mapsto e_{t_1}+e_{t_2}.
\end{align*}
Note that $\sigma$ satisfies the conditions in Definition \ref{def for quasi}. It defines quasi-cluster morphisms
\[
\sigma:\mathscr{A}'\rightarrow \mathscr{A} \quad \text{and} \quad \sigma:\mathscr{X}\rightarrow \mathscr{X}'.
\]
The next proposition is direct consequence of the construction of $\sigma$.
\begin{prop}\label{merging frozen}  The quasi-cluster morphism $\sigma:\mathscr{A}'\rightarrow \mathscr{A}$ embeds $\mathscr{A}'$ as a subvariety of $\mathscr{A}$ determined by the locus $\left\{A_i=A_j\right\}$.
\end{prop}

\bigskip

\section{Appendix II: Double Bott-Samelson Cells}
\label{Sec 3}
\subsection{Definition and Basic Properties}

Double Bott-Samelson (BS) cells, introduced in \cite{SWflag}, are moduli spaces of flags with prescribed relative positions encoded by positive braids. In this section we briefly recall their definition and basic properties following \textit{loc. cit}. Theorem \ref{3.6} establishes natural isomorphisms between the augmentation varieties of positive braid closures and the double BS cells associated with $\SL_n$.

Let $\B_\pm$ be a pair of opposite Borel subgroups of a Kac-Moody group $\G$ and let $\U_\pm:=\left[\B_\pm,\B_\pm\right]$ be the  maximal unipotent subgroups. There are flag varieties $\mathcal{B}_+:=\G/\B_+$ and $\mathcal{B}_-:=\B_-\backslash \G$. 
By replacing $\B_\pm$ with $\U_\pm$ we define \emph{decorated flag varieties} $\mathcal{A}_+:=\G/\U_+$ and $\mathcal{A}_-:=\U_-\backslash \G$. There are natural projections $\pi:\mathcal{A}_\pm\rightarrow \mathcal{B}_\pm$. If $\pi(\A)=\B$ then we say that $\A$ is a \emph{decoration} over $\B$.

We denote elements in $\mathcal{B}_+$ as $\B^i$ and elements in $\mathcal{B}_-$ as $\B_i$. The same convention is applied to $\mathcal{A}_\pm$  with the letter $\B$ replaced by $\A$.

Let $\T:=\B_+\cap \B_-$ and let $\W:=\left.N(\T)\right/\T$ be the Weyl group. Consider the Bruhat decompositions and  Birkhoff decomposition
\[
\G=\bigsqcup_{w\in \W}\B_+ w\B_+=\bigsqcup_{w\in \W}\B_-w\B_-=\bigsqcup_{w\in \W}\B_-w\B_+.
\]
We adopt the convention of writing 
\[
\left\{\begin{array}{ll}
    \xymatrix{x\B_+\ar[r]^w &y\B_+} & \text{if $x^{-1}y\in \B_+w\B_+$},  \\
    \xymatrix{\B_-x \ar[r]^w &\B_-y} & \text{if $xy^{-1}\in \B_-w\B_-$},\\
    \xymatrix{\B_-x\ar@{-}[r]^w & y\B_+} & \text{if $xy\in \B_-w\B_+$}.
\end{array}\right.
\]
We often omit $w$ in the notation when it is the identity. Moreover, when decorated flags are involved, the notations  only concern the underlying flags; for example, $\xymatrix{\B^i\ar[r]^w & \A^j}$ means $\xymatrix{\B^i\ar[r]^w & \pi\left(\A^j\right)}$.

For a positive braid word $\beta=s_{i_1}\dots s_{i_m}$, a chain 
$\xymatrix{\B^0\ar[r]^{s_{i_1}}  &\cdots \ar[r]^{s_{i_m}} &\B^m}$
will be abbreviated as $\xymatrix{\B^0\ar@{-->}[r]^\beta&\B^m}$.
By Theorem 2.18 of \cite{SWflag}, the chains of flags associated to different words of the positive braid $[\beta]$ have a natural one-to-one correspondence. In this sense, the chain $\xymatrix{\B^0\ar@{-->}[r]^\beta&\B^m}$ does not depend on the word $\beta$ chosen. 

\begin{defn} Let $\beta$ and $\gamma$ be positive braids. The \emph{half decorated double BS cell} $\conf_\beta^\gamma(\mathcal{C})$, viewed as a $\bbZ$-scheme, is a moduli space parametrizes $\G$-orbits of the chains of flags
\[
\begin{tikzpicture}[baseline=4ex]
\node (u0) at (0,1.3) [] {$\B_0$};
\node (u1) at (1.3,1.3) [] {$\A_m$};
\node (d0) at (0,0) [] {$\B^0$};
\node (d1) at (1.3,0) [] {$\B^l$};
\draw [dashed, ->] (d0) -- node [below] {$\beta$} (d1);
\draw [dashed, ->] (u0) -- node [above] {$\gamma$} (u1);
\draw (d0) -- (u0);
\draw (d1) -- (u1);
\end{tikzpicture}.
\]
If one forgets to choose a decoration $\A_m$ over $\B_m$, then the resulting space is denoted by $\conf_\beta^\gamma(\mathcal{B})$.
Denote by  $\pi$ the forgetful map from ${\rm Conf}_\beta^\gamma(\mathcal{C})$ to $\conf_\beta^\gamma(\mathcal{B})$.
\end{defn}

\begin{rmk} This version of double BS cells is slightly different from those in \cite{SWflag}: first, the two chains of flags swap places with the $\mathcal{B}_+$-chain at the bottom and the $\mathcal{B}_-$-chain at the top now; second, there is only one decoration $\A_m$ over $\B_m$ and the flag $\B^0$ is no longer decorated. 
\end{rmk}

For  $\xymatrix{\B_0 \ar@{-}[r] & \B^0 \ar[r]^{s_i} &\B^1}$, there is a unique flag $\B_{-1}$ such that $\xymatrix{\B^0 \ar@{-}[r]^{s_i} & \B_{-1} \ar[r]^{s_i} &\B_0}$. It then follows from $\xymatrix{\B_{-1} \ar@{-}[r]^{s_i} & \B^0\ar[r]^{s_i} & \B^1}$ that $\xymatrix{\B_{-1} \ar@{-}[r] & \B^1}$. This construction gives rise to the following reflection maps.
\begin{defn} The left reflection map $l_i: \conf^\gamma_{s_i\beta}(\mathcal{C})\rightarrow \conf^{s_i\gamma}_{\beta}(\mathcal{C})$ is an isomorphism mapping
\[
\vcenter{\vbox{\xymatrix{ & \B_0  \ar@{-}[dl]\ar@{-->}[r]^{\gamma} &\A_m \ar@{-}[d] \\ \B^0 \ar[r]_{s_i} & \B^1 \ar@{-->}[r]_{\beta}&\B^n}}}  \quad \longmapsto \quad \vcenter{\vbox{\xymatrix{\B_{-1} \ar[r]^{s_i} \ar@{-}[dr] & \B_0 \ar@{-->}[r]^{\gamma} &\A_m \ar@{-}[d] \\   & \B^1 \ar@{-->}[r]_\beta & \B^n}}} 
\]
Its inverse map $l^i:\conf^{s_i\gamma}_{\beta}(\mathcal{C}) \rightarrow \conf^{\gamma}_{s_i\beta}(\mathcal{C})$ is defined by an analogous process.
\end{defn}

Let $\varphi_i:\SL_2\rightarrow \G$ be the group homomorphism associated to the simple root $\alpha_i$. Define
\[
e_i(q) :=\varphi_i\begin{pmatrix} 1 & q \\ 0 & 1 \end{pmatrix}, \quad 
e_{-i}(q) :=\varphi_i\begin{pmatrix} 1 & 0 \\ q & 1 \end{pmatrix}, \quad
\overline{s}_i :=\varphi_i\begin{pmatrix} 0 & -1 \\ 1 & 0\end{pmatrix}, \quad 
\doverline{s}_i :=\varphi_i\begin{pmatrix} 0 & 1 \\ -1 & 0\end{pmatrix}.
\]
Consider a reduced expression $w=s_{i_1}\ldots s_{i_n}$ in $\W$. Let
\[
\overline{w}=\overline{s}_{i_1}\ldots \overline{s}_{i_n}, \hskip 8mm \doverline{w}=\doverline{s}_{i_1}\ldots \doverline{s}_{i_n}.
\]
The elements $\overline{w}$ and $\doverline{w}$ in $\G$ do not depend on the reduced expression chosen. 
We  set
\begin{equation}\label{R def}
R_i(q) := e_i(q)\overline{s}_i = \varphi_i\begin{pmatrix} q & -1 \\ 1 & 0 \end{pmatrix}.
\end{equation}

\begin{lem}\label{3.3} Fix a flag $\B^j$. The space of flags $\B^k$ such that $\xymatrix{\B^j\ar[r]^{s_i}&\B^k}$ is isomorphic to $\mathbb{A}^1$. In particular, if $\B^j=\B_+$, then $\B^k=R_i(q)\B_+$ for some unique $q\in \mathbb{A}^1$.
\end{lem}
\begin{proof} It suffices to prove the lemma for $\B^j=\B_+$. Let $\U_i:=\left\{e_i(t)\ \middle| \ t\in \mathbb{A}^1\right\}$ be the 1-dimensional unipotent subgroup associated to the simple root $\alpha_i$ and let $\Q_i:=\B_+\cap s_i\B_+ s_i$. By \cite[Lemma 6.1.3]{Kum}, we know that $\B_+=\U_i\Q_i$. Therefore
\[
\left.\B_+ s_i \B_+\right/\B_+=\left.\U_i\Q_i s_i\B_+\right/\B_+=\left.\U_is_i\Q_i\B_+\right/\B_+=\left.\U_is_i\B_+\right/\B_+.
\]
Hence $\B^k=R_i(q)\B_+$ for some unique $q\in \mathbb{A}^1$.
\end{proof}

We prove an important property of the double BS cells, following \cite[\S 2.4]{SWflag}. 

\begin{prop} 
\label{caosno}
The space $\conf^\gamma_\beta(\mathcal{C})$ is the non-vanishing locus of a polynomial in $\mathbb{A}^{l(\beta)+l(\gamma)}$. Consequently, it is a smooth affine variety.
\end{prop}
\begin{proof} 
It suffices to prove the lemma for $\conf^e_\beta(\mathcal{C})$; the general case will follow by using the reflections to shift the top $\gamma$ to the bottom.
Suppose $\beta$ is of length $l$. Every point of $\conf^e_\beta(\mathcal{C})$ admits  a unique representative as follows
\begin{equation}
\label{acnonoda}
\begin{gathered}
\begin{tikzpicture}
\node (u) at (3,1.5) [] {$\U_-$};
\node (d0) at (0,0) [] {$\B_+$};
\node (d1) at (1.5,0) [] {$\B^1$};
\node (d2) at (3,0) [] {$\B^2$};
\node (d3) at (4.5,0) [] {$\cdots$};
\node (d4) at (6,0) [] {$\B^l$};
\draw (u) -- (d0);
\draw [->] (d0) -- node [below] {$s_{i_1}$} (d1);
\draw [->] (d1) -- node [below] {$s_{i_2}$} (d2);
\draw [->] (d2) -- node [below] {$s_{i_3}$} (d3);
\draw [->] (d3) -- node [below] {$s_{i_l}$} (d4);
\draw (d4) -- (u);
\end{tikzpicture}
\end{gathered}
\end{equation}
Using Lemma \ref{3.3} recursively, we obtain parameters $(q_1,\dots,q_l) \in \mathbb{A}^l$ such that 
\[
\B^k= R_{i_1}(q_1)\cdots R_{i_k}(q_k) \B_+, \quad \quad  k=1,\ldots,l.
\]
By definition, we require that the rightmost pair $({\U_-, \B^{l}})$ is in general position.

Let $\omega_i$ be the $i$th fundamental weight. 
The \textit{$i$th principal minor} $\Delta_i: \G\rightarrow \mathbb{A}$ is  a regular function uniquely determined by the following two conditions:
(1) $\Delta_i(u_-gu_+)=\Delta_i(g)$, where $u_\pm \in \U_\pm$;  (2) $\Delta_i(h)=h^{\omega_i}$ for $h\in \T$. When $\G=\SL_n$, the function $\Delta_i$  coincides with \eqref{cfqhbvo}. Note that $g\in \B^-\B^+$ if and only if $\Delta_i(g)\neq 0$ for all $i$.
Therefore the pair $({\U_-, \B^{l}})$ is in general position  if and only if 
\[
f(q_1,\ldots,q_l):=\prod_{1\leq i \leq {\rm rk}\G} \Delta_i\left(R_{i_1}(q_1)\cdots R_{i_l}(q_l)\right)\neq 0. \qedhere
\]
\end{proof}

\subsection{Cluster Structures on Double Bott-Samelson Cells}\label{sec 3.2}

A pair of positive braids $(\beta,\gamma)$ can be regarded as a single braid in the product $\Br\times \Br$. We shall prove that every word of $(\beta,\gamma)$  gives rise to a cluster {\it seed} of $\conf^\gamma_\beta(\mathcal{C})$. First, each word determines a labeling of arrows and a triangulation on the configuration diagram. Then we require that every pair of flags that are connected by a diagonal in the triangulation are in general position. The subspace of $\conf^\gamma_\beta(\mathcal{C})$ that satisfy these general position conditions is an algebraic torus. The algebraic tori obtained from all words of $(\beta, \gamma)$ form a subset of the atlas of cluster charts.

In detail, let $\vec{t}$ be a word of $(\beta,\gamma)$. We label the arrows and draw the triangulation on the configuration diagram according to $\vec{t}$ as shown in Example \ref{3.10}. On top of the triangulation, we draw $\rank(\G)$ many parallel lines, each of which represents a simple root of $\G$. Triangles in the triangulation are either upward pointing or downward pointing (as shown below), and depending on the orientation and the labeling of the base, each triangle places a node at one of the lines, cutting it into segments called \emph{strings}. The segments from such cutting become the vertices of the quiver $Q_{\vec{t}}$, and the arrows in $Q_{\vec{t}}$ are drawn according to the pictures below, where the dashed arrows between different levels $i\neq j$ are weighted by weights that are related to Cartan numbers (see \cite{SWflag} for more details). In particular, in the simply-laced cases (which include $\SL_n$), the dashed arrows all have weight $1/2$. In the end, we delete the left most vertices (together with all incident arrows) and freeze the right most vertices on each level.

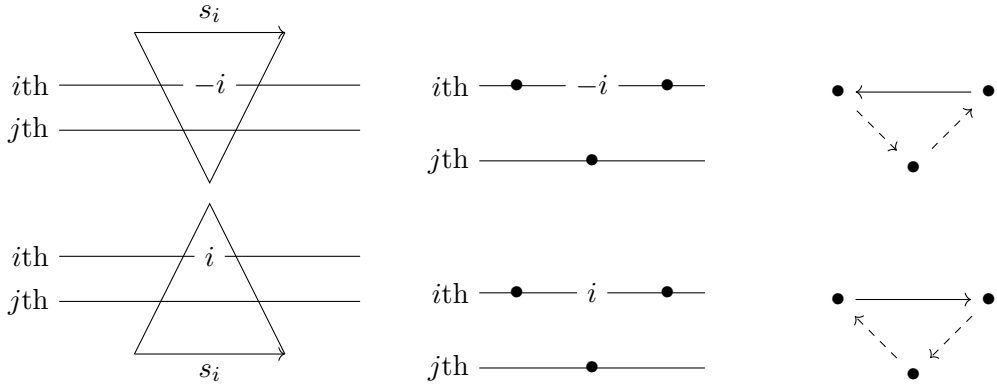
\[
\begin{tikzpicture}[scale=0.8,baseline=-15]
\draw (-1,1) -- (0,-1) -- (1,1);
\draw [->] (-1,1) -- node [above] {$s_i$} (1,1);
\node (0) at (0,0.3) [] {$-i$};
\draw (-2,-0.3) node [left] {$j$th} -- (2,-0.3);
\draw (-2,0.3) node [left] {$i$th} -- (0) -- (2,0.3); 
\end{tikzpicture}
\quad \quad 
\begin{tikzpicture}[scale=0.8]
\node (0) at (0,0.5) [] {$-i$};
\draw (-1.5,-0.5) node [left] {$j$th} -- (1.5,-0.5);
\draw (-1.5,0.5) node [left] {$i$th} -- (0) -- (1.5,0.5); 
\node at (-1,0.5) [] {$\bullet$};
\node at (1,0.5) [] {$\bullet$};
\node at (0,-0.5) [] {$\bullet$};
\end{tikzpicture}
\quad\quad\quad \quad 
\begin{tikzpicture}[scale=0.8]
\node (ul) at (-1,0.5) [] {$\bullet$};
\node (ur) at (1,0.5) [] {$\bullet$};
\node (d) at (0,-0.5) [] {$\bullet$};
\draw [<-] (ul) -- (ur);
\draw [dashed, <-] (ur) -- (d);
\draw [dashed, <-] (d) -- (ul);
\end{tikzpicture}
\]
\[
\begin{tikzpicture}[scale=0.8,baseline=-15]
\draw (-1,-1) -- (0,1) -- (1,-1);
\draw [->] (-1,-1) -- node [below] {$s_i$} (1,-1);
\node (0) at (0,0.3) [] {$i$};
\draw (-2,-0.3) node [left] {$j$th} -- (2,-0.3);
\draw (-2,0.3) node [left] {$i$th} -- (0) -- (2,0.3); 
\end{tikzpicture}
\quad \quad 
\begin{tikzpicture}[scale=0.8]
\node (0) at (0,0.5) [] {$i$};
\draw (-1.5,-0.5) node [left] {$j$th} -- (1.5,-0.5);
\draw (-1.5,0.5) node [left] {$i$th} -- (0) -- (1.5,0.5); 
\node at (-1,0.5) [] {$\bullet$};
\node at (1,0.5) [] {$\bullet$};
\node at (0,-0.5) [] {$\bullet$};
\end{tikzpicture}
\quad\quad\quad \quad 
\begin{tikzpicture}[scale=0.8]
\node (ul) at (-1,0.5) [] {$\bullet$};
\node (ur) at (1,0.5) [] {$\bullet$};
\node (d) at (0,-0.5) [] {$\bullet$};
\draw [->] (ul) -- (ur);
\draw [dashed, ->] (ur) -- (d);
\draw [dashed, ->] (d) -- (ul);
\end{tikzpicture}
\]

To define the cluster $\mathrm{K}_2$ coordinates, we first need to decorate the flags. 
 Two decorated flags $\xymatrix{x\U_+ \ar[r]^w & y\U_+}$ (resp. $\xymatrix{\U_-x \ar[r]^w & \U_-y}$) are said to be \emph{compatible} if $x^{-1}y\in \U_+\overline{w}\U_+$ (resp. $xy^{-1}\in \U_-\doverline{w}\U_-$). Two decorated flags $\xymatrix{\U_-x\ar@{-}[r] & y\U_+}$ is called a \emph{pinning} if $xy\in \U_-\U_+$.

\begin{lem} Given $\xymatrix{\B \ar[r]^w & \B'}$ (resp. $\xymatrix{\B'\ar[r]^w & \B}$ or $\xymatrix{\B \ar@{-}[r] & \B'}$), for every decoration $\A$ over $\B$, there exists a unique decoration $\A'$ over $\B'$, such that $\xymatrix{\A \ar[r]^w & \A'}$ are compatible (resp. $\xymatrix{\A' \ar[r]^w & \A}$ are compatible or $\xymatrix{\A \ar@{-}[r] & \A'}$ is a pinning).
\end{lem}

Using the above lemma, we can begin with the decoration $\A_m$ over $\B_m$ and induce decorations one-by-one over the rest flags following the $C$-shape path illustrated by the dashed circles below
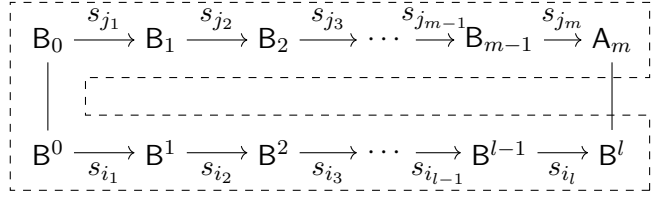
\[
\begin{tikzpicture}[baseline=4ex]
\node (u0) at (0,1.5) [] {$\B_0$};
\node (u1) at (1.5,1.5) [] {$\B_1$};
\node (u2) at (3,1.5) [] {$\B_2$};
\node (u3) at (4.5,1.5) [] {$\cdots$};
\node (u4) at (6,1.5) [] {$\B_{m-1}$};
\node (u5) at (7.5,1.5) [] {$\A_m$};
\node (d0) at (0,0) [] {$\B^0$};
\node (d1) at (1.5,0) [] {$\B^1$};
\node (d2) at (3,0) [] {$\B^2$};
\node (d3) at (4.5,0) [] {$\cdots$};
\node (d4) at (6,0) [] {$\B^{l-1}$};
\node (d5) at (7.5,0) [] {$\B^l$};
\draw [->] (d0) -- node [below] {$s_{i_1}$} (d1);
\draw [->] (d1) -- node [below] {$s_{i_2}$} (d2);
\draw [->] (d2) -- node [below] {$s_{i_3}$} (d3);
\draw [->] (d3) -- node [below] {$s_{i_{l-1}}$} (d4);
\draw [->] (d4) -- node [below] {$s_{i_l}$} (d5);
\draw [->] (u0) -- node [above] {$s_{j_1}$} (u1);
\draw [->] (u1) -- node [above] {$s_{j_2}$} (u2);
\draw [->] (u2) -- node [above] {$s_{j_3}$} (u3);
\draw [->] (u3) -- node [above] {$s_{j_{m-1}}$} (u4);
\draw [->] (u4) -- node [above] {$s_{j_m}$} (u5);
\draw (d0) -- (u0);
\draw (d5) -- (u5);
\draw [dashed] (8,-0.5) -- (-0.5,-0.5) -- (-0.5,2) -- (8,2) -- (8,1) -- (0.5,1) -- (0.5,0.5) -- (8,0.5) -- cycle;
\end{tikzpicture}
\]

The next proposition presents a \textit{standard representative} for every point in $\conf^\gamma_\beta(\mathcal{C})$. 

\begin{prop}\label{unique rep} Every point in $\conf_\beta^\gamma(\mathcal{C})$ admits a unique representative in the following form:
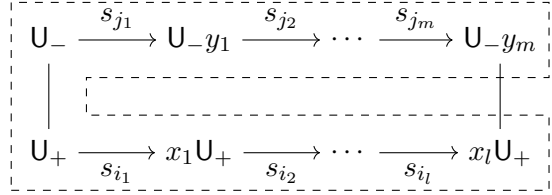
\[
\begin{tikzpicture}
\node (u0) at (0,1.5) [] {$\U_-$};
\node (u1) at (2,1.5) [] {$\U_-y_1$};
\node (u2) at (4,1.5) [] {$\cdots$};
\node (u3) at (6,1.5) [] {$\U_-y_m$};
\node (d0) at (0,0) [] {$\U_+$};
\node (d1) at (2,0) [] {$x_1\U_+$};
\node (d2) at (4,0) [] {$\cdots$};
\node (d3) at (6,0) [] {$x_l\U_+$};
\draw [->] (u0) -- node [above] {$s_{j_1}$} (u1);
\draw [->] (u1) -- node [above] {$s_{j_2}$} (u2);
\draw [->] (u2) -- node [above] {$s_{j_m}$} (u3);
\draw [->] (d0) -- node [below] {$s_{i_1}$} (d1);
\draw [->] (d1) -- node [below] {$s_{i_2}$} (d2);
\draw [->] (d2) -- node [below] {$s_{i_l}$} (d3);
\draw (u0) -- (d0);
\draw (u3) -- (d3);
\draw [dashed] (-0.5,-0.5) -- (6.65,-0.5) -- (6.65,0.5) -- (0.5,0.5) -- (0.5,1) -- (6.65,1) -- (6.65,2) -- (-0.5,2) -- cycle;
\end{tikzpicture}
\]
where
\[
x_k = R_{i_1}\left(q_1\right)R_{i_2}\left(q_2\right)\dots R_{i_k}\left(q_k\right), \qquad
y_k = R_{j_k}\left(p_k\right)\dots R_{j_2}\left(p_2\right) R_{j_1}\left(p_1\right).
\]
This gives an open embedding $\conf^\gamma_\beta(\mathcal{C})\hookrightarrow \mathbb{A}^m_{p_1,\dots, p_m}\times \mathbb{A}^l_{q_1,\dots, q_l}$.
\end{prop}
\begin{proof} Let us first verify that adjacent decorated flags along the top chain and the bottom chain are compatible. Let $x_0=y_0=e$. Note that
\[
\U_+x_{k-1}^{-1}x_k\U_+=\U_+e_{i_k}\left(q_k\right)\overline{s}_{i_k}\U_+=\U_+\overline{s}_{i_k}\U_+,
\]
\[
\U_-y_{k-1}y_k^{-1}\U_-=\U_-\left(e_{j_k}\left(p_k\right)\overline{s}_{j_k}\right)^{-1}\U_-=\U_-e_{-j_k}\left(-p_k\right)\doverline{s}_{j_k}\U_-=\U_-\doverline{s}_{j_k}\U_-.
\]

Since a compatible decoration on one end of any adjacent pair of flags along either of the horizontal chains can be uniquely determined by the decoration on the other end of the pair, the existence of uniqueness of such representative automatically follows from the fact that $\G$ acts freely and transitively on the space of pinnings. 
\end{proof}

Now for a fixed word $\vec{t}$ of $(\beta,\gamma)$, we get a quiver $Q_{\vec{t}}$ with vertices corresponding to strings, which necessarily cross certain diagonals (possibly more than one) in the triangulation. 
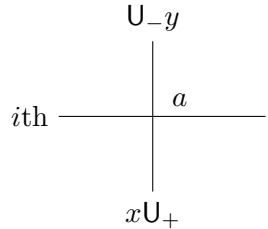
\[
\begin{tikzpicture}[scale=0.8]
\draw (-1.75,0) node[left] {$i$th} -- node [above right] {$a$} (1,0);
\draw (-0.5,-1) node [below] {$x\U_+$} -- (-0.5,1) node [above] {$\U_-y$};
\end{tikzpicture}
\]
The cluster $\mathrm{K}_2$ coordinate associated to the string $a$ is defined to be the $i$th principal minor of $xy$:
\[
A_a=\Delta_i(yx).
\]
The function $A_a$ is independent of the choice of diagonals if $a$ crosses more than one diagonals.

\begin{exmp}\label{3.10} Let $\G=\SL_3$, $\beta=s_2s_1s_2s_1$, and $\gamma = s_2s_1$. For $\Br\times \Br$, we use negative numbers for letters in the second factor. The word $\vec{t} = (2,-2,1,2,-1,1)$ for $(\beta,\gamma)$ gives rise to the following triangulation, string diagram, and quiver
\[
\begin{tikzpicture}
\foreach \i in {0,1,2}
    {
    \node (u\i) at (2+\i*2,2) [] {$\A_\i$};
    }
\foreach \i in {0,...,4}
    {
    \node (d\i) at (\i*2,0) [] {$\A^\i$};
    }
\draw [->] (u0) -- node [above] {$s_2$} (u1);
\draw [->] (u1) -- node [above] {$s_1$} (u2);
\draw [->] (d0) -- node [below] {$s_2$} (d1);
\draw [->] (d1) -- node [below] {$s_1$} (d2);
\draw [->] (d2) -- node [below] {$s_2$} (d3);
\draw [->] (d3) -- node [below] {$s_1$} (d4);
\draw (d0) -- (u0) -- (d1) --(u1)--(d2);
\draw (u1) -- (d3) -- (u2) -- (d4);
\draw [dashed] (-1,-0.5) --(2,2.5) -- (6,2.5) -- (7,1.5) -- (2,1.5) -- (1,0.5) -- (8,0.5) -- (9,-0.5) -- cycle;
\end{tikzpicture}
\]
\[
\begin{tikzpicture}
\foreach \i in {0,1,2}
    {
    \node [lightgray] (u\i) at (2+\i*2,2) [] {$\A_\i$};
    }
\foreach \i in {0,...,4}
    {
    \node [lightgray] (d\i) at (\i*2,0) [] {$\A^\i$};
    }
\draw [lightgray,->] (u0) -- node [above] {$s_2$} (u1);
\draw [lightgray,->] (u1) -- node [above] {$s_1$} (u2);
\draw [lightgray,->] (d0) -- node [below] {$s_2$} (d1);
\draw [lightgray,->] (d1) -- node [below] {$s_1$} (d2);
\draw [lightgray,->] (d2) -- node [below] {$s_2$} (d3);
\draw [lightgray,->] (d3) -- node [below] {$s_1$} (d4);
\draw [lightgray] (d0) -- (u0) -- (d1) --(u1)--(d2);
\draw [lightgray] (u1) -- (d3) -- (u2) -- (d4);
\node (u1) at (3.7,1.3) [] {$1$};
\node (u2) at (5.5,1.3) [] {$-1$};
\node (u3) at (6.3,1.3) [] {$1$};
\node (d1) at (1.2,0.7) [] {$2$};
\node (d2) at (2.3,0.7) [] {$-2$};
\node (d3) at (4.5,0.7) [] {$2$};
\draw (-1,0.7) node [left] {$2$nd} -- (d1) -- (d2) -- (d3) -- (9,0.7);
\draw (-1,1.3) node [left] {$1$st} -- (u1) -- (u2) -- (u3) -- (9,1.3);
\end{tikzpicture}
\]
\[
\begin{tikzpicture}
\node [lightgray] (u1) at (3.5,1.5) [] {$1$};
\node [lightgray](u2) at (5.5,1.5) [] {$-1$};
\node [lightgray](u3) at (6.5,1.5) [] {$1$};
\node [lightgray](d1) at (1,0.5) [] {$2$};
\node [lightgray](d2) at (2.5,0.5) [] {$-2$};
\node [lightgray](d3) at (4.5,0.5) [] {$2$};
\draw [lightgray](-1,0.5) node [left] {$2$nd} -- (d1) -- (d2) -- (d3) -- (9,0.5);
\draw [lightgray](-1,1.5) node [left] {$1$st} -- (u1) -- (u2) -- (u3) -- (9,1.5);
\node (1u) at (4.25,1.5) [] {$\bullet$};
\node (2u) at (6,1.5) [] {$\bullet$};
\node (3u) at (7.5,1.5) [] {$\square$};
\node (1d) at (1.75,0.5) [] {$\bullet$};
\node (2d) at (3.5,0.5) [] {$\bullet$};
\node (3d) at (7.5,0.5) [] {$\square$};
\draw [->] (2d) -- (1d);
\draw [->] (2d) -- (3d);
\draw [->] (2u) -- (1u);
\draw [->] (2u) -- (3u);
\draw [->] (1u) -- (2d);
\draw [->] (3d) -- (2u);
\draw [dashed, ->] (3u) -- (3d);
\end{tikzpicture}
\]
\end{exmp}

\begin{rmk} In \cite{SWflag} a cluster $\mathrm{K}_2$ structure is constructed on the decorated double BS cell $\conf^\gamma_\beta\left(\mathcal{A}_\sc\right)$ for a simply-connected group $\G$, which has frozen vertices on both sides of the quiver. The cluster $\mathrm{K}_2$ structure on $\conf^\gamma_\beta(\mathcal{C})$ is essentially obtained from that of $\conf^\gamma_\beta\left(\mathcal{A}_\sc\right)$ by setting all the frozen variables on the left to be 1 due to the pinning condition on $\xymatrix{\A_0 \ar@{-}[r] & \A^0}$.
\end{rmk}

The next Proposition provides an interpretation of left reflections in terms of standard representatives in Proposition \ref{unique rep}. It implies that the left reflections are cluster transformations.

\begin{prop}\label{leftreflection} The left reflection $\conf^\gamma_{s_i\beta}(\mathcal{C})\rightarrow \conf^{s_i\gamma}_\beta(\mathcal{C})$ can be expressed in terms of standard representatives as
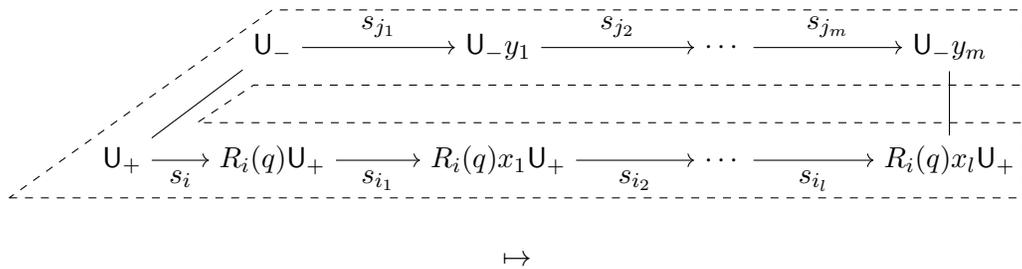
\[
\begin{tikzpicture}
\node (u0) at (0,1.5) [] {$\U_-$};
\node (u1) at (3,1.5) [] {$\U_-y_1$};
\node (u2) at (6,1.5) [] {$\cdots$};
\node (u3) at (9,1.5) [] {$\U_-y_m$};
\node (d-1) at (-2,0) [] {$\U_+$};
\node (d0) at (0,0) [] {$R_i(q)\U_+$};
\node (d1) at (3,0) [] {$R_i(q)x_1\U_+$};
\node (d2) at (6,0) [] {$\cdots$};
\node (d3) at (9,0) [] {$R_i(q)x_l\U_+$};
\draw [->] (u0) -- node [above] {$s_{j_1}$} (u1);
\draw [->] (u1) -- node [above] {$s_{j_2}$} (u2);
\draw [->] (u2) -- node [above] {$s_{j_m}$} (u3);
\draw [->] (d-1) -- node [below] {$s_i$} (d0);
\draw [->] (d0) -- node [below] {$s_{i_1}$} (d1);
\draw [->] (d1) -- node [below] {$s_{i_2}$} (d2);
\draw [->] (d2) -- node [below] {$s_{i_l}$} (d3);
\draw (u0) -- (d-1);
\draw (u3) -- (d3);
\draw [dashed] (-3.5,-0.5) -- (10,-0.5) -- (10,0.5) -- (-1,0.5) -- (-0.25,1) -- (10,1) -- (10,2) -- (0,2) -- cycle;
\end{tikzpicture}
\]
\[
\longmapsto
\]
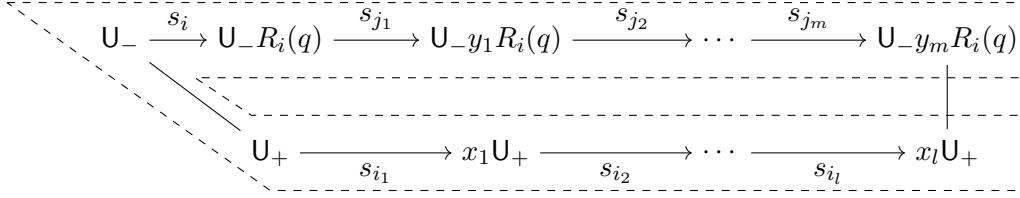
\[
\begin{tikzpicture}
\node (u0) at (0,0) [] {$\U_+$};
\node (u1) at (3,0) [] {$x_1\U_+$};
\node (u2) at (6,0) [] {$\cdots$};
\node (u3) at (9,0) [] {$x_l\U_+$};
\node (d-1) at (-2,1.5) [] {$\U_-$};
\node (d0) at (0,1.5) [] {$\U_-R_i(q)$};
\node (d1) at (3,1.5) [] {$\U_-y_1R_i(q)$};
\node (d2) at (6,1.5) [] {$\cdots$};
\node (d3) at (9,1.5) [] {$\U_-y_mR_i(q)$};
\draw [->] (u0) -- node [below] {$s_{i_1}$} (u1);
\draw [->] (u1) -- node [below] {$s_{i_2}$} (u2);
\draw [->] (u2) -- node [below] {$s_{i_l}$} (u3);
\draw [->] (d-1) -- node [above] {$s_i$} (d0);
\draw [->] (d0) -- node [above] {$s_{j_1}$} (d1);
\draw [->] (d1) -- node [above] {$s_{j_2}$} (d2);
\draw [->] (d2) -- node [above] {$s_{j_m}$} (d3);
\draw (u0) -- (d-1);
\draw (u3) -- (d3);
\draw [dashed] (-3.5,2) -- (10,2) -- (10,1) -- (-1,1) -- (-0.25,0.5) -- (10,0.5) -- (10,-0.5) -- (0,-0.5) -- cycle;
\end{tikzpicture}
\]
\end{prop}
\begin{proof} The left reflection does the following.
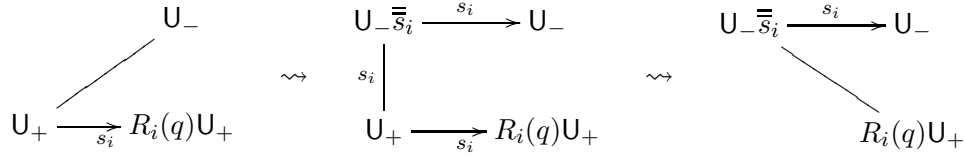
\[
\vcenter{\vbox{\xymatrix{ & \U_-  \ar@{-}[dl] \\ \U_+ \ar[r]_{s_i} & R_i(q)\U_+}}} \quad \rightsquigarrow \quad \vcenter{\vbox{\xymatrix{\U_-\doverline{s}_i \ar[r]^{s_i} \ar@{-}[d]_{s_i} & \U_-  \\   \U_+ \ar[r]_{s_i} & R_i(q)\U_+}}} \quad \rightsquigarrow \quad \vcenter{\vbox{\xymatrix{\U_-\doverline{s}_i \ar[r]^{s_i} \ar@{-}[dr] & \U_-  \\   & R_i(q)\U_+}}} 
\]
To restore to the standard representative, we need to act on the resulting configuration by $\left(R_i(q)\right)^{-1}$. Note that under the such action, $x\U_+\mapsto \left(R_i(q)\right)^{-1}x\U_+$ and $\U_-y\mapsto \U_-yR_i(q)$. It is not hard to see that such action will give the standard configuration as claimed in the proposition. 
\end{proof}

\subsection{An Open Embedding}\label{sec 3.3}

In this section we construct an open embedding $\psi:\conf^\gamma_\beta(\mathcal{C})\times \mathbb{G}_m\hookrightarrow \conf^\gamma_{s_i\beta}(\mathcal{C})$ whose image is the localization (freezing) at a cluster variable of the latter.

Recall from Lemma \ref{3.3} that the moduli space of $\B^{-1}$ that fits into the triangle in the picture on the left below is parametrized by the multiplicative group scheme $\mathbb{G}_m$. Note that the base change of $\mathbb{G}_m$ to any field $\mathds{k}$ is isomorphic to $\mathds{k}^\times$ as affine schemes over $\mathds{k}$.
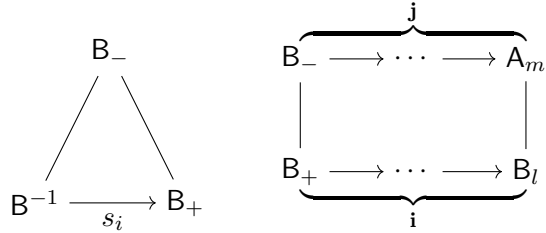
\[
\begin{tikzpicture}
\node (u) at (0,1) [] {$\B_-$};
\node (dl) at (-1,-1) [] {$\B^{-1}$};
\node (dr) at (1,-1) [] {$\B_+$};
\draw [->] (dl) -- node [below] {$s_i$} (dr);
\draw (dl) -- (u) --(dr);
\end{tikzpicture} \quad \quad \begin{tikzpicture}
\node (u0) at (0,1.5) [] {$\B_-$};
\node (u1) at (1.5,1.5) [] {$\cdots$};
\node (u2) at (3,1.5) [] {$\A_m$};
\node (d0) at (0,0) [] {$\B_+$};
\node (d1) at (1.5,0) [] {$\cdots$};
\node (d2) at (3,0) [] {$\B_l$};
\draw [->] (u0) -- (u1);
\draw [->] (u1) -- (u2);
\draw [->] (d0) -- (d1);
\draw [->] (d1) -- (d2);
\draw (u0) -- (d0);
\draw (d2) -- (u2);
\node at (1.5,-0.5) [] {$\underbrace{\hspace{3cm}}_\vec{i}$};
\node at (1.5,2) [] {$\overbrace{\hspace{3cm}}^\vec{j}$};
\end{tikzpicture}
\]
On the other hand, consider a standard representative and let us temporarily forget about the decorations on the pinning and the bottom chain, as shown in the picture on the right above. By gluing these two figures along the pinning $\xymatrix{\B_-\ar@{-}[r] & \B_+}$, we end up with a point in $\conf^\gamma_{s_i\beta}(\mathcal{C})$, which defines a morphism 
\begin{equation}\label{def open embedding}
\psi:\conf^\gamma_\beta(\mathcal{C})\times \mathbb{G}_m \rightarrow \conf^\gamma_{s_i\beta}(\mathcal{C}),
\end{equation}
It is easy to see that $\psi$ is an open embedding. 

\begin{prop}\label{image phi} The image of $\psi$ in $\conf^\gamma_{s_i\beta}(\mathcal{C})$ 
is the distinguished open subset corresponding to the localization (freezing) at the leftmost cluster variable $A_c$ in the picture below
\begin{equation}
\label{snvodqv}
\begin{gathered}
\begin{tikzpicture}[baseline=10]
\node (u0) at (2,1.5) [] {$\B_0$};
\node (u1) at (4,1.5) [] {$\cdots$};
\node at (4,1.5) [above ] {$\gamma$};
\node (u2) at (6,1.5) [] {$\A_m$};
\node (d0) at (0,0) [] {$\B^{-1}$};
\node (d1) at (2,0) [] {$\B^0$};
\node (d2) at (4,0) [] {$\cdots$};
\node at (4,0) [below] {$\beta$};
\node (d3) at (6,0) [] {$\B_l$};
\draw [->] (u0) -- (u1);
\draw [->] (d0) -- node [below] {$s_i$} (d1);
\draw [->] (d1) -- (d2);
\draw [->] (u1) -- (u2);
\draw [->] (d2) -- (d3);
\draw (d3) -- (u2);
\draw (u0) -- (d0);
\draw (u0) -- (d1);
\draw (1.75,0.75) node [left] {$i$\emph{th}} --  node [above right] {$c$} (3,0.75);
\end{tikzpicture}.
\end{gathered}
\end{equation}
\end{prop}
\begin{proof} There is a unique representative of \eqref{snvodqv} such that $\B_0=\B_-$, $\B^{-1}=\B_+$, and $\B^0=R_i(d)\B_+$.
The principal minors of $R_i(d)$ are
\[
\Delta_k(R_i(d))=\left\{ \begin{array}{ll} d & \mbox{if } i=k;\\
1 & \mbox{if } i\neq k. \\
\end{array}\right.
\]
Hence, the left cluster variable $A_c=d$.
By definition,  \eqref{snvodqv} is in the image of $\psi$ when $\B_0$ and  $\B^0$ are in general position, or equivalently when $d\neq 0$. In other words, the image of $\psi$ is precisely the non-vanishing locus of the cluster variable $A_c$. In cluster theory, localization of a cluster $\mathrm{K}_2$ variety at a cluster variable $A_c$ is again a cluster $\mathrm{K}_2$ variety, which can be obtained by freezing the vertex $c$. Therefore the image of $\psi$ is also a cluster $\mathrm{K}_2$ variety.
\end{proof}

Now we make $\conf^\gamma_\beta(\mathcal{C})\times\mathbb{G}_m$ into a cluster $\mathrm{K}_2$ variety by adding an extra frozen variable $d$ corresponding to the $\mathbb{G}_m$ factor. There should not be no arrows connecting $c$ and the unfrozen variables of  $\conf^\gamma_\beta(\mathcal{C})$ because the extra $\mathbb{G}_m$ factor will not affect their mutations. However, there is freedom of adding arrows connecting $c$ and the frozen variables of $\conf^\gamma_\beta(\mathcal{C})$. The next proposition shows that these arrows can be uniquely determined by requiring $\psi$ to be a quasi-cluster isomorphism onto its image.

\begin{prop}\label{quasi iso} The space $\conf_\beta^\gamma(\mathcal{C})\times \mathbb{G}_m$ can be equipped with a unique cluster $\mathrm{K}_2$ structure which extends the cluster $\mathrm{K}_2$ structure on $\conf_\beta^\gamma(\mathcal{C})$ by adding one extra frozen vertex $c$ and possibly arrows between $c$ and the original frozen part, such that $\psi$ becomes a quasi-cluster isomorphism onto its image.
\end{prop}
\begin{proof} Suppose we start with a standard representative in the image of $\psi$ as follows.
\[
\begin{tikzpicture}
\node (u0) at (0,2) [] {$\U_-$};
\node (u1) at (2,2) [] {$\U_-y_1$};
\node (u2) at (4,2) [] {$\cdots$};
\node (u3) at (6,2) [] {$\U_-y_m$};
\node (d-1) at (-3,0) [] {$\U_+$};
\node (d0) at (0,0) [] {$x_0\U_+$};
\node (d1) at (2,0) [] {$x_1\U_+$};
\node (d2) at (4,0) [] {$\cdots$};
\node (d3) at (6,0) [] {$x_l\U_+$};
\draw [->] (u0) -- node [above] {$s_{j_1}$} (u1);
\draw [->] (u1) -- node [above] {$s_{j_2}$} (u2);
\draw [->] (u2) -- node [above] {$s_{j_m}$} (u3);
\draw [->] (d0) -- node [below] {$s_{i_1}$} (d1);
\draw [->] (d1) -- node [below] {$s_{i_2}$} (d2);
\draw [->] (d2) -- node [below] {$s_{i_l}$} (d3);
\draw (u0) -- (d0);
\draw (u3) -- (d3);
\draw (u0) -- (d-1);
\draw [->] (d-1) -- node [below] {$s_i$} (d0);
\draw [dashed] (-5,-0.5) -- (7,-0.5) -- (7,0.5) -- (-1.5,0.5) -- (0,1.5) -- (7,1.5) -- (7,2.5) -- (-0.5,2.5) -- cycle;
\draw (-0.5,1) node [left] {$i$th} -- node [above right] {$c$} (0.5,1) ;
\end{tikzpicture}
\]
From the last proposition we know that $x_0\U_+=R_i(d)\U_+$ for some non-zero $d$.

To obtain the preimage of this representative under $\psi$, we need to delete the flag $\U_+$ at the lower left corner and re-scale the decorations along the bottom chain as follows.
\[
\begin{tikzpicture}
\node (u0) at (0,1.5) [] {$\U_-$};
\node (u1) at (2,1.5) [] {$\U_-y_1$};
\node (u2) at (4,1.5) [] {$\cdots$};
\node (u3) at (6,1.5) [] {$\U_-y_m$};
\node (d0) at (0,0) [] {$x_0h_0\U_+$};
\node (d1) at (2,0) [] {$x_1h_1\U_+$};
\node (d2) at (4,0) [] {$\cdots$};
\node (d3) at (6,0) [] {$x_lh_l\U_+$};
\draw [->] (u0) -- node [above] {$s_{j_1}$} (u1);
\draw [->] (u1) -- node [above] {$s_{j_2}$} (u2);
\draw [->] (u2) -- node [above] {$s_{j_m}$} (u3);
\draw [->] (d0) -- node [below] {$s_{i_1}$} (d1);
\draw [->] (d1) -- node [below] {$s_{i_2}$} (d2);
\draw [->] (d2) -- node [below] {$s_{i_l}$} (d3);
\draw (u0) -- (d0);
\draw (u3) -- (d3);
\draw [dashed] (-0.8,-0.5) -- (6.75,-0.5) -- (6.75,0.5) -- (0.5,0.5) -- (0.5,1) -- (6.75,1) -- (6.75,2) -- (-0.8,2) -- cycle;
\end{tikzpicture}
\]
Here $h_k\in \T$ are such that $(\U_-,~ x_0h_0\U_+)$ is a pinning and $(x_{k-1} h_{k-1} \U_+, ~ x_{k}h_k \U_+)$ are compatible.  

Set $\lambda_0^\vee = -\alpha_i^\vee$. Define co-characters $\lambda^\vee_k$ of $\T$ for $1\leq k\leq l$ by the recursive relation
\[
\lambda^\vee_k:=s_{i_k}\left(\lambda^\vee_{k-1}\right).
\]
 Note that $x_0=R_i(d)$. An easy calculation shows that $x_0h_0\in \U_-\U_+$ if and only if $h_0= d^{\lambda_0^\vee}$.  Since $(x_{k-1}\U_+,~ x_k\U_+)$ is a compatible pair, by definition we get $x_{k-1}^{-1}x_k\in \U_+\overline{s}_{i_k}\U_+$. Therefore,
\[
\U_+\left(x_{k-1}h_{k-1}\right)^{-1}\cdot x_kh_k\U_+ = \U_+\overline{s}_{i_k}\cdot s_{i_k}(h_{k-1}^{-1})h_k\U_+.
\]
The pair $(x_{k-1}h_{k-1}\U_+, ~x_kh_k\U_+)$ is compatible if and only if $h_k=s_{i_k}(h_{k-1})$. By induction we get $h_k= d^{\lambda_k^\vee}$ for $0\leq k\leq l$.

Next we investigate the pull-back of cluster $\mathrm{K}_2$ coordinates of $\conf^\gamma_{s_i\beta}(\mathcal{C})$ under $\psi$. Fix a word $\vec{t}$ for $(\beta,\gamma)$ and consider the word $(i,\vec{t})$ for $\left(s_i\beta,\gamma\right)$. Let $Q_{\vec{t}}$ be the quiver associated to $\vec{t}$, and $Q_{i, \vec{t}}$  the quiver associated to $(i,\vec{t})$ with the leftmost vertex $c$ frozen.

Recall that
\[
\psi^*\left(A_c\right)=\Delta_i(x_0)=\Delta_i(R_i(d))=d.
\]
We define $d$ to be the cluster variable $A'_c$ for the new frozen vertex $c$.

For any other string (vertex) $a$ associated to $(i,\vec{t})$ as the left picture below, there is a corresponding string (vertex) $a$ associated to $\vec{t}$ as the right right below.
\[
\begin{tikzpicture}[scale=.6]
\draw (-1,0) node [left] {$h$th level} -- (1,0) node [right] {$a$};
\draw (0,-1) node [below] {$x_k\U_+$} -- (0,1) node [above] {$\U_-y_j$};
\end{tikzpicture}\quad \quad\quad \quad 
\begin{tikzpicture}[scale=.6]
\draw (-1,0) node [left] {$h$th level} -- (1,0) node [right] {$a$};
\draw (0,-1) node [below] {$x_kd^{\lambda^\vee_k}\U_+$} -- (0,1) node [above] {$\U_-y_j$};
\end{tikzpicture}
\]
Let $\delta_a:=-\inprod{\lambda^\vee_k}{\omega_h}\in \mathbb{Z}$; then
\[
\psi^*\left(A_a\right)=\Delta_{h}\left(y_jx_k\right)=\Delta_{h}\left(y_jx_kh_k\right)d^{-\inprod{\lambda_k^\vee}{\omega_h}} =A'_aA'^{\delta_a}_c.
\]
In addition we define $\delta_c:=-\langle \lambda_0^\vee, \omega_i\rangle =1$. Let $I$ denote the vertices of $Q_{i, \vec{t}}$ and let $\varepsilon_{ij}$ be the exchange matrix encoded by $Q_{i, \vec{t}}$. 
The set $I'=I-\{c\}$ consists of vertices of $Q_{\vec{t}}$ and $I^\uf$ consists of unfrozen vertices of $Q_{\vec{t}}$. We claim that for any $a\in I^\uf$, we have 
\begin{equation}
\label{andnjonadjo}
\sum_{b\in I}\varepsilon_{ab}\delta_b=0.
\end{equation}
To see this, recall that there is projection map 
\[p:~~\conf^\gamma_\beta(\mathcal{C})\times\mathbb{G}_m\longrightarrow \conf^\gamma_\beta(\mathcal{C})\stackrel{\pi}{\longrightarrow} \conf^\gamma_\beta(\mathcal{B})\] 
As in  \cite[\S 3]{SWflag}, $\conf^\gamma_\beta(\mathcal{B})$ is  equipped with the cluster Poission variables $\{X_a'\}_{a\in I^\uf}$ such that 
\begin{equation}
\label{acsdonoj}
p^*\left(X'_a\right)=\prod_{b\in I'} A'^{\varepsilon_{ab}}_b.
\end{equation}
Consider  the composition
\[
p':= \conf^\gamma_{s_i\beta}(\mathcal{C}) \stackrel{\pi}{\longrightarrow} \conf^\gamma_{s_i\beta}(\mathcal{B}) \longrightarrow \conf^\gamma_{\beta}(\mathcal{B})
\]
Here the second map is rational, obtained by forgetting the flag $\B^{-1}$. Note that $\B^{-1}$ only changes the decorations on the other flags. Therefore we have $p=p'\circ \psi.$ Therefore for $a\in I^\uf$ we have 
\[
p^*\left(X'_a\right)=\psi^*\circ p'^*\left(X'_a\right)=\psi^*\left(\prod_{b\in I}A^{\varepsilon_{ab}}_b\right)=A'^{\varepsilon_{ac}}_c\prod_{b\in I'}A'^{\varepsilon_{ab}\delta_b}_cA'^{\varepsilon_{ab}}_b= A_c'^{\sum_{b\in I}\varepsilon_{ab}\delta_b}\cdot \prod_{b\in I'} A'^{\varepsilon_{ab}}_b.
\]
Comparing it with \eqref{acsdonoj}, we arrive at the identity \eqref{andnjonadjo}. 

Note that identity \eqref{andnjonadjo} satisfies the assumptions stated in Proposition \ref{changing frozen vertex}. Therefore we know that there is a unique way to extend the quiver of $\conf^\gamma_\beta(\mathcal{C})$ so that $\psi$ becomes a quasi-cluster isomorphism onto its image.
\end{proof}

\begin{exmp}
We continue from Example \ref{3.10}. Consider the map $\phi_1:\conf_\beta^\gamma(\mathcal{C})\times\mathbb{G}_m\rightarrow \conf_{s_1\beta}^\gamma(\mathcal{C})$. Let $d=A'_c$ be the coordinate for the $\mathbb{G}_m$ factor. Then in the preimage,
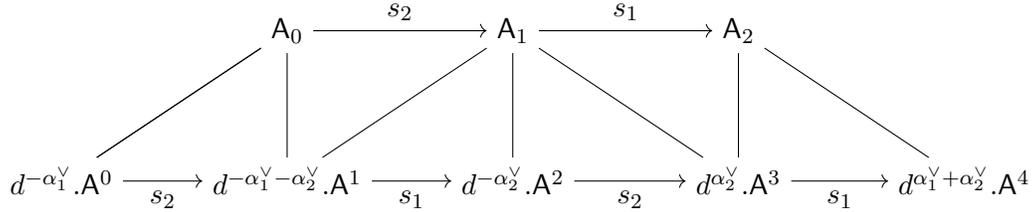
\[
\begin{tikzpicture}
\foreach \i in {0,1,2}
    {
    \node (u\i) at (3+\i*3,2) [] {$\A_\i$};
    }
\node (d0) at (0,0) [] {$d^{-\alpha_1^\vee}.\A^0$};
\node (d1) at (3,0) [] {$d^{-\alpha_1^\vee-\alpha_2^\vee}.\A^1$};
\node (d2) at (6,0) [] {$d^{-\alpha_2^\vee}.\A^2$};
\node (d3) at (9,0) [] {$d^{\alpha_2^\vee}.\A^3$};
\node (d4) at (12,0) [] {$d^{\alpha_1^\vee+\alpha_2^\vee}.\A^4$};
\draw [->] (u0) -- node [above] {$s_2$} (u1);
\draw [->] (u1) -- node [above] {$s_1$} (u2);
\draw [->] (d0) -- node [below] {$s_2$} (d1);
\draw [->] (d1) -- node [below] {$s_1$} (d2);
\draw [->] (d2) -- node [below] {$s_2$} (d3);
\draw [->] (d3) -- node [below] {$s_1$} (d4);
\draw (u0) -- (d0);
\draw (d0) -- (u0) -- (d1) --(u1)--(d2);
\draw (u1) -- (d3) -- (u2) -- (d4);
\end{tikzpicture}
\]
The change of decorations gives rise the the pull-backs $\phi_1^*\left(A_c\right)=A'_c$ and $\phi_1^*\left(A_a\right)=A'_ad^{\delta_a}=A'_aA'^{\delta_a}_c$ for $a\neq c$. The integers $\delta_a$ assigned to the vertices $a$ are as follows.
\[
\begin{tikzpicture}
\node [lightgray] (u0) at (0,2) [] {$1$};
\node [lightgray] (u1) at (3.5,2) [] {$1$};
\node [lightgray](u2) at (5.5,2) [] {$-1$};
\node [lightgray](u3) at (6.5,2) [] {$1$};
\node [lightgray](d1) at (1,0.5) [] {$2$};
\node [lightgray](d2) at (2.5,0.5) [] {$-2$};
\node [lightgray](d3) at (4.5,0.5) [] {$2$};
\draw [lightgray](-1,0.5) node [left] {$2$nd} -- (d1) -- (d2) -- (d3) -- (9,0.5);
\draw [lightgray](-1,2) node [left] {$1$st} -- (u0) -- (u1) -- (u2) -- (u3) -- (9,2);
\node (0u) at (1.75,2) [] {$\boxed{1}$};
\node (1u) at (4.25,2) [] {$0$};
\node (2u) at (6,2) [] {$0$};
\node (3u) at (7.5,2) [] {$\boxed{-1}$};
\node (1d) at (1.75,0.5) [] {$1$};
\node (2d) at (3.5,0.5) [] {$1$};
\node (3d) at (7.5,0.5) [] {$\boxed{-1}$};
\draw [->] (0u) -- (1u);
\draw [->] (1d) -- (0u);
\draw [->] (2d) -- (1d);
\draw [->] (2d) -- (3d);
\draw [->] (2u) -- (1u);
\draw [->] (2u) -- (3u);
\draw [->] (1u) -- (2d);
\draw [->] (3d) -- (2u);
\draw [dashed, ->] (3u) -- (3d);
\end{tikzpicture}
\]
Using these  $\delta_a$ we conclude that the cluster structure  on $\conf^\gamma_\beta(\mathcal{C})\times \mathbb{G}_m$ is given by the following quiver, where the right most vertex is the extra frozen vertex $c$.
\[
\begin{tikzpicture}
\node [lightgray] (u1) at (3.5,2) [] {$1$};
\node [lightgray](u2) at (5.5,2) [] {$-1$};
\node [lightgray](u3) at (6.5,2) [] {$1$};
\node [lightgray](d1) at (1,0.5) [] {$2$};
\node [lightgray](d2) at (2.5,0.5) [] {$-2$};
\node [lightgray](d3) at (4.5,0.5) [] {$2$};
\draw [lightgray](-1,0.5) node [left] {$2$nd} -- (d1) -- (d2) -- (d3) -- (9,0.5);
\draw [lightgray](-1,2) node [left] {$1$st} -- (u1) -- (u2) -- (u3) -- (9,2);
\node (1u) at (4.25,2) [] {$\bullet$};
\node (2u) at (6,2) [] {$\bullet$};
\node (3u) at (7.5,2) [] {$\square$};
\node (1d) at (1.75,0.5) [] {$\bullet$};
\node (2d) at (3.5,0.5) [] {$\bullet$};
\node (3d) at (7.5,0.5) [] {$\square$};
\node (c) at (9,1.25) [] {$\square$};
\draw [->] (2d) -- (1d);
\draw [->] (2d) -- (3d);
\draw [->] (2u) -- (1u);
\draw [->] (2u) -- (3u);
\draw [->] (1u) -- (2d);
\draw [->] (3d) -- (2u);
\draw [dashed, ->] (3u) -- (3d);
\draw [dashed, ->] (c) -- (3d);
\draw [dashed,->] (c) -- (3u);
\end{tikzpicture}
\]
\end{exmp}

\bigskip

\bibliographystyle{alphaurl-a}

\bibliography{biblio}

\begin{thebibliography}{LLMSS20}

\bibitem[Arn76]{Arn}
V.~I. Arnold.
\newblock Local normal forms of functions.
\newblock {\em Invent. Math.}, 35:87--109, 1976.
\newblock \href {https://doi.org/10.1007/BF01390134}
  {\path{doi:10.1007/BF01390134}}.

\bibitem[Baa13]{Baa}
Sebastian Baader.
\newblock Positive braids of maximal signature.
\newblock {\em Enseign. Math.}, 59(3-4):351--358, 2013.
\newblock \href {https://doi.org/10.4171/LEM/59-3-8}
  {\path{doi:10.4171/LEM/59-3-8}}.

\bibitem[BFZ05]{BFZ}
Arkady Berenstein, Sergey Fomin, and Andrei Zelevinsky.
\newblock Cluster algebras. {III}. {U}pper bounds and double {B}ruhat cells.
\newblock {\em Duke Math. J.}, 126(1):1--52, 2005.
\newblock \href {http://arxiv.org/abs/math/0305434}
  {\path{arXiv:math/0305434}}, \href
  {https://doi.org/10.1215/S0012-7094-04-12611-9}
  {\path{doi:10.1215/S0012-7094-04-12611-9}}.

\bibitem[Cas20]{casals2020lagrangian}
Roger Casals.
\newblock Lagrangian skeleta and plane curve singularities.
\newblock Preprint, 2020.
\newblock \href {http://arxiv.org/abs/2009.06737} {\path{arXiv:2009.06737}}.

\bibitem[CG20]{CasalsGao}
Roger Casals and Honghao Gao.
\newblock Infinitely many {L}agrangian fillings.
\newblock {\em To appear in Ann. Math (2)}, 2020.
\newblock Preprint.
\newblock \href {http://arxiv.org/abs/2001.01334} {\path{arXiv:2001.01334}}.

\bibitem[Cha10]{chantraine2010}
Baptiste Chantraine.
\newblock Lagrangian concordance of {L}egendrian knots.
\newblock {\em Algebr. Geom. Topol.}, 10(1):63--85, 2010.
\newblock \href {http://arxiv.org/abs/math/0611848}
  {\path{arXiv:math/0611848}}, \href {https://doi.org/10.2140/agt.2010.10.63}
  {\path{doi:10.2140/agt.2010.10.63}}.

\bibitem[Cha15]{chantraine2015}
Baptiste Chantraine.
\newblock Lagrangian concordance is not a symmetric relation.
\newblock {\em Quantum Topol.}, 6(3):451--474, 2015.
\newblock \href {https://doi.org/10.4171/QT/68} {\path{doi:10.4171/QT/68}}.

\bibitem[Che02]{Che02}
Yuri Chekanov.
\newblock Differential algebra of {L}egendrian links.
\newblock {\em Invent. Math.}, 150(3):441--483, 2002.
\newblock \href {http://arxiv.org/abs/math/9709233}
  {\path{arXiv:math/9709233}}, \href {https://doi.org/10.1007/s002220200212}
  {\path{doi:10.1007/s002220200212}}.

\bibitem[CN21]{CasalsNg}
Roger Casals and Lenhard Ng.
\newblock Braid loops with infinite monodromy on the legendrian contact dga.
\newblock Preprint, 2021.
\newblock \href {http://arxiv.org/abs/2101.02318} {\path{arXiv:2101.02318}}.

\bibitem[CZ20]{CZ}
Roger Casals and Eric Zaslow.
\newblock Legendrian weaves: N-graph calculus, flag moduli and applications,
  2020.
\newblock \href {http://arxiv.org/abs/2007.04943} {\path{arXiv:2007.04943}}.

\bibitem[EENS13]{EENS}
Tobias Ekholm, John Etnyre, Lenhard Ng, and Michael Sullivan.
\newblock Knot contact homology.
\newblock {\em Geom. Topol.}, 17(2):975--1112, 2013.
\newblock \href {http://arxiv.org/abs/1109.1542} {\path{arXiv:1109.1542}},
  \href {https://doi.org/10.2140/gt.2013.17.975}
  {\path{doi:10.2140/gt.2013.17.975}}.

\bibitem[EGH00]{EGH}
Yakov Eliashberg, Alexander Givental, and Helmut Hofer.
\newblock Introduction to symplectic field theory.
\newblock In {\em Visions in Mathematics}, Special Volume, Part II, pages
  560--673. Birkh{\"{a}}user Basel, 2000.
\newblock GAFA 2000 (Tel Aviv, 1999).
\newblock \href {http://arxiv.org/abs/math/0010059}
  {\path{arXiv:math/0010059}}, \href
  {https://doi.org/10.1007/978-3-0346-0425-3_4}
  {\path{doi:10.1007/978-3-0346-0425-3_4}}.

\bibitem[EHK16]{EHK12}
Tobias Ekholm, Ko~Honda, and Tam\'{a}s K\'{a}lm\'{a}n.
\newblock Legendrian knots and exact {L}agrangian cobordisms.
\newblock {\em J. Eur. Math. Soc. (JEMS)}, 18(11):2627--2689, 2016.
\newblock \href {http://arxiv.org/abs/1212.1519} {\path{arXiv:1212.1519}},
  \href {https://doi.org/10.4171/JEMS/650} {\path{doi:10.4171/JEMS/650}}.

\bibitem[Ekh07]{ekholm2007}
Tobias Ekholm.
\newblock Morse flow trees and {L}egendrian contact homology in 1-jet spaces.
\newblock {\em Geom. Topol.}, 11:1083--1224, 2007.
\newblock \href {http://arxiv.org/abs/math/0509386}
  {\path{arXiv:math/0509386}}, \href {https://doi.org/10.2140/gt.2007.11.1083}
  {\path{doi:10.2140/gt.2007.11.1083}}.

\bibitem[EL17]{EL}
Tobias Ekholm and Yanki Lekili.
\newblock Duality between {L}agrangian and {L}egendrian invariants, 2017.
\newblock \href {http://arxiv.org/abs/1701.01284} {\path{arXiv:1701.01284}}.

\bibitem[EN18]{ENsurvey}
John Etnyre and Lenhard Ng.
\newblock Legendrian contact homology in $\mathbb{R}^3$, 2018.
\newblock \href {http://arxiv.org/abs/1811.10966} {\path{arXiv:1811.10966}}.

\bibitem[ENS02]{ENS02}
John Etnyre, Lenhard Ng, and Joshua Sabloff.
\newblock Invariants of {L}egendrian knots and coherent orientations.
\newblock {\em J. Symplectic Geom.}, 1(2):321--367, 2002.
\newblock URL: \url{http://projecteuclid.org/euclid.jsg/1092316653}, \href
  {http://arxiv.org/abs/math/0101145} {\path{arXiv:math/0101145}}.

\bibitem[EP96]{EP}
Yakov. Eliashberg and Leonid Polterovich.
\newblock Local {L}agrangian {$2$}-knots are trivial.
\newblock {\em Ann. of Math. (2)}, 144(1):61--76, 1996.
\newblock \href {https://doi.org/10.2307/2118583} {\path{doi:10.2307/2118583}}.

\bibitem[EV18]{EV}
John Etnyre and Vera V\'{e}rtesi.
\newblock Legendrian satellites.
\newblock {\em Int. Math. Res. Not. IMRN}, 2018(23):7241--7304, 2018.
\newblock \href {http://arxiv.org/abs/1608.05695} {\path{arXiv:1608.05695}},
  \href {https://doi.org/10.1093/imrn/rnx106} {\path{doi:10.1093/imrn/rnx106}}.

\bibitem[FG06]{FGteich}
Vladimir Fock and Alexander Goncharov.
\newblock Moduli spaces of local systems and higher {T}eichm\"{u}ller theory.
\newblock {\em Publ. Math. Inst. Hautes \'{E}tudes Sci.}, 103:1--211, 2006.
\newblock \href {http://arxiv.org/abs/math/0311149}
  {\path{arXiv:math/0311149}}, \href
  {https://doi.org/10.1007/s10240-006-0039-4}
  {\path{doi:10.1007/s10240-006-0039-4}}.

\bibitem[FG09]{FGensemble}
Vladimir Fock and Alexander Goncharov.
\newblock Cluster ensembles, quantization and the dilogarithm.
\newblock {\em Ann. Sci. \'{E}c. Norm. Sup\'{e}r. (4)}, 42(6):865--930, 2009.
\newblock \href {http://arxiv.org/abs/math/0311245}
  {\path{arXiv:math/0311245}}, \href {https://doi.org/10.24033/asens.2112}
  {\path{doi:10.24033/asens.2112}}.

\bibitem[FR11]{FR}
Dmitry Fuchs and Dan Rutherford.
\newblock Generating families and {L}egendrian contact homology in the standard
  contact space.
\newblock {\em J. Topol.}, 4(1):190--226, 2011.
\newblock \href {http://arxiv.org/abs/0807.4277} {\path{arXiv:0807.4277}},
  \href {https://doi.org/10.1112/jtopol/jtq033}
  {\path{doi:10.1112/jtopol/jtq033}}.

\bibitem[Fra16]{Fra}
Chris Fraser.
\newblock Quasi-homomorphisms of cluster algebras.
\newblock {\em Adv. in Appl. Math.}, 81:40--77, 2016.
\newblock \href {http://arxiv.org/abs/1509.05385} {\path{arXiv:1509.05385}},
  \href {https://doi.org/10.1016/j.aam.2016.06.005}
  {\path{doi:10.1016/j.aam.2016.06.005}}.

\bibitem[Fuc03]{F}
Dmitry Fuchs.
\newblock Chekanov-{E}liashberg invariant of {L}egendrian knots: existence of
  augmentations.
\newblock {\em J. Geom. Phys.}, 47(1):43--65, 2003.
\newblock \href {https://doi.org/10.1016/S0393-0440(01)00013-4}
  {\path{doi:10.1016/S0393-0440(01)00013-4}}.

\bibitem[FZ02]{FZI}
Sergey Fomin and Andrei Zelevinsky.
\newblock Cluster algebras. {I}. {F}oundations.
\newblock {\em J. Amer. Math. Soc.}, 15(2):497--529, 2002.
\newblock \href {http://arxiv.org/abs/math/0104151}
  {\path{arXiv:math/0104151}}, \href
  {https://doi.org/10.1090/S0894-0347-01-00385-X}
  {\path{doi:10.1090/S0894-0347-01-00385-X}}.

\bibitem[FZ03]{FZII}
Sergey Fomin and Andrei Zelevinsky.
\newblock Cluster algebras. {II}. {F}inite type classification.
\newblock {\em Invent. Math.}, 154(1):63--121, 2003.
\newblock \href {http://arxiv.org/abs/math/0208229}
  {\path{arXiv:math/0208229}}, \href
  {https://doi.org/10.1007/s00222-003-0302-y}
  {\path{doi:10.1007/s00222-003-0302-y}}.

\bibitem[FZ07]{FZIV}
Sergey Fomin and Andrei Zelevinsky.
\newblock Cluster algebras. {IV}. {C}oefficients.
\newblock {\em Compos. Math.}, 143(1):112--164, 2007.
\newblock \href {http://arxiv.org/abs/math/0602259}
  {\path{arXiv:math/0602259}}, \href
  {https://doi.org/10.1112/S0010437X06002521}
  {\path{doi:10.1112/S0010437X06002521}}.

\bibitem[GHKK18]{GHKK}
Mark Gross, Paul Hacking, Sean Keel, and Maxim Kontsevich.
\newblock Canonical bases for cluster algebras.
\newblock {\em J. Amer. Math. Soc.}, 31(2):497--608, 2018.
\newblock \href {http://arxiv.org/abs/1411.1394} {\path{arXiv:1411.1394}},
  \href {https://doi.org/10.1090/jams/890} {\path{doi:10.1090/jams/890}}.

\bibitem[GLS11]{GLS}
C.~Geiss, B.~Leclerc, and J.~{Schr\"{o}er}.
\newblock Kac-{M}oody groups and cluster algebras.
\newblock {\em Adv. Math.}, 228(1):329--433, 2011.

\bibitem[GPS18]{GPSsectorial}
Sheel Ganatra, John Pardon, and Vivek Shende.
\newblock Sectorial descent for wrapped {F}ukaya categories, 2018.
\newblock \href {http://arxiv.org/abs/1809.03427} {\path{arXiv:1809.03427}}.

\bibitem[GR91]{GR}
I.~M. Gelfand and V.~S. Retakh.
\newblock Determinants of matrices over noncommutative rings.
\newblock {\em Funktsional. Anal. i Prilozhen.}, 25(2):13--25, 96, 1991.
\newblock \href {https://doi.org/10.1007/BF01079588}
  {\path{doi:10.1007/BF01079588}}.

\bibitem[GS18]{GS2}
Alexander Goncharov and Linhui Shen.
\newblock Donaldson-{T}homas transformations of moduli spaces of {G}-local
  systems.
\newblock {\em Adv. Math.}, 327:225--348, 2018.
\newblock \href {http://arxiv.org/abs/1602.06479} {\path{arXiv:1602.06479}},
  \href {https://doi.org/10.1016/j.aim.2017.06.017}
  {\path{doi:10.1016/j.aim.2017.06.017}}.

\bibitem[GS19]{GS3}
Alexander Goncharov and Linhui Shen.
\newblock Quantum geometry of moduli spaces of local systems and representation
  theory.
\newblock Preprint, 2019.
\newblock \href {http://arxiv.org/abs/1904.10491} {\path{arXiv:1904.10491}}.

\bibitem[{K}{\'{a}}l05]{Kalman}
Tam\'{a}s {K}{\'{a}}lm\'{a}n.
\newblock Contact homology and one parameter families of {L}egendrian knots.
\newblock {\em Geom. Topol.}, 9:2013--2078, 2005.
\newblock \href {http://arxiv.org/abs/math/0407347}
  {\path{arXiv:math/0407347}}, \href {https://doi.org/10.2140/gt.2005.9.2013}
  {\path{doi:10.2140/gt.2005.9.2013}}.

\bibitem[K{\'{a}}l06]{Kalman06}
Tam\'{a}s K{\'{a}}lm{\'{a}}n.
\newblock Braid-positive {L}egendrian links.
\newblock {\em Int. Math. Res. Not.}, pages Art ID 14874, 29, 2006.
\newblock \href {http://arxiv.org/abs/math/0608457}
  {\path{arXiv:math/0608457}}, \href {https://doi.org/10.1155/IMRN/2006/14874}
  {\path{doi:10.1155/IMRN/2006/14874}}.

\bibitem[Kar20]{karlsson2017note}
Cecilia Karlsson.
\newblock A note on coherent orientations for exact {L}agrangian cobordisms.
\newblock {\em Quantum Topol.}, 11(1):1--54, 2020.
\newblock \href {http://arxiv.org/abs/1707.04219} {\path{arXiv:1707.04219}},
  \href {https://doi.org/10.4171/QT/132} {\path{doi:10.4171/QT/132}}.

\bibitem[Kel13]{Kelperiod}
Bernhard Keller.
\newblock The periodicity conjecture for pairs of {D}ynkin diagrams.
\newblock {\em Ann. of Math. (2)}, 177(1):111--170, 2013.
\newblock \href {http://arxiv.org/abs/1001.1531} {\path{arXiv:1001.1531}},
  \href {https://doi.org/10.4007/annals.2013.177.1.3}
  {\path{doi:10.4007/annals.2013.177.1.3}}.

\bibitem[Kel17]{KelDT}
Bernhard Keller.
\newblock Quiver mutation and combinatorial {DT}-invariants.
\newblock {\em Discrete Mathematics and Theoretical Computer Science}, 2017.
\newblock \href {http://arxiv.org/abs/1709.03143} {\path{arXiv:1709.03143}}.

\bibitem[Kum02]{Kum}
Shrawan Kumar.
\newblock {\em Kac-{M}oody groups, their flag varieties and representation
  theory}, volume 204 of {\em Progress in Mathematics}.
\newblock Birkh\"{a}user Boston, Inc., Boston, MA, 2002.
\newblock \href {https://doi.org/10.1007/978-1-4612-0105-2}
  {\path{doi:10.1007/978-1-4612-0105-2}}.

\bibitem[LLMSS20]{lee2018frieze}
Kyungyong Lee, Li~Li, Matthew Mills, Ralf Schiffler, and Alexandra Seceleanu.
\newblock Frieze varieties: a characterization of the finite-tame-wild
  trichotomy for acyclic quivers.
\newblock {\em Adv. Math.}, 367:107130, 33, 2020.
\newblock \href {http://arxiv.org/abs/1803.08459} {\path{arXiv:1803.08459}},
  \href {https://doi.org/10.1016/j.aim.2020.107130}
  {\path{doi:10.1016/j.aim.2020.107130}}.

\bibitem[LS15]{LS}
Kyungyong Lee and Ralf Schiffler.
\newblock Positivity for cluster algebras.
\newblock {\em Ann. of Math. (2)}, 182(1):73--125, 2015.
\newblock \href {http://arxiv.org/abs/1306.2415} {\path{arXiv:1306.2415}},
  \href {https://doi.org/10.4007/annals.2015.182.1.2}
  {\path{doi:10.4007/annals.2015.182.1.2}}.

\bibitem[Nad09]{Nadler}
David Nadler.
\newblock Microlocal branes are constructible sheaves.
\newblock {\em Selecta Math. (N.S.)}, 15(4):563--619, 2009.
\newblock \href {http://arxiv.org/abs/math/0612399}
  {\path{arXiv:math/0612399}}, \href
  {https://doi.org/10.1007/s00029-009-0008-0}
  {\path{doi:10.1007/s00029-009-0008-0}}.

\bibitem[Ng03]{NG2003}
Lenhard Ng.
\newblock Computable {L}egendrian invariants.
\newblock {\em Topology}, 42(1):55--82, 2003.
\newblock \href {http://arxiv.org/abs/math/0011265}
  {\path{arXiv:math/0011265}}, \href
  {https://doi.org/10.1016/S0040-9383(02)00010-1}
  {\path{doi:10.1016/S0040-9383(02)00010-1}}.

\bibitem[NRSSZ15]{NRSSZ}
Lenhard Ng, Dan Rutherford, Vivek Shende, Steven Sivek, and Eric Zaslow.
\newblock Augmentations are sheaves.
\newblock Preprint, 2015.
\newblock \href {http://arxiv.org/abs/1502.04939} {\path{arXiv:1502.04939}}.

\bibitem[NZ12]{NZ}
Tomoki Nakanishi and Andrei Zelevinsky.
\newblock On tropical dualities in cluster algebras.
\newblock In {\em Algebraic groups and quantum groups}, volume 565 of {\em
  Contemp. Math.}, pages 217--226. Amer. Math. Soc., Providence, RI, 2012.
\newblock \href {http://arxiv.org/abs/1101.3736} {\path{arXiv:1101.3736}},
  \href {https://doi.org/10.1090/conm/565/11159}
  {\path{doi:10.1090/conm/565/11159}}.

\bibitem[Pan17]{Pan_2017}
Yu~Pan.
\newblock Exact {L}agrangian fillings of {L}egendrian {$(2,n)$} torus links.
\newblock {\em Pacific J. Math.}, 289(2):417--441, 2017.
\newblock \href {http://arxiv.org/abs/1607.03167} {\path{arXiv:1607.03167}},
  \href {https://doi.org/10.2140/pjm.2017.289.417}
  {\path{doi:10.2140/pjm.2017.289.417}}.

\bibitem[Sab05]{Sab}
Joshua Sabloff.
\newblock Augmentations and rulings of {L}egendrian knots.
\newblock {\em Int. Math. Res. Not.}, 2005(19):1157--1180, 2005.
\newblock \href {http://arxiv.org/abs/math/0409032}
  {\path{arXiv:math/0409032}}, \href {https://doi.org/10.1155/IMRN.2005.1157}
  {\path{doi:10.1155/IMRN.2005.1157}}.

\bibitem[She14]{Sh}
L.~Shen.
\newblock Stasheff polytopes and the coordinate ring of the cluster
  $\mathcal{X}$-variety of type ${A}_n$.
\newblock {\em Selecta Math. (N.S.)}, 20(3):929--959, 2014.
\newblock \href{https://arxiv.org/abs/1104.3528}{arXiv:1104.3528}.

\bibitem[Siv11]{Siv10}
Steven Sivek.
\newblock A bordered {C}hekanov-{E}liashberg algebra.
\newblock {\em J. Topol.}, 4(1):73--104, 2011.
\newblock \href {http://arxiv.org/abs/1004.4929} {\path{arXiv:1004.4929}},
  \href {https://doi.org/10.1112/jtopol/jtq035}
  {\path{doi:10.1112/jtopol/jtq035}}.

\bibitem[STWZ19]{STWZ}
Vivek Shende, David Treumann, Harold Williams, and Eric Zaslow.
\newblock Cluster varieties from {L}egendrian knots.
\newblock {\em Duke Math. J.}, 168(15):2801--2871, 2019.
\newblock \href {http://arxiv.org/abs/1512.08942} {\path{arXiv:1512.08942}},
  \href {https://doi.org/10.1215/00127094-2019-0027}
  {\path{doi:10.1215/00127094-2019-0027}}.

\bibitem[STZ17]{STZ}
Vivek Shende, David Treumann, and Eric Zaslow.
\newblock Legendrian knots and constructible sheaves.
\newblock {\em Invent. Math.}, 207(3):1031--1133, 2017.
\newblock \href {http://arxiv.org/abs/1402.0490} {\path{arXiv:1402.0490}},
  \href {https://doi.org/10.1007/s00222-016-0681-5}
  {\path{doi:10.1007/s00222-016-0681-5}}.

\bibitem[SW19]{SWflag}
Linhui Shen and Daping Weng.
\newblock Cluster structures on double {B}ott-{S}amelson cells.
\newblock To appear in {\it Forum of Mathematics, Sigma}, 2019.
\newblock \href {http://arxiv.org/abs/1904.07992} {\path{arXiv:1904.07992}}.

\bibitem[Syl19]{Sylvan}
Zachary Sylvan.
\newblock On partially wrapped {F}ukaya categories.
\newblock {\em J. Topol.}, 12(2):372--441, 2019.
\newblock \href {http://arxiv.org/abs/1604.02540} {\path{arXiv:1604.02540}},
  \href {https://doi.org/10.1112/topo.12088} {\path{doi:10.1112/topo.12088}}.

\bibitem[Wen16]{weng}
Daping Weng.
\newblock {D}onaldson-{T}homas transformation of {G}rassmannian.
\newblock Preprint, 2016.
\newblock \href {http://arxiv.org/abs/1603.00972} {\path{arXiv:1603.00972}}.

\end{thebibliography}

\end{document}